\documentclass[10pt]{article}
\usepackage{amsfonts,color}
\usepackage{amssymb}
\usepackage{footnote}
\usepackage{mathtools}
\usepackage{amsthm,wasysym}
\usepackage[english]{babel}
\usepackage{bbm}
\usepackage{colonequals}
\usepackage{setspace}
\usepackage{titlesec}

\usepackage{graphicx}
\usepackage{amsfonts}
\usepackage{hyperref}
\usepackage{subcaption}
\usepackage{amsmath}
\usepackage{nicefrac}
\usepackage{gensymb}
\usepackage{enumerate}
\usepackage[nameinlink,capitalize]{cleveref}
\usepackage{chngcntr}

\usepackage{cite}
\usepackage[nottoc,numbib]{tocbibind}

\usepackage{pgf,tikz,pgfplots}
\pgfplotsset{compat=1.14}
\usepackage{mathrsfs}
\usetikzlibrary{arrows}
\usetikzlibrary{arrows.meta}

\allowdisplaybreaks[1]

\makeindex

\newtheoremstyle{break}
{\topsep}{\topsep}%
{\itshape}{}%
{\bfseries}{}%
{\newline}{}%
\theoremstyle{break}

\counterwithin{theorem}{section}
\counterwithin{lemma}{section}
\counterwithin{remark}{section}
\counterwithin{observation}{section}
\counterwithin{definition}{section}

\counterwithin{figure}{section}
\counterwithin{equation}{section}

\makeatletter
\let\@fnsymbol\@arabic
\makeatother

\crefname{Problem}{Problem.}{Problem.}

\DeclareMathOperator{\at}{\bigg\vert}
\newcommand{\vb}[1]{\mathbf{#1}}
\newcommand{\bm}[1]{\boldsymbol{#1}}

\DeclareMathOperator*{\argmin}{arg\,min}

\DeclareMathOperator{\sym}{\mathrm{sym}}
\DeclareMathOperator{\spa}{\mathrm{span}}

\DeclareMathOperator{\dev}{\mathrm{dev}}
\DeclareMathOperator{\skw}{\mathrm{skew}}
\DeclareMathOperator{\Anti}{\mathrm{Anti}}
\DeclareMathOperator{\axl}{\mathrm{axl}}
\DeclareMathOperator{\tr}{\mathrm{tr}}
\DeclareMathOperator{\id}{\mathrm{id}}
\DeclareMathOperator{\cof}{\mathrm{cof}}
\newcommand{\jump}[1]{\ensuremath{[\![#1]\!]} }
\newcommand{\one}{\bm{\mathbbm{1}}}
\newcommand{\con}[2]{\langle {#1} , \, {#2} \rangle}
\newcommand{\norm}[1]{\| {#1} \|}

\newcommand{\der}[2]{ \dfrac{\mathrm{d}{#1}}{\mathrm{d}{#2}}}
\newcommand{\pder}[2]{ \dfrac{\partial{#1}}{\partial{#2}}}

\newcommand{\dd}{\mathrm{d}}
\newcommand{\D}{\mathrm{D}}
\newcommand{\Dcov}{\mathrm{D}^{\mathrm{cov}}_t}

\DeclareMathOperator{\di}{\mathrm{div}}
\DeclareMathOperator{\Di}{\mathrm{Div}}

\DeclareMathOperator{\curl}{\mathrm{curl}}
\DeclareMathOperator{\Curl}{\mathrm{Curl}}

\newcommand{\so}{\mathfrak{so}}

\newcommand{\Nedtwo}{\mathcal{N}_{II}}
\newcommand{\CG}{\mathcal{CG}}
\newcommand{\DG}{\mathcal{DG}}

\newcommand{\Q}{\mathit{Q}}
\newcommand{\Le}{{\mathit{L}^2}}

\newcommand{\Hone}{\mathit{H}^1}

\newcommand{\Hc}[1]{\mathit{H}(\mathrm{curl}{#1})}
\newcommand{\Hct}[1]{\mathit{H}(\mathrm{curl}_{t}{#1})}

\newcommand{\body}{\mathcal{B}}

\newcommand{\vol}{V}
\newcommand{\surf}{A}
\newcommand{\curv}{s}
\newcommand{\tansurf}{T_p}

\renewcommand{\O}{\mathcal{O}}
\newcommand{\R}{\mathbb{R}}

\newcommand{\C}{\mathit{C}}

\newcommand{\Cm}{\mathbb{C}}

\newcommand{\lama}{\lambda_{\mathrm{M}}}
\newcommand{\mua}{\mu_{\mathrm{M}}}
\newcommand{\lame}{\lambda_{\mathrm{e}}}
\newcommand{\mue}{\mu_{\mathrm{e}}}
\newcommand{\muc}{\mu_{\mathrm{c}}}
\newcommand{\Lc}{L_\mathrm{c}}
\newcommand{\Ce}{\mathbb{C}_{\mathrm{e}}}
\newcommand{\De}{\mathbb{D}_{\mathrm{e}}}

\newcommand{\Lm}{\mathbb{L}}

\newcommand{\wrt}{\text{w.r.t.}}

\newcommand{\ud}{\vb{u}}
\newcommand{\Rtheta}{\bm{\Theta}}
\newcommand{\rtheta}{\bm{\theta}}
\newcommand{\Pt}{\bm{P}}
\newcommand{\Pn}{\bm{Q}}

\newcommand{\mm}{\mathrm{mm}}
\newcommand{\MPa}{\mathrm{MPa}}
\newcommand{\Nwtn}{\mathrm{N}}

\newcommand{\AS}[1]{{\color{black} #1}}

\title{Intrinsic mixed-dimensional beam-shell-solid couplings in linear Cosserat continua via tangential differential calculus}

\author{\normalsize{Adam Sky}\thanks{Corresponding author: Adam Sky, Institute of Computational Engineering and Sciences, Department of Engineering, Faculty of Science, Technology and Medicine, University of Luxembourg, 6 Avenue de la Fonte, L-4362 Esch-sur-Alzette, Luxembourg, email: adam.sky@uni.lu}
    , \quad
	\normalsize{Jack S. Hale}\thanks{Jack S. Hale, Institute of Computational Engineering and Sciences, Department of Engineering, Faculty of Science, Technology and Medicine, University of Luxembourg, 6 Avenue de la Fonte, L-4362 Esch-sur-Alzette, Luxembourg, email: jack.hale@uni.lu}
    , \quad
    \normalsize{Andreas Zilian}\thanks{Andreas Zilian, Institute of Computational Engineering and Sciences, Department of Engineering, Faculty of Science, Technology and Medicine, University of Luxembourg, 6 Avenue de la Fonte, L-4362 Esch-sur-Alzette, Luxembourg, email: andreas.zilian@uni.lu}
    , \quad
	\normalsize{St\'ephane P. A. Bordas}\thanks{St\'ephane P. A. Bordas, Institute of Computational Engineering and Sciences, Department of Engineering, Faculty of Science, Technology and Medicine, University of Luxembourg, 6 Avenue de la Fonte, L-4362 Esch-sur-Alzette, Luxembourg, email: stephane.bordas@alum.northwestern.edu}
    \\ 
    \normalsize{and} \quad
	\normalsize{Patrizio Neff}\thanks{Patrizio Neff, Chair for Nonlinear 
		Analysis and Modelling, Faculty of Mathematics, Universit\"{a}t Duisburg-Essen,
		Thea-Leymann Str. 9, 45127 Essen, Germany, email: patrizio.neff@uni-due.de}
}

\begin{document}

\maketitle

\begin{abstract}
We present an approach to the coupling of mixed-dimensional continua by employing the mathematically enriched linear Cosserat micropolar model. The kinematical reduction of the model to lower dimensional domains leaves its fundamental degrees of freedom intact. Consequently, the degrees of freedom intrinsically agree even at the interface with a domain of a different dimensionality. Thus, this approach circumvents the need for intermediate finite elements or mortar methods. We introduce the derivations of all models of various dimensions using tangential differential calculus. The coupling itself is then achieved by defining a mixed-dimensional action functional with consistent Sobolev trace operators. Finally, we present numerical examples involving a three-dimensional silicone-rubber block reinforced with a curved graphite shell on its lower surface, a three-dimensional silver block reinforced with a graphite plate and beams, and lastly, intersecting silver shells reinforced with graphite beams. 
\\
\vspace*{0.25cm}
\\

\end{abstract}

\section{Introduction}
The design of structural parts in engineering practices often entails combinations of physically large and small components, for example in sandwich structures \cite{Bence,DICARA2024117856}, compound composites \cite{GORTHOFER2020109456}, or fibre-reinforced materials \cite{YANG2024105545,Alessio,MULLER201636,BUCK2015159}. From a modelling point of view, it is often possible to consider these joint components as a single continuous body, and thus model its behaviour through a single continuum theory. However, from a computational perspective, the latter approach may prove inefficient or even unfeasible for a given computational power, as it requires the discretisation of the domain to bridge the scale-gap between the small scale components and the large scale components. This task is difficult in terms of both, generating satisfactory finite element meshes for high-fidelity simulations \cite{Zou1,ZOU2024103653,Sevilla} (see also \cref{fig:conv}), as well as in the required computational effort that arises from the very fine mesh which is needed for scale-transition. Consequently, it is commonplace to regard designs containing small and large parts as mixed-dimensional continua. Let a $d$-dimensional body $\body \subset \R^d$ represent a large scale component, the smaller scale components can be modelled as continua on its $k$-codimensions $\Xi \subset \R^{d - k}$. This approach allows for coarser discretisations of both the large and small scale components, while compensating for inaccuracies in the behaviour of the small scale components with idealised models, which are fine-tuned for the expected behaviour in the small scale. However, since multiple continuum models are now used to model a single structural part, it becomes necessary to couple them in a single computational framework. The latter often proves challenging, as the differing continuum models may entail different types of degrees of freedom. A common example for the aforementioned difficulty and the focus of this work are kinematical couplings. Namely, in linear elasticity one is often compelled to couple the kinematical displacement field $\ud:\vol \subset \R^3 \to \R^3$ of a three-dimensional body $\vol \subset \R^3$ with kinematical fields of codimensional models, such as shells or beams. For linear elasticity, the shell models are the Cosserat \cite{ghiba2023essaydeformationmeasuresisotropic,Ghiba2021,Ghiba2023,Ghiba2020I,Ghiba2020II,Mircea,NEBEL2023116309}, and the classical Naghdi \cite{neunteufel_hellanherrmannjohnson_2019,NS21,SCHOLLHAMMER2019172,HuBordas2020} and Koiter shells \cite{Schollhammer2019,ghiba2023essaydeformationmeasuresisotropic,hale_simple_2018,MichaelThesis}, or their flat correspondents, the Reissner--Mindlin \cite{sky2023reissnermindlin,pechstein_tdnns_2017,NeffReissner} and Kirchhoff--Love plates \cite{Nguyen2021,ARF2023116198}. The beam models are the well-known Timoshenko–Ehrenfest \cite{Banerjee,Yuan,VO2022114883} and Euler--Bernoulli \cite{KaiserKirch,BORKOVIC2023115848,Meier2019} beams, and in the geometrically nonlinear case also the Cosserat beam \cite{Choi2022,Harsch,WEEGER2017100}, possibly with enhanced Saint-Venant torsion kinematics \cite{SIMO1991371,Gruttmann,Pi}. It is precisely in this scenario that the coupling becomes difficult, as the shell and beam models introduce rotation degrees of freedom $\rtheta:\Xi \subset \R^{3-k}\to \R^{3}$, which the three-dimensional model simply does not have. Thus, the coupling requires some additional treatment, e.g., using intermediate finite elements \cite{BOURNIVAL2010838,Klarmann2022}, static condensation \cite{Shim2002,Song}, Nitsche's method \cite{nguyena2013nitschesmethodmethodmixed,HANSBO2022114707}, or mortar approaches \cite{Steinbrecher2020,Steinbrecher2022,Yamamoto2019}. \AS{Additionally and related to the approach in this work, in \cite{BURMAN201851} the authors coupled the Euler--Bernoulli beam and Kirchhoff--Love plate using CutFEM, by exploiting that only one kinematic field is present in the mixed-dimensional problem.} 
Specifically in mechanical continuum models it is clear that the difficulty in the coupling stems from the lacking agreement in degrees of freedom. As such, this work considers an alternative approach with an enriched continuum model \cite{SKY2022115298,Neff2014,Eringen1999,Mindlin1964,Forest2006}, such that degrees of freedom intrinsically agree for all possible dimensions. Namely, this work employs the Cosserat micropolar theory as its starting point.    

The Cosserat micropolar continuum \cite{boon2024mixed,Munch2011,Neff2009,Neff2015,Shirani}, a concept introduced by the Cosserat brothers in 1909 \cite{Cosserat}, represents a generalised continuum model that extends beyond the classical Cauchy continuum theory. Unlike the Cauchy continuum, which considers material points characterised solely by their position, the Cosserat continuum also takes into account possible independent rotations of each material point, allowing for the introduction of more intricate kinematics. Effectively, the Cosserat theory turns each material point into a non-deformable solid microbody, implying the existence of a local micro-moment of inertia. Another consequence is the introduction of couple-forces $\bm{M}: \vol \subset \R^3 \to \so(3)$, relating to higher order tractions. Concisely, the Cosserat model naturally encompasses both a displacement field $\ud :\vol \subset \R^3 \to \R^3$ and a rotation field $\rtheta: \vol \subset \R^3 \to \R^3$, even in three-dimensional bodies $\vol \subset \R^3$. As a result, any dimensional reduction of the model to shells or beams leaves both the translational and rotational degrees of freedom intact. Therefore, the model is ideal for mixed-dimensional designs, as the degrees of freedom intrinsically agree even on interfaces of differing dimensionality $\jump{\ud}|_\Xi = \jump{\rtheta}|_\Xi = 0$. In fact, the coupling procedure is reduced to combining the bulk energy functional with lower dimensional energy functionals based on Sobolev trace operators \cite{Hiptmair} applied to the three-dimensional fields on codimensional domains $I_\vol(\ud,\rtheta) + I_\surf(\tr_\surf \ud,\tr_\surf \rtheta) + I_\curv(\tr_\curv \ud,\tr_\curv \rtheta)$.  

The derivation of reduced dimensional models is an involved procedure, especially if the lower dimensional domain is curved. This is the case since derivatives on curved domains are naturally expressed using differential geometry, implying the usage of co- and contravariant derivatives as well as Christoffel symbols. This makes the derivation process difficult to interpret, and the translation to a finite element software challenging. Alternatively, a recent approach called tangential differential calculus (TDC) \cite{Delfour,Gurtin1975,Hansbo2015} allows to circumvent the need for curvilinear coordinates and Christoffel symbols by introducing equivalent differentiation operators based on projections \cite{Schollhammer2019,SCHOLLHAMMER2019172,HANSBO20141,Hansbo2014}. To clarify, a gradient on a curved surface can be computed as the tangential projection of the gradient with respect to the Cartesian coordinates onto the surface $\nabla_t \lambda = \Pt \nabla \lambda$. If the surface is parameterised $\vb{r}:\omega \subset \R^2 \to \surf \subset \R^3$, then the tangential gradient $\nabla_t \lambda$ and the gradient given by the chain-rule $\nabla \lambda$ are equivalent. We give a short summary of the relation of tangential differential calculus to classical differential geometry via tensor calculus \cite{Itskov} in \cref{ap:curves,ap:shellshifter,ap:tdc,ap:tensorcal,ap:weingarten,ap:tangradshell}. While tangential differential calculus can also be used to introduce the geometry of the domain implicitly, e.g., with a level-set function \cite{FRIES2023116223,Kaiser}, we employ this framework for its advantageous applicability to automated solvers of partial differential equations. Specifically, the reinterpretation of co- and contravariant gradients as projected gradients of the global Cartesian system allows to directly define energy functionals or variational forms in the Cartesian sense, making it straightforward to script these functionals in frameworks like NGSolve \cite{Sch1997,Sch2014,Gangl2021}, FEniCS \cite{baratta_2023_10447666,Kuchta}, or Firedrake \cite{FiredrakeUserManual}.

In this work we introduce the linear isotropic micropolar Cosserat model in three dimensions, and subsequently reduce it to shell-, plate- and beam models by means of kinematical assumptions and integration. The models all entail the same kinematical degrees of freedom, namely, displacements $\ud$ and rotations $\rtheta$, and thus intrinsically agree on interfaces of mixed dimensionality $\Xi \subset \R^{d-k}$. 
The coupling itself is then achieved by restricting the bulk fields to codimensional domains using consistent Sobolev trace operators, yielding a mixed-dimensional action functional 
\begin{align}
    \boxed{
    \begin{aligned}
    &I(\ud,\rtheta) - L(\ud,\rtheta)   \to  \min \quad \wrt \quad \{\ud,\rtheta\} \, , &
    I(\ud, \rtheta) &= I_\vol(\ud,\rtheta) + I_\surf(\tr_\surf \ud,\tr_\surf \rtheta) + I_\curv(\tr_\curv \ud,\tr_\curv \rtheta) \, ,
    \\
    && L(\ud,\rtheta) &= L_\vol(\ud, \rtheta) + L_\surf(\tr_\surf \ud, \tr_\surf \rtheta) + L_\curv(\tr_\curv \ud, \tr_\curv \rtheta) \, .
    \end{aligned}
    }
\end{align}
The derivations in this work are done using tangential differential calculus, thus allowing for direct usage in automated solvers of partial differential equations. We employ NGSolve\footnote{www.ngsolve.org} to compute numerical solutions of three coupling examples\footnote{\AS{https://github.com/Askys/NGSolve-Beam-shell-solid} \label{fn:github}}, and finally, we discuss conclusions and outlook. 

We emphasise that \textbf{the resulting reduced models in this work do not coincide with the standard linear Naghdi shell \cite{SCHOLLHAMMER2019172} or the Reissner--Mindlin plate \cite{sky2023reissnermindlin} formulations, nor with the traditional definition of the Cosserat rod \cite{Harsch}} in the case of the beam formulations. This is because the aforementioned models are derived from the energy functional of the corresponding Cauchy continuum theory, as opposed to the energy functional of the Cosserat micropolar continuum that we employ here. However, all the models are referred to as Cosserat type by virtue of their descendance from the three-dimensional Cosserat model.

\subsection{Notation}
The following notation is used throughout this work.
Exceptions to these rules are made clear in the precise context.
\begin{itemize}
    \item Vectors are defined as bold lower-case letters $\vb{v}, \, \bm{\xi} \in \R^d$.
    \item Second order tensors are denoted with bold capital letters $\bm{T}\in \R^{d \times d}$.
    \item Higher-order tensors are designated by the blackboard-bold format $\mathbb{C} \in \R^{d \times d \times d \dots}$.
    \item We denote the Cartesian basis as $\{\vb{e}_1, \, \vb{e}_2, \, \vb{e}_3\}$.
    \item Summation over indices follows the standard rule of repeating indices. Latin indices represent summation over the full dimension, whereas Greek indices define summation over the co-dimension.
    \item The angle-brackets are used to define scalar products of arbitrary dimensions $\con{\vb{v}}{\vb{u}} = v_i u_i$, $\con{\bm{T}}{\bm{F}} = T_{ij}F_{ij}$.
    \item The matrix product is used to indicate all partial-contractions between a higher-order and a lower-order tensor $\bm{T}\vb{v} = T_{ij} v_j \vb{e}_i$, $\mathbb{C}\bm{T} = C_{ijkl}T_{kl}\vb{e}_i \otimes \vb{e}_j$.
    \item The second-order identity tensor is defined via $\one$, such that $\one \vb{v} = \vb{v}$. 
    \item The trace operator reads $\tr \bm{T} = \con{\bm{T}}{\one}$.
    \item A general physical body of some arbitrary dimension $d$ is denoted with $\body \subset \R^d$.
    \item Volumes, surfaces and curves of the physical domain are identified via $\vol$, $\surf$ and $\curv$, respectively. Their counterparts on the reference domain are $\Omega \subset \R^3$, $\omega \subset \R^2$ and $\gamma \subset \R$. 
    \item Tangential and normal vectors on the physical domain are designated by $\vb{t}$ and $\vb{n}$, respectively. On reference domain the respective vectors read $\bm{\tau}$ and $\bm{\nu}$. 
    \item We define the constant space of skew-symmetric second order tensors as $\so(d) = \{ \bm{T} \in \R^{d \times d} \; | \; \bm{T} = -\bm{T}^{T} \}$.
   \item The space is associated with the operators $\skw \bm{T} = (1/2)(\bm{T} - \bm{T}^T) \in \so(d)$, $\Anti \vb{v} = \vb{v} \times \one \in \so(3)$, and its inverse $\axl (\Anti \vb{v}) = \vb{v}$.
   \item The nabla operator is defined as $\nabla = \vb{e}_i \partial_i$.
   \item The left-gradient is given via $\nabla$, such that $\nabla \lambda = \nabla \otimes \lambda$.
    \item The right-gradient is defined for vectors and higher order tensors via $\D$, such that $\D \vb{v} = \vb{v} \otimes \nabla$.
    \item We define the vectorial divergence as $\di \vb{v} = \con{\nabla}{\vb{v}} = \tr (\D \vb{v})$.
    \item The tensor divergence is given by $\Di \bm{T} = \bm{T} \nabla$, implying a single contraction acting row-wise 
    \item The vectorial curl operator reads $\curl \vb{v} = \nabla \times \vb{v}$
    \item For tensors the operator is given by $\Curl \bm{T} = -\bm{T} \times \nabla$, acting row-wise. 
    \item The jump operator is denoted with $\jump{\cdot}$.
\end{itemize}

Further, we introduce the following Hilbert spaces and their respective norms 
    \begin{align}
\Le(\body) &= \{ u : \body \to \mathbb{R} \; | \; \| u \|_\Le < \infty  \} \, , & \|u\|_\Le^2 &= \int_\body u^2 \, \dd \body \, , \notag \\
    \Hone(\body) &= \{ u \in \Le(\body) \; | \; \nabla u \in [\Le(\body)]^d \} \, , & \norm{u}^2_{\Hone} &= \norm{u}^2_\Le + \norm{\nabla u}^2_\Le \, ,  \\[0.5em]
    \Hc{,\body} &= \{ \ud \in [\Le(\body)]^{d} \; | \; \curl \ud \in [\Le(\body)]^{d} \} \, , & \norm{\ud}^2_{\Hc{}} &= \norm{\ud}^2_\Le + \norm{\curl \ud}^2_\Le \, , 
\end{align}
where $\body \subset \R^d$. 

\section{The linear Cosserat micropolar model}
\begin{figure}
		\centering
		\definecolor{asl}{rgb}{0.4980392156862745,0.,1.}
		\definecolor{asb}{rgb}{0.,0.4,0.6}
		\begin{tikzpicture}[line cap=round,line join=round,>=triangle 45,x=1.0cm,y=1.0cm]
			\clip(-0.5,-0.5) rectangle (16,4.5);
			
			\fill [asb, opacity=0.1] plot [smooth cycle] coordinates {(1,3) (3,4) (7, 2) (10,3) (12,1) (10,0) (5,1) (2,1)};
			
			\begin{scope}
				\clip(5,-0.5) rectangle (12.5,1.5);
				\draw [asl, dashed] plot [smooth cycle] coordinates {(1,3) (3,4) (7, 2) (10,3) (12,1) (10,0) (5,1) (2,1)};
			\end{scope}
		    \begin{scope}
		    	\clip(5,1.5) rectangle (12.5,4.5);
		    	\draw [asl, dashed] plot [smooth cycle] coordinates {(1,3) (3,4) (7, 2) (10,3) (12,1) (10,0) (5,1) (2,1)};
		    \end{scope}
			\begin{scope}
				\clip(0,-0.5) rectangle (5,4.5);
				\draw [asb] plot [smooth cycle] coordinates {(1,3) (3,4) (7, 2) (10,3) (12,1) (10,0) (5,1) (2,1)};
			\end{scope}
			
			\draw [-to,color=black,line width=1.pt] (0,0) -- (1,0);
			\draw [-to,color=black,line width=1.pt] (0,0) -- (0,1);
			\draw (1,0) node[color=black,anchor=west] {$x$};
			\draw (0,1) node[color=black,anchor=south] {$y$};
			
			\draw [-to,color=asb,line width=1.pt] (10.,1.2) -- (10.3,0.9);
			\draw [-to,color=asb,line width=1.pt] (9.6,1.2) -- (9.3,0.9);
			\draw [-to,color=asb,line width=1.pt] (9.8,1.2) -- (9.8,0.7);
			
			\draw [-to,color=asl,line width=1.pt] (11.35,2) -- (12,2.55);
			\draw (11.8,2.3) node[color=asl,anchor=north] {$\vb{n}$};
			\draw (9.8,1.2) node[color=asb,anchor=south] {$\vb{f}$};
			\draw (2.8,2.3) node[color=asb,anchor=south] {$\bm{M}$};
			\draw [-to,asb,domain=0:180,line width=1.pt] plot ({0.5*cos(\x-180)+2.8}, {0.5*sin(\x-180)+2.5});
			
			\draw (6.5,1.15) node[color=asb,anchor=south] {$\vol$};
			\draw (4.3,3.55) node[color=asb,anchor=west] {$\surf_D$};
			\draw (10.6,2.95) node[color=asl,anchor=south] {$\surf_N$};
			
			\draw [black,domain=0:360,densely dotted] plot ({0.3*cos(\x)+11.3}, {0.3*sin(\x)+1});
			\draw [black,domain=0:360,densely dotted] plot ({0.78*cos(\x)+14}, {0.78*sin(\x)+1});
			\draw [color=black, densely dotted] (11.3,1.3) -- (14,1.78);
			\draw [color=black, densely dotted] (11.3,0.7) -- (14,0.22);
			

   \draw [color=black, line width=0.7] (13.35, 1+0.375) -- (14.65, 1+0.375);
   \draw [color=black, line width=0.7] (14-0.375, 1+1.3/2) -- (14-0.375, 1-1.3/2);
   \draw [color=black, line width=0.7] (13.25, 1) -- (14.75, 1);
   \draw [color=black, line width=0.7] (13.35, 1-0.375) -- (14.65, 1-0.375);
   \draw [color=black, line width=0.7] (14, 1+0.75) -- (14, 1-0.75);
   \draw [color=black, line width=0.7] (14+0.375, 1+1.3/2) -- (14+0.375, 1-1.3/2);
			
			\fill [asb,domain=0:360, opacity=0.1] plot ({0.78*cos(\x)+14}, {0.78*sin(\x)+1});
			
			
			
			
			
			
		\end{tikzpicture}
		\caption{The domain $\vol \subset \R^3$ of the linear Cosserat micropolar model with Neumann $\surf_N$ and Dirichlet $\surf_D$ boundaries under internal forces $\vb{f}$ and couple-forces $\bm{M}$. The model can be intuitively understood as the superposition of a linear elastic Navier--Cauchy model with an infinitesimally thin fibre-matrix that reinforces the material with respect to bending and torsion.}
		\label{fig:domaincoss}
\end{figure}
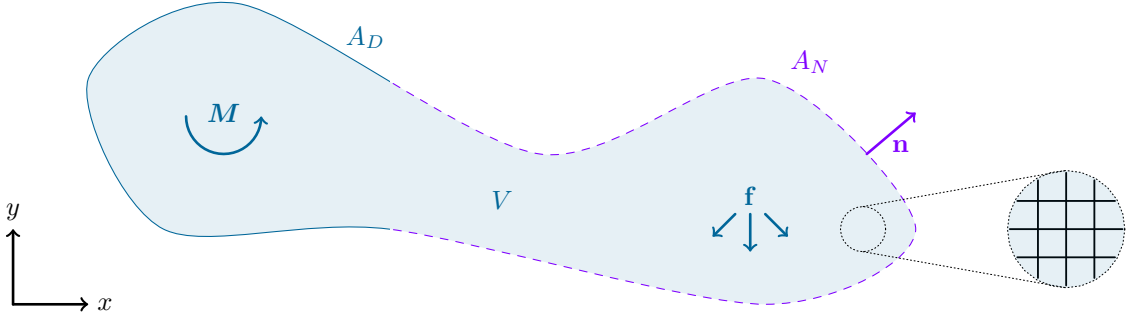

The internal energy of a three-dimensional Cosserat medium is given in dislocation form (Curl instead of gradient) \cite{Ghiba,Jeong1,Jeong2,Jeong3} by the functional
\begin{align}
    I_\vol(\ud, \rtheta) = \dfrac{1}{2}\int_\vol \con{\sym \D \ud}{\Ce \sym \D \ud} + 2 \muc \norm{\skw(\D \ud - \Rtheta)}^2 + \mue \Lc^2 \con{\Curl \Rtheta}{\Lm \Curl \Rtheta} \, \dd \vol \, ,
\end{align}
where $\ud:\vol \subset \R^3 \to \R^3$ is the displacement of the body $\vol \subset \R^3$. Its infinitesimal micro-rotation $\Rtheta:\vol \subset \R^3 \to \so(3)$ is given by the skew-symmetric tensor  
\begin{align}
    &\Rtheta = \Anti \rtheta = \rtheta \times \one \, , && \rtheta = \axl \Rtheta = -\dfrac{1}{2} \mathbb{E} \Rtheta  \, , 
\end{align}
where $\mathbb{E} \in \R^{3 \times 3 \times 3}$ is the third-order Levi-Civita permutation tensor. The infinitesimal micro-rotation tensor is characterised by the axial vector $\rtheta:\vol \subset \R^3 \to \R^3$.
Consequently, its Curl can be written as
\begin{align}
    \Curl \Rtheta = \Curl (\Anti \rtheta) = (\di \rtheta) \one - (\D \rtheta)^T \, ,
    \label{eq:curlskw}
\end{align}
using Nye's formula \cite{Lewintan2021,refId0}.
The material tensor $\Ce:\R^{3 \times 3} \to \R^{3 \times 3}$ is positive definite, and assumed to be isotropic hereinafter
\begin{align}
    &\Ce \bm{T} = 2\mue \bm{T} + \lame (\tr \bm{T})\one \,, && \bm{T}  \in \R^{3 \times 3}  \, ,
\end{align}
where $\mue > 0$ and $\lame \geq 0$ are Lam\'e constants.
The coupling between the infinitesimal macro-rotation $\skw \D \ud$ and the infinitesimal micro-rotation $\Rtheta$ is governed by the Cosserat couple modulus $\muc > 0$. 
The characteristic length-scale parameter is denoted by $\Lc > 0$, and $\Lm:\R^{3 \times 3} \to \R^{3 \times 3}$ is a positive definite isotropic fourth-order tensor of dimensionless weights 
\begin{align}
    &\Lm \bm{T} = a_1 \dev \sym \bm{T} + a_2 \skw \bm{T} + \dfrac{a_3}{3} (\tr \bm{T}) \one \, , && a_1,a_2,a_3 \geq 0 \, , && \bm{T}  \in \R^{3 \times 3}  \, .
\end{align}
The external work for the medium is given by 
\begin{align}
    L_\vol(\ud,\rtheta) = \int_\vol \con{\ud}{\vb{f}} + \con{\Rtheta}{\bm{M}} \, \dd \vol \, ,
\end{align}
where $\vb{f}:\vol\subset \R^3 \to \R^3$ are the body forces and $\bm{M}:\vol\subset \R^3 \to \so(3)$ are couple-forces. For simplicity, we do not consider external fluxes in this work, such that any Neumann boundary is always homogeneous. The balance of energy is expressed as the minimisation problem
\begin{align}
    &I_\vol(\ud,\rtheta) - L_\vol(\ud,\rtheta)   \to  \min \quad \wrt \quad \{\ud,\rtheta\} \, , && \{\widetilde{\ud},\widetilde{\rtheta}\} = \argmin_{\ud,\rtheta} (I_\vol - L_\vol) \, .
    \label{eq:energybalance}
\end{align}
In order to find minimisers $\{\widetilde{\ud},\widetilde{\rtheta}\}$ we take variations with respect to the displacements and infinitesimal rotations. The variation with respect to the displacements $\ud$ yields
\begin{align}
    \delta_\ud (I_\vol-L_\vol) = \int_\vol \con{\sym \D \delta\ud}{\Ce \sym \D \ud} + 2 \muc \con{\skw\D  \delta \ud}{\skw(\D \ud - \Rtheta)} - \con{\delta \ud}{\vb{f}} \, \dd \vol = 0  \, .
\end{align}
The variation with respect to the infinitesimal rotation $\Rtheta$ reads
\begin{align}
    \delta_{\Rtheta} (I_\vol-L_\vol) = \int_\vol -2 \muc \con{\delta \Rtheta}{\skw(\D \ud - \Rtheta)}  + \mue \Lc^2 \con{\Curl \delta\Rtheta}{\Lm \Curl \Rtheta} -\con{\delta \Rtheta}{\bm{M}} \,\dd \vol = 0 \, .
\end{align}
Combining the two results in the weak form
\begin{align}
    \int_\vol \con{\sym \D \delta \ud}{\Ce \sym \D \ud} &+ 2 \muc \con{\skw(\D\delta \ud - \delta\Rtheta)}{\skw(\D \ud - \Rtheta)}   \notag \\  &+ \mue \Lc^2 \con{\Curl \delta\Rtheta}{\Lm \Curl \Rtheta}  \, \dd \vol = \int_\vol \con{\delta \ud}{\vb{f}} + \con{\delta \Rtheta}{\bm{M}} \, \dd \vol \, .
\end{align}
Now, by splitting the boundary between Dirichlet and Neumann $\partial\vol = \surf_D \cup \surf_N$ and applying partial integration we obtain the boundary value problem
\begin{subequations}
    \begin{align}
        -\Di[\Ce \sym\D \ud + 2 \muc \skw(\D \ud - \Rtheta)] &= \vb{f} && \text{in} \quad \vol \, ,  \\
        -2\muc \skw(\D \ud - \Rtheta) + \mue \Lc^2 \Curl (\Lm \Curl \Rtheta) &= \bm{M} && \text{in} \quad \vol \, , \label{eq:momentbalance} \\
        [\Ce \sym\D \ud + 2 \muc \skw(\D \ud - \Rtheta)] \vb{n} &= 0 && \text{on} \quad \surf_N^\ud \, , \\
        -\mue \Lc^2 (\Lm\Curl \Rtheta) (\Anti\vb{n}) &= 0 && \text{on} \quad \surf_N^{\Rtheta} \, , \\
        \ud &= \widetilde{\ud} && \text{on} \quad \surf_D^{\ud} \, , \\
        \Rtheta(\Anti \vb{n}) &= \widetilde{\Rtheta}(\Anti \vb{n}) && \text{on} \quad \surf_D^{\Rtheta} \, .
    \end{align}
\end{subequations}
Accordingly, the non-symmetric stress tensor is 
\begin{align}
    \bm{\sigma} = \Ce \sym \D \ud + 2\muc \skw(\D \ud - \Rtheta) \, .
    \label{eq:cosseratstress}
\end{align}
The domain of the Cosserat model with boundary conditions is depicted in \cref{fig:domaincoss}. 
We note that setting $\muc = 0$ decouples the equations, leaving $\mue \Lc^2 \Curl (\Lm \Curl \Rtheta) = \bm{M}$ as a standalone Curl-Curl problem. In contrast, the coefficient is redundant in the geometrically exact Cosserat model \cite{Ghiba}, and can be set to zero $\muc = 0$ in the related relaxed micromorphic model \cite{Neff2014,SKYNOVEL,SKY2022115298,GOURGIOTIS2024112700} without decoupling the kinematical fields.

\textbf{Hereinafter, we employ material norms for better readability of the shell and beam models, and compactness of notation}. Namely, we define the norms
\begin{align}
    \norm{\sym \D \ud}_{\Ce}^2 = \con{\sym \D \ud}{\Ce \sym \D \ud} = 2\mue\norm{\sym \D \ud}^2 + \lame [\tr(\sym \D \ud)]^2  \, ,
\end{align}
and 
\begin{align}
    \norm{(\di \rtheta) \one - (\D \rtheta)^T}_{\Lm}^2  = \con{\Curl \Rtheta}{\Lm \Curl \Rtheta} =  a_1 \norm{\dev \sym \Curl\Rtheta}^2 + a_2 \norm{\skw \Curl\Rtheta}^2 + \dfrac{a_3}{3} (\tr \Curl\Rtheta)^2 \, ,
\end{align}
where we used the identity \cref{eq:curlskw} for skew-symmetric tensors. Thus, the internal energy functional reads
\begin{align}
    \boxed{
    \begin{aligned}
    I_\vol(\ud, \rtheta) = \dfrac{1}{2}\int_\vol \norm{\sym \D \ud}_{\Ce}^2  + 2 \muc \norm{\skw(\D \ud - \Rtheta)}^2 + \mue \Lc^2 \norm{(\di \rtheta) \one - (\D \rtheta)^T}_{\Lm}^2 \, \dd \vol \, .
    \label{eq:energy3d}
\end{aligned}
    }
\end{align}
The functional is composed of three energy densities pertaining to the displacement, rotational coupling, and dislocation, respectively.  

\subsection{Linear isotropic Cauchy materials}
The Cosserat micropolar model intrinsically incorporates the classical Navier--Cauchy linear elasticity model into the structure of its energy functional. Namely, a highly homogeneous Cauchy medium is characterised by a vanishing characteristic length-scale parameter $\Lc \to 0$, and the absence of couple-forces $\bm{M} = 0$ \cite{NeffWieners}. Thus, \cref{eq:momentbalance} reads
\begin{align}
    - 2\muc \skw(\D \ud - \Rtheta) = 0 \qquad \overset{\muc > 0}{\iff} \qquad \Rtheta = \skw \D \ud \, ,
\end{align}
such that the internal energy functional yields
\begin{align}
     \dfrac{1}{2}\int_\vol \norm{\sym \D \ud}_{\Ce}^2  + 2 \muc \norm{\skw(\D \ud - \Rtheta)}^2 + \mue \Lc^2 \norm{(\di \rtheta) \one - (\D \rtheta)^T}_{\Lm}^2 \, \dd \vol \quad \mapsto \quad 
    \dfrac{1}{2}\int_\vol \norm{\sym \D \ud}_{\Cm_{\mathrm{M}}}^2  \, \dd \vol \, ,
\end{align}
where $\Cm_{\mathrm{M}}: \R^{3 \times 3} \to \R^{3 \times 3}$ is the material tensor of linear isotopic Cauchy materials
\begin{align}
    &\Cm_{\mathrm{M}} \bm{T} = 2\mua\bm{T} + \lama(\tr\bm{T})\one \, , && \bm{T} \in \R^{3 \times 3} \, ,
\end{align}
with the Lam\'e constants $\mua > 0$ and $\lama \geq 0$. In summary, it suffices to set $\Lc \to 0$, $\bm{M} = 0$, $\muc > 0$ and $\Ce = \Cm_{\mathrm{M}}$ in order to obtain the behaviour of the Navier--Cauchy model from the Cosserat model.

\section{Linear isotropic Cosserat shells in three dimensions}
\begin{figure}
    \centering
    \definecolor{asb}{rgb}{0.4980392156862745,0.,1.}
\definecolor{asl}{rgb}{0.,0.4,0.6}
\begin{tikzpicture}[scale = 0.4, line cap=round,line join=round,>=triangle 45,x=1.0cm,y=1.0cm]
\clip(1.5,0) rectangle (40,11.5);
\fill[line width=0.7pt,color=asl,fill=asl,fill opacity=0.10000000149011612] (2.,3.) -- (8.,1.) -- (14.,3.) -- (8.,5.) -- cycle;
\draw [line width=0.7pt,color=asl] (2.,3.)-- (8.,1.);
\draw [line width=0.7pt,color=asl] (8.,1.)-- (14.,3.);
\draw [line width=0.7pt,color=asl] (14.,3.)-- (8.,5.);
\draw [line width=0.7pt,color=asl] (8.,5.)-- (2.,3.);
\draw [line width=0.7pt,dashed,color=asl] (2.,4.)-- (8.,6.);
\draw [line width=0.7pt,dashed,color=asl] (8.,6.)-- (14.,4.);
\draw [line width=0.7pt,dashed,color=asl] (2.,4.)-- (8.,2.);
\draw [line width=0.7pt,dashed,color=asl] (14.,4.)-- (8.,2.);
\draw [line width=0.7pt,dashed,color=asl] (2.,4.)-- (2.,2.);
\draw [line width=0.7pt,dashed,color=asl] (2.,2.)-- (8.,0.);
\draw [line width=0.7pt,dashed,color=asl] (8.,0.)-- (14.,2.);
\draw [line width=0.7pt,dashed,color=asl] (14.,2.)-- (14.,4.);
\draw [-to,line width=1pt] (2.,3.) -- (5.3,4.1);
\draw [-to,line width=1pt] (2.,3.) -- (5.,2.);
\draw [-to,line width=1pt] (2.,3.) -- (2.,6.);
\draw [-to,line width=1pt] (20.,3.) -- (20.,6.);
\draw [-to,line width=1pt] (20.,3.) -- (23.,4.);
\draw [-to,line width=1pt] (20.,3.) -- (23.,2.);
\draw [line width=0.7pt, asl, dashed] plot [smooth cycle] coordinates {(26,7) (32,11) (38, 9) (40,5) (38,3) (34,5) (30,7)};
\draw [line width=0.7pt, asl] plot [smooth cycle] coordinates {(26,7-1) (32,11-1) (38, 9-1) (40,5-1) (38,3-1) (34,5-1) (30,7-1)};
\fill [asl, opacity=0.1] plot [smooth cycle] coordinates {(26,7-1) (32,11-1) (38, 9-1) (40,5-1) (38,3-1) (34,5-1) (30,7-1)};
\draw [line width=0.7pt, asl, dashed] plot [smooth cycle] coordinates {(26,7-2) (32,11-2) (38, 9-2) (40,5-2) (38,3-2) (34,5-2) (30,7-2)};
\draw [line width=0.7pt,dashed,color=asl] (26.,5.)-- (26.,7.);
\draw [line width=0.7pt,dashed,color=asl] (37.99999998282016,0.9999999956690431)-- (38.,3.);
\draw [line width=0.7pt,dashed,color=asl] (39.999999930006084,2.9999997424275477)-- (40.,5.);
\draw [-to,line width=1pt] (15.,8.) -- (19.,8.);
\draw [-to,line width=1pt,color=asb] (4.,3.) -- (4.,4.666666666666667);
\draw [-to,line width=1pt,color=asb] (30.973993776118732,5.740820382058072) -- (31.401161525645392,6.579144483189179);
\draw [-to,line width=1pt,color=asb] (20.,3.) -- (30.973993776118732,5.740820382058072);
\draw (5,2.01) node[anchor=north west] {$\xi$};
\draw (5,4.1) node[anchor=south west] {$\eta$};
\draw (2,6) node[anchor=south] {$\zeta$};
\draw [color=asl](6.8,4.2) node[anchor=north west] {$\omega \subset \R^2$};
\draw [color=asl](11,6.4) node[anchor=north west] {$\Omega \subset \R^3$};
\draw (23,2) node[anchor=north west] {$x$};
\draw (23,4.2) node[anchor=south west] {$y$};
\draw (20,6) node[anchor=south] {$z$};
\draw (17,8) node[anchor=south] {$\mathbf{x}:\Omega \to V$};
\draw [color=asb](24.5,4.3) node[anchor=north west] {$\mathbf{r}:\omega \to A$};
\draw [color=asb](31.3,7.5) node[anchor=north west] {$\mathbf{n}$};
\draw [color=asb](4,4.8) node[anchor=south] {$\boldsymbol{\nu}$};
\draw [color=asl](34,7.3) node[anchor=north west] {$A \subset \R^3$};
\draw [color=asl](36.5,11) node[anchor=north west] {$V \subset \R^3$};
\begin{scriptsize}
\draw [color=asb] (4.,3.)-- ++(-2.5pt,-2.5pt) -- ++(5.0pt,5.0pt) ++(-5.0pt,0) -- ++(5.0pt,-5.0pt);
\draw [color=asb] (30.973993776118732,5.740820382058072)-- ++(-2.5pt,-2.5pt) -- ++(5.0pt,5.0pt) ++(-5.0pt,0) -- ++(5.0pt,-5.0pt);
\end{scriptsize}
\end{tikzpicture}
    \caption{Mapping of the flat reference domain to a shell in three-dimensional space $\vb{x}:\Omega \subset \R^3 \to \vol \subset \R^3$. The middle surface of the shell is mapped via $\vb{r}:\omega \subset \R^2 \to \surf \subset \R^3$, and the thickness via $\zeta \vb{n}$ with $\zeta \in [-h/2, h/2]$.}
    \label{fig:shell}
\end{figure}
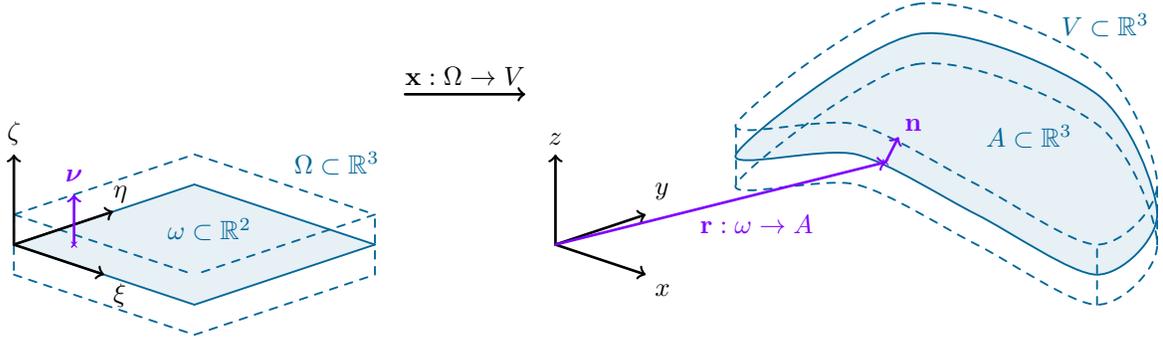

Let the domain of a thin curved shell be given by the mapping
\begin{align}
    &\vb{x}:\Omega \subset \R^3 \to \vol \subset \R^3 \, , && \Omega = \omega \times [-h/2,h/2] \, , && \vol = \surf \times [-h/2,h/2] \, ,
\end{align}
where $h \ll |\surf|$ is the thickness of the shell,
we can write the mapping explicitly as
\begin{align}
    &\vb{x}(\xi,\eta,\zeta) = \vb{r} + \zeta\vb{n} \, , && \vb{r} = \vb{r}(\xi,\eta) \, , && \zeta \in [-h/2,/h/2] \, ,  
\end{align}
such that $\vb{r}$ maps the middle surface of the shell $\vb{r}:\omega \subset \R^2 \to \surf \subset \R^3$, and $\vb{n} = \vb{n}(\xi,\eta)$ is the unit normal vector of the surface, see \cref{fig:shell}.
With the surface normal we can define corresponding tangential and normal projection operators
\begin{align}
    &\Pt = \one - \vb{n} \otimes \vb{n} \, , &&  \Pn = \vb{n} \otimes \vb{n} \, , && \Pt:\R^3 \to \tansurf \surf  \, , && \Pn:\R^3 \to \R^3 \setminus \tansurf \surf \, ,
\end{align}
where $\tansurf \surf$ is the space of the tangential vectors of the surface $\surf$. Using the latter we split the symmetrised gradient of the displacements between its tangential and normal components 
\begin{align}
    \sym \D \ud &=  (\Pt + \Pn) \sym  \D \ud (\Pt + \Pn) =  \Pt(\sym\D \ud)\Pt + \Pt(\sym\D \ud)\Pn  + \Pn(\sym\D \ud)\Pt + \Pn(\sym\D \ud)\Pn \, .     
\end{align}
Consequently, the first energy density component in \cref{eq:energy3d} can be expressed as
\begin{align}
    \norm{\sym \D \ud}_{\Ce}^2 = &\norm{\Pt(\sym \D \ud)\Pt + \Pn(\sym \D \ud)\Pn}_{\Ce}^2 + 4\mue \norm{\Pn (\sym\D \ud)\Pt}^2  \, ,
    \label{eq:decompDu}
\end{align}
where $\norm{\Pn (\sym\D \ud)\Pt} = \norm{\Pt(\sym\D \ud) \Pn}$ due to symmetry.
Analogously, we split the infinitesimal rotation tensor
\begin{align}
    \Rtheta = (\Pt + \Pn) \Rtheta (\Pt + \Pn) = \Pt \Rtheta \Pt + \Pt \Rtheta \Pn + \Pn \Rtheta \Pt \, ,  
\end{align}
where $\Pn \Rtheta \Pn = \con{\Rtheta}{\Pn} \Pn = 0$ due to its skew-symmetry $\Rtheta : \vol \subset \R^3 \to \so(3)$.
As such, the second energy density term is expanded to  
\begin{align}
    2 \muc \norm{\skw(\D \ud - \Rtheta)}^2 = 2\muc \norm{ \Pt(\skw \D \ud  - \Rtheta) \Pt }^2 + 4 \muc \norm{\Pn (\skw \D \ud  -  \Rtheta) \Pt }^2  \, ,
    \label{eq:decompSkewDuO}
\end{align}
where $\norm{\Pn (\skw \D \ud  -  \Rtheta) \Pt } = \norm{\Pt (\skw  \D\ud  -  \Rtheta) \Pn }$ and $\Pn (\skw \D \ud) \Pn = (1/2) \Pn (\D \ud -  [\D \ud]^T)\Pn = 0$ hold due to skew-symmetry. 
Putting it all together, the internal energy functional reads
\begin{align}
    I_\surf(\ud,\rtheta) =\dfrac{1}{2} \int_\vol &\norm{\Pt(\sym \D \ud)\Pt + \Pn(\sym \D \ud)\Pn}_{\Ce}^2 + 4\mue \norm{\Pn (\sym\D \ud)\Pt}^2 + 4 \muc \norm{\Pn (\skw \D \ud  -  \Rtheta) \Pt }^2 
    \notag \\
    & + 2\muc \norm{ \Pt(\skw \D \ud  - \Rtheta) \Pt }^2 + \mue \Lc^2\norm{(\di \rtheta) \one - (\D \rtheta)^T}_{\Lm}^2  \, \dd \vol  \, . 
\end{align}
Integration of the energy over the volume is split across the middle surface and the thickness of the shell via
\begin{align}
    \int_\vol (\cdot) \, \dd \vol = \int_\surf \int_{-h/2}^{h/2} (\cdot) (1 -2H \, \zeta + K \, \zeta^2) \, \dd \zeta \dd \surf \, ,
\end{align}
using the shell-shifter from Steiner's formula \cite{Ghiba2021,Ghiba2023,Ghiba2020I,Ghiba2020II}, see also \cref{ap:shellshifter}. The factor is composed of the mean and Gauß curvatures 
\begin{align}
    &H = \dfrac{1}{2} \tr \bm{W} \, , &&
    K = \det \bm{W} \, ,
\end{align}
being invariants of the Weingarten curvature tensor 
\begin{align}
    \bm{W} = - \D_t \vb{n} \, ,
\end{align}
whose derivation can be found in \cref{ap:weingarten}. Thus, the energy functional is finally given as 
\begin{align}
    I_\surf(\ud,\rtheta) = \dfrac{1}{2}\int_\surf \int_{-h/2}^{h/2} &(\norm{\Pt(\sym \D \ud)\Pt + \Pn(\sym \D \ud)\Pn}_{\Ce}^2 + 4\mue \norm{\Pn (\sym\D \ud)\Pt}^2 + 4 \muc \norm{\Pn (\skw \D \ud  -  \Rtheta) \Pt }^2 
    \notag \\
    & + 2\muc \norm{ \Pt(\skw \D \ud  - \Rtheta) \Pt }^2  + \mue \Lc^2\norm{(\di \rtheta) \one - (\D \rtheta)^T}_{\Lm}^2 ) (1 -2H \, \zeta + K \, \zeta^2) \, \dd \zeta \dd \surf  \, . 
    \label{eq:shellenergypart1}
\end{align}

\subsection{The shell model} \label{sec:shell}
At this point we make standard engineering assumptions for the kinematics of the shell. Firstly, the infinitesimal rotation is constant throughout the thickness of the shell $\Rtheta \neq \Rtheta(\zeta)$. Secondly, the displacement vector takes the form
\begin{align}
    &\ud(\xi,\eta,\zeta) = \vb{v} + \zeta \Rtheta \vb{n} \, , && \vb{v} = \vb{v}(\xi,\eta) \,, && \Rtheta = \Rtheta(\xi,\eta) \, , 
\end{align}
where $\vb{v}:\omega \subset \R^2 \to \R^3$ is the displacement of the middle surface. In other words, no out-of-plane stretching is permitted, and displacements parallel to the middle surface $\Pt \vb{v}:\omega \subset \R^2 \to \tansurf \surf$ are given by the translation of the middle surface $\vb{v}$ together with its rotation $\zeta \Rtheta \vb{n}$. Observe that $\Rtheta = \Anti \rtheta$ is still given by the three-dimensional axial vector $\rtheta:\omega \subset \R^2 \to \R^3$, such that constant drill throughout the thickness is possible. 
Using the projection operators we define tangential gradients \textbf{with respect to the parametrisation of the middle surface} $\vb{r}:\omega \subset \R^2 \to \surf \subset \R^3$
\begin{align}
    \nabla_t(\cdot) = \Pt \nabla (\cdot) \, , &&  \D_t (\cdot) = [\D(\cdot)]\Pt \, .
\end{align}
The gradient of the displacement field can now be split across its tangential gradient and normal components via \cref{eq:splittangrad}
\begin{align}
    \D \ud = \D \vb{v} + \D(\zeta \Rtheta \vb{n}) = \D_t \vb{v} + (\Rtheta \vb{n}) \otimes \nabla\zeta + \zeta\D_t(\Rtheta \vb{n}) \, .
\end{align}
Thus, for the tangential part we find
\begin{align}
    (\D\ud) \Pt = \D_t \ud = \D_t \vb{v} + \zeta [\D_t (\Rtheta \vb{n})] =   \D_t \vb{v} + \zeta [\D_t (\rtheta \times \vb{n})] = \D_t \vb{v} - \zeta [(\Anti \vb{n}) \D_t \rtheta + \Rtheta \bm{W}] \, ,  
\end{align}
where $\bm{W}= -\D_t \vb{n}$ is the Weingarten curvature tensor,
which is naturally tangential $\bm{W} = \Pt \bm{W} \Pt$ as per \cref{ap:weingarten}.
Since $(\Anti \vb{n})\vb{n} = 0$, we can further replace $(\Anti \vb{n})\D_t \rtheta$ with $(\Anti \vb{n})\Dcov \rtheta$ where
\begin{align}
    \Dcov \rtheta = \Pt \D_t \rtheta \, ,
\end{align}
is the covariant gradient, finally yielding
\begin{align}
    \D_t \ud = \D_t \vb{v} - \zeta [(\Anti \vb{n}) \Dcov \rtheta + \Rtheta \bm{W}] \, .  
\end{align}
The normal part of the gradient is simply
\begin{align}
    &(\D \ud) \Pn = (\Rtheta \vb{n}) \otimes \nabla \zeta = \Rtheta \vb{n} \otimes \vb{n} = \Rtheta \bm{Q} \,  ,
\end{align}
and for the micro-dislocation $\Curl \Rtheta$ we find
\begin{align}
    &\Curl \Rtheta = (\di_t \rtheta) \one - (\D_t \rtheta)^T \, , && \di_t \rtheta = \tr(\D_t \rtheta) \, ,
\end{align}
seeing as $\Rtheta = (\Rtheta \circ \vb{r})(\xi,\eta)$.
Next, we observe that $\Pt (\sym \D \ud) \Pt =  \sym (\Pt[\D \ud] \Pt)$ yields
\begin{align}
    \Pt (\sym \D \ud) \Pt = \sym(\Dcov \vb{v} -  \zeta [(\Anti \vb{n}) \Dcov \rtheta + \Pt \Rtheta \bm{W}]) \, . 
\end{align}
Analogously, we find
\begin{align}
    \Pt (\skw \D \ud - \Rtheta) \Pt = \skw(\Dcov \vb{v} -  \zeta [(\Anti \vb{n}) \Dcov \rtheta + \Pt \Rtheta \bm{W}]) - \Pt\Rtheta\Pt \, .
\end{align}
Further, there holds 
\begin{align}
    \Pn (\sym \D \ud)\Pt = \frac{1}{2}(\Pn [\D \ud]\Pt + \Pn [\D \ud]^T\Pt) = \dfrac{1}{2}(\Pn \D_t \vb{v} - \zeta \Pn \Rtheta \bm{W}) - \dfrac{1}{2} \Pn \Rtheta \Pt \, ,
\end{align}
such that the skew-symmetric counterpart satisfies
\begin{align}
     \Pn (\skw \D \ud - \Rtheta)\Pt = \dfrac{1}{2}(\Pn \D_t \vb{v} - \zeta \Pn \Rtheta \bm{W}) + \dfrac{1}{2} \Pn \Rtheta \Pt -  \Pn \Rtheta \Pt = \dfrac{1}{2}(\Pn \D_t \vb{v} - \zeta \Pn \Rtheta \bm{W}) - \dfrac{1}{2} \Pn \Rtheta \Pt  \, ,
\end{align}
where we use that $\Rtheta^T = -\Rtheta$. Lastly, the normal-normal part of the symmetrised gradient yields
\begin{align}
    \Pn (\sym \D \ud )\Pn =  \sym (\Pn \D \ud\Pn ) = \Pn\Rtheta \Pn = 0\, ,
\end{align}
since $\Pn\Rtheta\Pn = \con{\Rtheta}{\Pn} \Pn = 0$. Consequently we can simplify the energy density terms of \cref{eq:shellenergypart1}. We start by observing that
\begin{align}
    \norm{\Pt(\sym \D \ud)\Pt + \Pn(\sym \D \ud)\Pn}_{\Ce}^2 = \norm{\sym(\Dcov \vb{v} -  \zeta [(\Anti \vb{n}) \Dcov \rtheta + \Pt \Rtheta \bm{W}])}_{\De}^2 \, .
\end{align}
Here, we get that the normal-normal strain component vanishes. Following the classical Naghdi shell and corresponding Reissner--Mindlin plate formulations, this leads to an asymptotically improper model since an out-of-plane stress $\Pn \bm{\sigma}\Pn$ arises from transverse contractions via the Poisson ratio \cite{MichaelThesis}. Thus, we adopt the plane-stress assumption, such that the normal-normal component of the strain is condensed a priori and the corresponding stress component is set to vanish $\Pn \bm{\sigma}\Pn = 0$. For the planar material tensor $\De$ this implies that its material constants are adapted to those of plane stress
\begin{align}
    & \norm{\cdot}_{\De}^2 = 2\mue^*\norm{\cdot}^2 + \lame^*[\tr(\cdot)]^2  \, ,
 && \mue^* = \dfrac{E_e}{2(1+\nu_e)} = \mue \, , && \lame^* = \dfrac{E_e\,\mue}{1-\nu_e^2} = \dfrac{2 \lame\, \mue}{\lame + 2 \mue} \, ,
\end{align}
where $E_e$ and $\nu_e$ are Young's modulus and the Poisson ratio, respectively.
Observing that $\Pn (\sym \D \ud)\Pt = \Pn (\skw \D \ud - \Rtheta)\Pt$, the next two energy density terms are combined into
\begin{align}
     4\mue \norm{\Pn (\sym\D \ud)\Pt}^2 + 4 \muc \norm{\Pn (\skw \D \ud  -  \Rtheta) \Pt }^2 = (\mue + \muc) \norm{\Pn(\D_t \vb{v} - \zeta  \Rtheta \bm{W} -   \Rtheta \Pt)}^2 \, .
\end{align}
For the skew-symmetric tangential-tangential energy term we find
\begin{align}
     2\muc \norm{ \Pt(\skw \D \ud  - \Rtheta) \Pt }^2 = 2\muc \norm{\skw(\Dcov \vb{v} -  \zeta [(\Anti \vb{n}) \Dcov \rtheta + \Pt \Rtheta \bm{W}]) - \Pt\Rtheta\Pt}^2 \, .
\end{align}
Finally, the energy density of the micro-dislocation reads
\begin{align}
    \mue \Lc^2\norm{(\di_t \rtheta) \one - (\D_t \rtheta)^T}_{\Lm}^2 \, .
\end{align}
Thus, the total energy is given by
\begin{align}
    I_\surf(\vb{v},\rtheta) = \dfrac{1}{2}\int_\surf \int_{-h/2}^{h/2} &(\norm{\sym(\Dcov \vb{v} -  \zeta [(\Anti \vb{n}) \Dcov \rtheta + \Pt \Rtheta \bm{W}])}_{\De}^2 + (\mue + \muc) \norm{\Pn(\D_t \vb{v} - \zeta  \Rtheta \bm{W} -   \Rtheta \Pt)}^2 
    \notag \\
    & + 2\muc \norm{\skw(\Dcov \vb{v} -  \zeta [(\Anti \vb{n}) \Dcov \rtheta + \Pt \Rtheta \bm{W}]) - \Pt\Rtheta\Pt}^2 
    \label{eq:volenergshell}
    \\
    \phantom{\dfrac{1}{2}} & \phantom{+ 2\muc \|} + \mue \Lc^2\norm{(\di_t \rtheta) \one - (\D_t \rtheta)^T}_{\Lm}^2)   (1 -2H \, \zeta + K \, \zeta^2) \, \dd \zeta \dd \surf  \, . \notag
\end{align}
Having derived the total energy, we employ asymptotic analysis to obtain the final shell model by eliminating $\zeta$-terms of order three or higher $\O(\zeta^3)$, see \cref{ap:asym}. These terms result in constants of order $\O(h^4)$ in the thickness $h$, which is imperatively required to be very small $h \ll |\surf|$ for a well-defined formulation \cite{MichaelThesis}. Thus, these terms are expected to generate insignificantly small energies for a thin shell and are therefore omitted.
Further, we exploit that linear terms of the form of $(\cdot)\zeta$ vanish by the symmetry of integration over the thickness of the shell $\int_{-h/2}^{h/2}(\cdot)\zeta\,\dd \zeta = 0$. 
We integrate the energy densities multiplied with the shell-shifter $1 -2H \, \zeta + K \, \zeta^2$ over the thickness in three steps. Firstly, we find the integral over the constant $1$ to be
\begin{align}
    & h \norm{\sym\Dcov \vb{v}}_{\De}^2 + \dfrac{h^3}{12} \norm{\sym ([\Anti \vb{n}] \Dcov \rtheta + \Pt \Rtheta \bm{W})}_{\De}^2 + 
    (\mue + \muc) (h\norm{\Pn(\D_t \vb{v} -   \Rtheta \Pt)}^2 + \dfrac{h^3}{12}\norm{\Pn\Rtheta \bm{W}}^2) 
    \notag \\
    & + 2\muc (h\norm{\skw\Dcov \vb{v} - \Pt\Rtheta\Pt}^2 + \dfrac{h^3}{12}\norm{\skw( [\Anti \vb{n}] \Dcov \rtheta + \Pt \Rtheta \bm{W})}^2) + \mue \Lc^2h\norm{(\di_t \rtheta) \one - (\D_t \rtheta)^T}_{\Lm}^2   \, . 
\end{align}
Secondly, the linear term of the shell-shifter $-2H \, \zeta$ leads to
\begin{align}
     \dfrac{H h^3}{3} (& \con{\sym\Dcov \vb{v}}{\De\sym([\Anti \vb{n}] \Dcov \rtheta + \Pt \Rtheta \bm{W})} 
     + (\mue + \muc) \con{\Pn(\D_t \vb{v} -   \Rtheta \Pt)}{\Pn\Rtheta \bm{W} }  
    \notag \\
    & + 2\muc \con{\skw\Dcov \vb{v} - \Pt\Rtheta\Pt}{\skw([\Anti \vb{n}] \Dcov \rtheta + \Pt \Rtheta \bm{W})} ) \, . 
\end{align}
Thirdly, the quadratic term $K \, \zeta^3$ yields
\begin{align}
    \dfrac{K h^3}{12} (& \norm{\sym\Dcov \vb{v}}_{\De}^2 + (\mue + \muc) \norm{\Pn(\D_t \vb{v} -   \Rtheta \Pt)}^2 
    \\
    & + 2\muc \norm{\skw\Dcov \vb{v}  - \Pt\Rtheta\Pt}^2 
     + \mue \Lc^2\norm{(\di_t \rtheta) \one - (\D_t \rtheta)^T}_{\Lm}^2)    \, . \notag
\end{align}
A common simplification at this point can be made for shells with a relatively small curvature, $H/ |\surf| \ll 1$ and $K/ |\surf| \ll 1$. Therein, also the terms $H h^3$ and $K h^3$ become insignificantly small $\O(h^4)$ and are thus omitted from the internal energy, yielding the final internal energy functional of the curved shell 
\begin{align}
    \boxed{
    \begin{aligned}
    I_\surf(\vb{v},\rtheta) = \dfrac{1}{2}\int_\surf  &  h (\norm{\sym\Dcov \vb{v}}_{\De}^2 + 2\muc \norm{\skw\Dcov \vb{v} - \Pt\Rtheta\Pt}^2+ (\mue + \muc) \norm{\Pn(\D_t \vb{v} -   \Rtheta \Pt)}^2 )  \\
    &+ \dfrac{h^3}{12} (\norm{\sym ([\Anti \vb{n}] \Dcov \rtheta + \Pt \Rtheta \bm{W})}_{\De}^2 + 2\muc \norm{\skw( [\Anti \vb{n}] \Dcov \rtheta + \Pt \Rtheta \bm{W})}^2 ) 
     \\
    & \phantom{+ \dfrac{h^3}{12} ( \|} + (\mue + \muc)   \dfrac{h^3}{12}\norm{\Pn\Rtheta \bm{W}}^2  + \mue \Lc^2h\norm{(\di_t \rtheta) \one - (\D_t \rtheta)^T}_{\Lm}^2 \, \dd \surf  \, . 
\end{aligned}
    }
\end{align}
The energy density terms in the functional can be interpreted with respect to their mechanical action. The first term of the functional gives the so-called membrane energy. By \cref{eq:curlDcov} it is apparent that the second term governs in-plane torque and is therefore the drill energy density \cite{SaemDrill}, which does \textbf{not} appear in the Naghdi shell \cite{SCHOLLHAMMER2019172} or Reissner--Mindlin plate \cite{sky2023reissnermindlin} models. The third term is given by the out-of-plane coefficients of the in-plane tensors, yielding the shear energy. The next three terms represent the bending energy, and the final term is the dislocation energy. 
Now, in order to construct a minimisation functional we simply redefine the forces and couple-forces as surface quantities 
\begin{align}
    L_\surf(\vb{v},\rtheta) = \int_\surf \con{\vb{v}}{\vb{q}} + \con{\Rtheta}{\bm{M}} \, \dd \surf \, ,  
\end{align}
such that $\vb{q}:\surf \to \R^3$ and $\bm{M}:\surf \to \so(3)$. The balance of energy is now expressed analogously to \cref{eq:energybalance} via 
\begin{align}
    &I_\surf(\vb{v},\rtheta) - L_\surf(\vb{v},\rtheta)   \to  \min \quad \wrt \quad \{\vb{v},\rtheta\} \, , && \{\widetilde{\vb{v}},\widetilde{\rtheta}\} = \argmin_{\vb{v},\rtheta} (I_\surf - L_\surf) \, ,
    \label{eq:energybalanceshell}
\end{align}
such that minimisers are found by standard variation with respect to $\vb{v}$ and $\Rtheta$.

\subsection{The plate model}
A plate is defined as a flat shell model. Therefore, its kinematical assumptions are similar as for the curved shell but we have that $\bm{W} = 0$.
Accordingly, we get that the mean curvature $H = (1/2)\tr \bm{W} = 0$ and the Gaussian curvature $K = \det \bm{W} = 0$ also vanish. Consequently, the internal energy functional is reduced to 
\begin{align}
    \boxed{
    \begin{aligned}
    I_\surf(\vb{v},\rtheta) = \dfrac{1}{2}\int_\surf  &  h (\norm{\sym\Dcov \vb{v}}_{\De}^2 + 2\muc \norm{\skw\Dcov \vb{v} - \Pt\Rtheta\Pt}^2+ (\mue + \muc) \norm{\Pn(\D_t \vb{v} -   \Rtheta \Pt)}^2 )  \\
    &+ \dfrac{h^3}{12} (\norm{\sym ([\Anti \vb{n}] \Dcov \rtheta)}_{\De}^2 + 2\muc \norm{\skw( [\Anti \vb{n}] \Dcov \rtheta)}^2 ) 
    \\ \phantom{\dfrac{h^3}{12}}& \phantom{+ \dfrac{h^3}{12} (\|} + \mue \Lc^2 h\norm{(\di_t \rtheta) \one - (\D_t \rtheta)^T}_{\Lm}^2 \, \dd \surf  \, . 
\end{aligned}
    }
\end{align}
The form can be even further simplified if the plate is embedded in the $x-y$ or $y-z$ planes.

\subsection{A membrane-shell model} \label{sec:membrane-shell}
In the case of an extremely thin shell $h / |\surf| \ll 1$, the term $h^3$ becomes negligible $\O(h^3)$, such that the energy functional reduces to  
\begin{align}
    \boxed{
    \begin{aligned}
    I_\surf(\vb{v},\rtheta) = \dfrac{1}{2}\int_\surf  h (&\norm{\sym\Dcov \vb{v}}_{\De}^2 + 2\muc \norm{\skw\Dcov \vb{v} - \Pt\Rtheta\Pt}^2
     \\
    &+ (\mue + \muc) \norm{\Pn(\D_t \vb{v} -   \Rtheta \Pt)}^2  
    + \mue \Lc^2\norm{(\di_t \rtheta) \one - (\D_t \rtheta)^T}_{\Lm}^2) \, \dd \surf  \, . 
\end{aligned}
    }
\end{align}
Interestingly, the functional naturally discards any explicit energy terms related to curvature. However, the Weingarten tensor $\bm{W}$ and its invariants still influence the energy implicitly, which can be directly observed by taking derivatives of the kinematic fields defined on the curvilinear coordinates of the shell.

\section{Linear three-dimensional Cosserat beams}
\begin{figure}
    \centering
    \definecolor{xfqqff}{rgb}{0.4980392156862745,0.,1.}
\definecolor{qqwwzz}{rgb}{0.,0.4,0.6}
\begin{tikzpicture}[scale = 0.275, line cap=round,line join=round,>=triangle 45,x=1.0cm,y=1.0cm]
\clip(9,6.5) rectangle (70,21);
\draw [line width=0.7pt,dashed,color=qqwwzz,fill=qqwwzz,fill opacity=0.05000000074505806] (18.,12.) circle (0.9883379552976526cm);
\draw [line width=0.7pt,dashed,color=qqwwzz,fill=qqwwzz,fill opacity=0.05000000074505806] (56.991816300346315,16.97026995335137) circle (1.cm);
\draw [line width=0.7pt,color=qqwwzz] (10.,14.)-- (30.,9.);
\draw [line width=0.7pt,dashed,color=qqwwzz] (10.,13.)-- (30.,8.);
\draw [line width=0.7pt,dashed,color=qqwwzz] (10.,14.) circle (1.cm);
\draw [line width=0.7pt,dashed,color=qqwwzz] (30.,9.) circle (1.cm);
\draw [line width=0.7pt,dashed,color=qqwwzz] (9.,14.)-- (29.,9.);
\draw [line width=0.7pt,dashed,color=qqwwzz] (10.,15.)-- (30.,10.);
\draw [-to, line width=1pt] (42.,8.) -- (46.,7.);
\draw [-to, line width=1pt] (42.,8.) -- (46.,9.);
\draw [-to, line width=1pt] (42.,8.) -- (42.,12.);
\draw[line width=0.7pt,color=qqwwzz, smooth,samples=100,domain=0.0:0.3890404923960828] plot[parametric] function{-0.5648396714375444*t**(3.0)+23.896112147974378*t+45.99999999999999,-21.013000507111983*t**(3.0)+6.565362214806673*t+16.000000000000004};
\draw[line width=0.7pt,color=qqwwzz, smooth,samples=100,domain=0.3890404923960828:0.7559204984484887] plot[parametric] function{-12.904253994854885*t**(3.0)+14.401595472784676*t**(2.0)+18.293308353953034*t+46.72657251560823,11.314676134580687*t**(3.0)-37.73032571611638*t**(2.0)+21.243986709669176*t+14.096473566273833};
\draw[line width=0.7pt,color=qqwwzz, smooth,samples=100,domain=0.7559204984484887:1.0] plot[parametric] function{20.296904312505433*t**(3.0)-60.8907129375163*t**(2.0)+75.20830765680509*t+32.38550096820579,16.485528648217176*t**(3.0)-49.456585944651536*t**(2.0)+30.108107186560165*t+11.8629501098742};
\draw [line width=0.7pt,dashed,color=qqwwzz] (46.,16.) circle (1.cm);
\draw [line width=0.7pt,dashed,color=qqwwzz] (67.,9.) circle (1.cm);
\draw[line width=0.7pt,dashed,color=qqwwzz, smooth,samples=100,domain=0.0:0.3890404923960828] plot[parametric] function{-0.5648396714375444*t**(3.0)+23.896112147974378*t+45.99999999999999,-21.01300050711187*t**(3.0)+6.565362214806641*t+14.999999999999998};
\draw[line width=0.7pt,dashed,color=qqwwzz, smooth,samples=100,domain=0.3890404923960828:0.7559204984484887] plot[parametric] function{-12.904253994854885*t**(3.0)+14.401595472784676*t**(2.0)+18.293308353953034*t+46.72657251560823,11.314676134580461*t**(3.0)-37.730325716115985*t**(2.0)+21.243986709668988*t+13.096473566273852};
\draw[line width=0.7pt,dashed,color=qqwwzz, smooth,samples=100,domain=0.7559204984484887:1.0] plot[parametric] function{20.296904312505433*t**(3.0)-60.8907129375163*t**(2.0)+75.20830765680509*t+32.38550096820579,16.485528648217333*t**(3.0)-49.456585944652005*t**(2.0)+30.108107186560623*t+10.862950109874053};
\draw[line width=0.7pt,dashed,color=qqwwzz, smooth,samples=100,domain=0.0:0.3890404923960828] plot[parametric] function{-0.5648396714375444*t**(3.0)+23.896112147974378*t+45.99999999999999,-21.01300050711201*t**(3.0)+6.565362214806689*t+17.0};
\draw[line width=0.7pt,dashed,color=qqwwzz, smooth,samples=100,domain=0.3890404923960828:0.7559204984484887] plot[parametric] function{-12.904253994854885*t**(3.0)+14.401595472784676*t**(2.0)+18.293308353953034*t+46.72657251560823,11.314676134580711*t**(3.0)-37.73032571611644*t**(2.0)+21.243986709669212*t+15.096473566273831};
\draw[line width=0.7pt,dashed,color=qqwwzz, smooth,samples=100,domain=0.7559204984484887:1.0] plot[parametric] function{20.296904312505433*t**(3.0)-60.8907129375163*t**(2.0)+75.20830765680509*t+32.38550096820579,16.485528648217176*t**(3.0)-49.45658594465153*t**(2.0)+30.10810718656015*t+12.862950109874202};
\draw[line width=0.7pt,dashed,color=qqwwzz, smooth,samples=100,domain=0.0:0.3890404923960828] plot[parametric] function{-0.564839671437503*t**(3.0)+23.89611214797436*t+45.0,-21.013000507111983*t**(3.0)+6.565362214806673*t+16.000000000000004};
\draw[line width=0.7pt,dashed,color=qqwwzz, smooth,samples=100,domain=0.3890404923960828:0.7559204984484887] plot[parametric] function{-12.904253994854995*t**(3.0)+14.401595472784853*t**(2.0)+18.29330835395294*t+45.72657251560824,11.314676134580687*t**(3.0)-37.73032571611638*t**(2.0)+21.243986709669176*t+14.096473566273833};
\draw[line width=0.7pt,dashed,color=qqwwzz, smooth,samples=100,domain=0.7559204984484887:1.0] plot[parametric] function{20.296904312505546*t**(3.0)-60.89071293751664*t**(2.0)+75.20830765680539*t+31.385500968205704,16.485528648217176*t**(3.0)-49.456585944651536*t**(2.0)+30.108107186560165*t+11.8629501098742};
\draw [-to, line width=1pt] (10.,14.) -- (14.,13.);
\draw [-to, line width=1pt] (10.,14.) -- (14.,15.);
\draw [-to, line width=1pt] (10.,14.) -- (10.,18.);
\draw [-to, line width=1pt,color=xfqqff] (42.,8.) -- (49.27761757117354,16.846618269559436);
\draw [-to, line width=1pt,color=xfqqff] (49.27761757117354,16.846618269559436) -- (51.53636252791112,17.26070689448472);
\draw [-to, line width=1pt,color=xfqqff] (49.27761757117354,16.846618269559436) -- (51.,18.);
\draw [-to, line width=1pt,color=xfqqff] (49.27761757117354,16.846618269559436) -- (48.,19.);
\draw [line width=0.7pt,dashed,color=qqwwzz] (49.27761757117354,16.846618269559436) circle (0.985190229365948cm);
\draw (14,12.65) node[anchor=south west] {$s$};
\draw (14,16) node[anchor=north west] {$\eta$};
\draw (10,18) node[anchor=south] {$\zeta$};
\draw [color=xfqqff](45.5,13.197794885098638) node[anchor=north west] {$\mathbf{r}:[0,l] \to \R^3$};
\draw [color=xfqqff](51.5,17.72616728254068) node[anchor=north west] {$\mathbf{t}$};
\draw [color=xfqqff](50.8,19.5) node[anchor=north west] {$\mathbf{n}$};
\draw [color=xfqqff](46.8,20.4) node[anchor=north west] {$\mathbf{c}$};
\draw (46.2,7.6) node[anchor=north west] {$x$};
\draw (46.1,10) node[anchor=north west] {$y$};
\draw (42,12) node[anchor=south] {$z$};
\draw [color=qqwwzz](18.25,12.75) node[anchor=south west] {$\omega$};
\draw [-to, line width=1pt] (32.,16.) -- (38.,16.);
\draw (35,16) node[anchor=south] {$\mathbf{x}:\Omega \to V$};
\draw [color=qqwwzz](24.1,13.5) node[anchor=north west] {$\Omega$};
\draw [color=qqwwzz](63.1561396881974,16.86577652702669) node[anchor=north west] {$V$};
\draw [color=qqwwzz](57.25,17.75) node[anchor=south west] {$\omega$};
\end{tikzpicture}
    \caption{A three-dimensional beam mapped from the reference domain $\Omega \subset \R^3$ to the physical domain $\vol\subset \R^3$ via $\vb{x}: \Omega \subset \R^3 \to \vol \subset \R^3$. The cross-section $\omega \subset \R^2$ of the beam is left unchanged by the mapping. The centroid line of the beam is mapped by $\vb{r}:[0,l] \to \R^3$.}
    \label{fig:beam}
\end{figure}
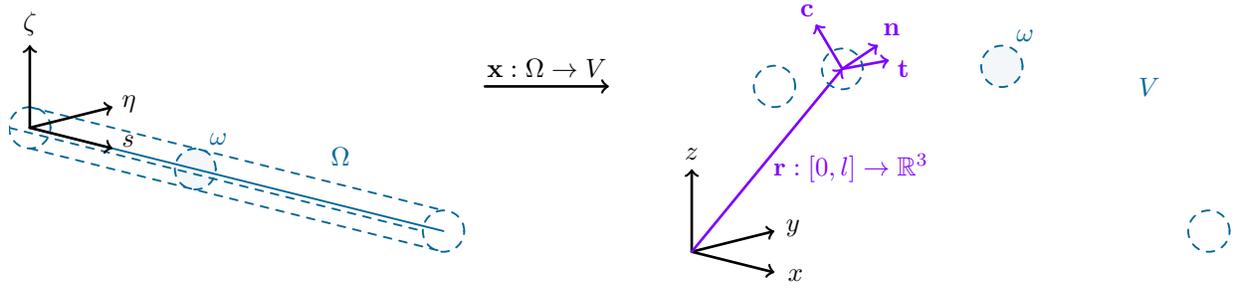

Let $\vb{r} = \vb{r}(s)$ map the centroid curve of some three-dimensional beam, parameterised by the arc-length parameter $s \in [0,l]$ and equipped with the accompanying normal vectors $\vb{n}(s) \perp \vb{t}(s)$ and $\vb{c}(s) = \vb{t} \times \vb{n}$, the beam domain is given by 
\begin{align}
    &\vb{x}(s,\eta,\zeta) = \vb{r} + \eta \vb{n} + \zeta \vb{c} \, , && \vb{x}: [0,l] \times\omega \subset \R^3 \to \vol \subset \R^3 \, , 
\end{align}
where $l \geq 0\setminus\{0\}$ is the length of the beam and $\{\eta,\zeta\}$ define the cross-section of the beam $\omega \subset \R^2$, see \cref{fig:beam}.
The initial orientation of the cross-section of the beam is given by the choice of $\vb{n}$, and \textbf{we exclude twists of its cross-section in the mapping} $\tau = 0$, see \cref{ap:curves}.
The unit tangent vector of the centroid curve is given by
\begin{align}
    &\vb{t} = \vb{r}_{,s} \, , && \norm{\vb{t}} = 1 \, ,
\end{align}
such that we can define the tangential and normal projection operators
\begin{align}
    &\Pt = \vb{t} \otimes \vb{t} \, , && \Pn = \one - \Pt = \one - \vb{t} \otimes \vb{t} = \vb{n} \otimes \vb{n} + \vb{c} \otimes \vb{c} \, , && \Pt : \R^3 \to \tansurf \curv \, , && \Q : \R^3 \to \R^3 \setminus \tansurf \curv  \, ,
\end{align}
where $\tansurf \curv = \spa\{\vb{t}\}$ is the space of tangential vectors to the curve.
Following the same procedure as for the shell models, we employ the decomposition of $\sym\D \ud$ and $\skw(\D\ud-\Rtheta)$ using the projectors. The decomposition of $\sym\D \ud$ yields the same result as in \cref{eq:decompDu}, but the projection operators are now those of the beam. Observe that $\Pt = \vb{t} \otimes \vb{t}$ implies $\Pt [\skw(\D \ud) - \Rtheta] \Pt = \con{\skw \D \ud - \Rtheta}{\Pt} \Pt = 0$ due to skew-symmetry. Thus, \cref{eq:decompSkewDuO} changes to 
\begin{align}
    2 \muc \norm{\skw(\D \ud - \Rtheta)}^2 = 2\muc \norm{ \Pn(\skw \D \ud  - \Rtheta) \Pn }^2 + 4 \muc \norm{\Pn (\skw \D \ud  -  \Rtheta) \Pt }^2  \, ,
\end{align}
and the internal energy of the beam reads
\begin{align}
    I_\curv(\ud,\rtheta) =\dfrac{1}{2} \int_s \int_\omega (&\norm{\Pt(\sym \D \ud)\Pt + \Pn(\sym \D \ud)\Pn}_{\Ce}^2 + 4\mue \norm{\Pn (\sym\D \ud)\Pt}^2 + 4 \muc \norm{\Pn (\skw \D \ud  -  \Rtheta) \Pt }^2 
    \notag \\
    & + 2\muc \norm{ \Pn(\skw \D \ud  - \Rtheta) \Pn }^2 + \mue \Lc^2\norm{(\di \rtheta) \one - (\D \rtheta)^T}_{\Lm}^2) (1- \kappa_n\eta - \kappa_c\zeta)  \,  \dd \omega \dd s  \, ,
\end{align}
where we use \cref{eq:beamvoldecomp} to decompose the integral over the volume of the beam across its cross-section surface $\omega \subset \R^2$ and length $[0,l]\subset \R$.

\subsection{The beam model}
The kinematic of the beam is given via
\begin{align}
    &\vb{u}(s,\xi,\eta) = \vb{v} + \eta \Rtheta \vb{n} + \zeta \Rtheta \vb{c} = \vb{v} + \Rtheta(\eta \vb{n} + \zeta \vb{c}) \, , && \vb{v} = \vb{v}(s) \, , && \Rtheta = \Rtheta(s) \, ,
\end{align}
where $\vb{v}:[0,l] \to \R^3$ is the translation of the centroid curve and $\Rtheta \neq \Rtheta(\eta,\zeta)$ is assumed to be constant throughout the cross-section of the beam. 
In the framework of tangential differential calculus, the gradient of the displacement field can be expressed as
\begin{align}
    \D \ud &= \D_t\vb{v} + \Rtheta_{,s}(\eta \vb{n} + \zeta \vb{c}) \otimes \vb{t} + \Rtheta(\eta \vb{n} + \zeta \vb{c})_{,s} \otimes \vb{t} + \Rtheta\D_n(\eta \vb{n} + \zeta \vb{c}) 
    \notag \\
    &= \D_t \vb{v} - \Anti(\eta \vb{n} + \zeta \vb{c})\D_t \rtheta - (\kappa_n \eta + \kappa_c \zeta)\Rtheta\Pt + \Rtheta\Pn  \, ,
\end{align} 
\textbf{for a twist-free cross-section}, as per \cref{ap:curves}. Consequently, its symmetrised tangential-tangential part reads
\begin{align}
    \Pt (\sym\D \ud)\Pt = \sym(\Dcov \vb{v} - \Pt\Anti(\eta \vb{n} + \zeta \vb{c})\D_t \rtheta) \, ,
\end{align}
since $\Pt \Rtheta \Pt = \con{\Rtheta}{\Pt}\Pt = 0$. Its normal-normal part reads
\begin{align}
    \Pn (\sym\D \ud)\Pn = \Pn (\sym \Rtheta) \Pn = 0 \, .
\end{align}
Analogously, we find
\begin{align}
    \Pt(\skw \D \ud - \Rtheta) \Pt = \skw(\Dcov \vb{v} - \Pt\Anti(\eta \vb{n} + \zeta \vb{c})\D_t \rtheta) \,,
\end{align}
and 
\begin{align}
    \Pn(\skw \D \ud - \Rtheta) \Pn = \Pn \Rtheta \Pn - \Pn \Rtheta \Pn = 0 \, .
\end{align}
Due to $\Pn (\sym\D \ud)\Pn = \Pn(\skw \D \ud - \Rtheta) \Pn = 0$ we get that the perpendicular strains vanish, similarly to the case of shells \cite{MichaelThesis}. Thus, the model is adjusted to the beam equivalent of plane-stress, being
\begin{align}
    \norm{\Pt (\sym \D \ud) \Pt}_{\Ce}^2 = E_e\norm{\Pt (\sym \D \ud) \Pt}^2 \, , && E_e = \dfrac{\mue(3\lame + 2 \mue)}{\lame + \mue} \, ,
\end{align}
via Young's modulus $E_e$, such that the perpendicular strains are eliminated from the energy a priori.
Next, the mixed normal-tangential part of the symmetrised gradient reads
\begin{align}
    \Pn (\sym\D \ud)\Pt = \dfrac{1}{2}\Pn(\D_t \vb{v} - \Anti(\eta \vb{n} + \zeta \vb{c})\D_t \rtheta - (\kappa_n \eta + \kappa_c \zeta)\Rtheta\Pt - \Rtheta\Pt) \, . 
\end{align}
The corresponding skew-symmetric part also yields
\begin{align}
    \Pn(\skw \D \ud - \Rtheta) \Pt = \dfrac{1}{2}\Pn(\D_t \vb{v} - \Anti(\eta \vb{n} + \zeta \vb{c})\D_t \rtheta - (\kappa_n \eta + \kappa_c \zeta)\Rtheta\Pt - \Rtheta\Pt) \, . 
\end{align}
Consequently, the internal energy of the beam is given by
\begin{align}
    I_\curv (\vb{v},\rtheta) =\dfrac{1}{2} \int_s \int_\omega (&E_e\norm{\sym(\Dcov \vb{v} - \Pt[\eta(\Anti \vb{n}) + \zeta (\Anti\vb{c})]\D_t \rtheta)}^2  
    \notag \\
    & + (\mue + \muc) \norm{\Pn(\D_t \vb{v} - [\eta(\Anti \vb{n}) + \zeta (\Anti\vb{c})]\D_t \rtheta - (\kappa_n \eta + \kappa_c \zeta)\Rtheta\Pt - \Rtheta\Pt)}^2
    \\
    \phantom{\dfrac{1}{2}}& \phantom{+ (\mue + \muc)\|} + \mue \Lc^2\norm{(\di_t \rtheta) \one - (\D_t \rtheta)^T}_{\Lm}^2) (1- \kappa_n\eta - \kappa_c\zeta)  \,  \dd \omega \dd s  \, , \notag
\end{align}
where the differential operators of the dislocation density are inherently tangential. Observe that since the parametrisation is with respect to the centroid-line of the beam, integrals of the form $\int_\omega \eta \, \dd \omega = \int_\omega \zeta \, \dd \omega = 0$ vanish. Further, we restrict the formulation to cross-sections that are symmetric with respect to at least the $\eta$- or the $\zeta$-axis. Thus, mixed terms of the form $\int_\omega \eta \zeta \, \dd \omega = 0$ also vanish. By definition we have that
\begin{align}
     &\surf = \int_\omega \,\dd \omega \, , && I_\eta = \int_\omega \zeta^2 \, \dd \omega \, , && I_\zeta = \int_\omega \eta^2 \, \dd \omega \, , && I_p = \int_\omega \zeta^2 + \eta^2 \, \dd \omega = I_\eta + I_\zeta \, , 
\end{align}
represent the cross-section surface and the classical second order moments of inertia, respectively. Lastly, assuming the curvature of the beam is relatively small $(|\kappa_n| + |\kappa_c|)/|l \ll 1$, using asymptotic analysis we omit quadratic terms multiplied by curvatures, cubic and higher order terms. The latter implies that in the integration, the shifter term $1- \kappa_n \eta - \kappa_c \zeta$ is essentially reduced to $1$. 
Consequently, the integration of the first energy term yields
\begin{align}
    E_e\surf\norm{\sym\Dcov \vb{v}}^2 + E_eI_\zeta\norm{\sym( \Pt[\Anti \vb{n} ]\D_t \rtheta)}^2 + E_eI_\eta \norm{\sym(\Pt[\Anti\vb{c}]\D_t \rtheta)}^2 \, ,
\end{align}
and for the second term we find
\begin{align}
    (\mue + \muc)(A\norm{\Pn(\D_t \vb{v} - \Rtheta\Pt)}^2
    + I_\zeta\norm{ \Pn([\Anti \vb{n}]\D_t \rtheta + \kappa_n \Rtheta\Pt )}^2
    + I_\eta\norm{\Pn([\Anti\vb{c}]\D_t \rtheta + \kappa_c\Rtheta\Pt )}^2
    )  \, .
\end{align}
Finally, the dislocation energy is just multiplied by the surface area $A$, such that the internal energy of the beam is given by
\begin{align}
    \boxed{
    \begin{aligned}
    I_\curv(\vb{v},\rtheta) =\dfrac{1}{2} \int_s &E_e\surf\norm{\sym\Dcov \vb{v}}^2 + E_eI_\zeta\norm{\sym( \Pt[\Anti \vb{n} ]\D_t \rtheta)}^2 + E_eI_\eta \norm{\sym(\Pt[\Anti\vb{c}]\D_t \rtheta)}^2  
    \\
    & + (\mue + \muc)A\norm{\Pn(\D_t \vb{v} - \Rtheta\Pt)}^2
    + (\mue + \muc)I_\zeta\norm{ \Pn([\Anti \vb{n}]\D_t \rtheta + \kappa_n \Rtheta\Pt )}^2 
    \\
    \phantom{\dfrac{1}{2}}& \phantom{+ (\mue + \muc)A \|} + (\mue + \muc)I_\eta\norm{\Pn([\Anti\vb{c}]\D_t \rtheta + \kappa_c\Rtheta\Pt )}^2 + \mue \Lc^2 A \norm{(\di_t \rtheta) \one - (\D_t \rtheta)^T}_{\Lm}^2  \,   \dd s  \, . 
\end{aligned}
    }
\end{align}
By their mechanical action, the first term is the membrane energy, the next two terms are the bending energy, and the fourth term is the shear energy density. The two subsequent terms represent warp torsion energy since
\begin{align}
     \Pn(\Anti\vb{n})\D_t \rtheta + \kappa_n\Pn\Rtheta\Pt &= \kappa_n\Pn\Rtheta\Pt  - (\di_t \rtheta) \vb{c} \otimes \vb{t}  \, ,  \notag \\
    \Pn(\Anti\vb{c})\D_t \rtheta + \kappa_c\Pn\Rtheta\Pt &= (\di_t \rtheta) \vb{n} \otimes \vb{t} + \kappa_c\Pn\Rtheta\Pt   \, , 
    \label{eq:warpenergy}
\end{align}
couple rotations perpendicular to the curvature with the rotational intensity. One can observe that $\di_t \rtheta$ relates to the change in torque $\theta^t$ via 
\begin{align}
    &\rtheta = \theta^t \vb{t} + \theta^n \vb{n} + \theta^c \vb{c} \, , && \di_t \rtheta = \con{\D_t \rtheta}{\Pt} = \con{(\theta^t_{,s}  -  \kappa_n \rtheta_n  - \kappa_c \rtheta_c) \vb{t} \otimes \vb{t} }{\Pt} = \theta^t_{,s}  -  \kappa_n \rtheta_n  - \kappa_c \rtheta_c  \, ,    
\end{align} 
as the remaining terms of the gradient $\D_t \rtheta$ are not tangential-tangential and thus eliminated, compare \cref{ap:curves}.
Lastly, the endmost term is the dislocation energy.   
Now, for the minimisation functional we simply reformulate the forces and couple-forces as curve quantities 
\begin{align}
    L_\curv(\vb{v},\rtheta) = \int_\curv \con{\vb{v}}{\vb{q}} + \con{\Rtheta}{\bm{M}} \, \dd \curv \, ,  
\end{align}
such that $\vb{q}:[0,l] \to \R^3$ and $\bm{M}:[0,l] \to \so(3)$. The balance of energy is now expressed as in \cref{eq:energybalanceshell}.

\subsection{The straight beam model}
If the beam is not curved $\kappa_n = \kappa_c = 0$, then the energy functional simplifies to
\begin{align}
    \boxed{
    \begin{aligned}
    I_\curv(\vb{v},\rtheta) =\dfrac{1}{2} \int_s &E_e\surf\norm{\sym\Dcov \vb{v}}^2 + E_eI_\zeta\norm{\sym( \Pt[\Anti \vb{n} ]\D_t \rtheta)}^2 + E_eI_\eta \norm{\sym(\Pt[\Anti\vb{c}]\D_t \rtheta)}^2  
    \\
    & + (\mue + \muc)A\norm{\Pn(\D_t \vb{v} - \Rtheta\Pt)}^2
    + (\mue + \muc)I_p\norm{ \di_t \rtheta  }^2 
    + \mue \Lc^2 A \norm{(\di_t \rtheta) \one - (\D_t \rtheta)^T}_{\Lm}^2  \,   \dd s  \, ,
\end{aligned}
    }
\end{align}
where we used \cref{eq:warpenergy} to obtain $\norm{ \Pn(\Anti \vb{n})\D_t \rtheta }^2 = \norm{ \Pn(\Anti \vb{c})\D_t \rtheta }^2 = \norm{ \di_t \rtheta }^2$, enabling the combination of $I_\eta$ and $I_\zeta$ to the polar second order moment of inertia $I_p = I_\eta + I_\zeta$.
The form can be further simplified if the beam is embedded in the $x$-axis with its cross-section aligned to the $y-z$-plane.

\subsection{A micro-beam model}
A special Cosserat beam model can be derived under the assumption that the surface of the beam is very small $|\surf | / l \ll 1$, implying that the second order moments of inertia $I_\eta \to 0$ and $I_\zeta \to 0$ vanish. In this scenario the energy functional reduces to
\begin{align}
    \boxed{
    \begin{aligned}
    I_\curv(\vb{v},\rtheta) =\dfrac{1}{2} \int_s &E_e\surf\norm{\sym\Dcov \vb{v}}^2 + (\mue + \muc)A\norm{\Pn(\D_t \vb{v} - \Rtheta\Pt)}^2  
     + \mue \Lc^2 A \norm{(\di_t \rtheta) \one - (\D_t \rtheta)^T}_{\Lm}^2  \,   \dd s  \, . 
\end{aligned}
    }
\end{align}
Analogously to the membrane-shell model, we get that the micro-beam model is naturally independent of any explicit curvature energy terms.  

\section{\AS{Discretisation}}
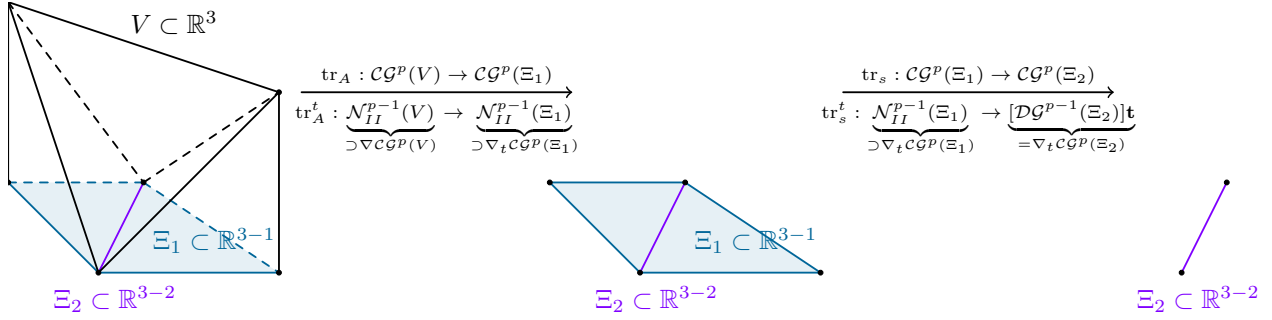
\begin{figure}
    \centering
    \definecolor{xfqqff}{rgb}{0.4980392156862745,0.,1.}
\definecolor{qqwwzz}{rgb}{0.,0.4,0.6}
\begin{tikzpicture}[scale = 0.6, line cap=round,line join=round,>=triangle 45,x=1.0cm,y=1.0cm]
\clip(4,3) rectangle (32,10);
\fill[line width=0.7pt,color=qqwwzz,fill=qqwwzz,fill opacity=0.10000000149011612] (4.,6.) -- (7.,6.) -- (6.,4.) -- cycle;
\fill[line width=0.7pt,color=qqwwzz,fill=qqwwzz,fill opacity=0.10000000149011612] (7.,6.) -- (10.,4.) -- (6.,4.) -- cycle;
\fill[line width=0.7pt,color=qqwwzz,fill=qqwwzz,fill opacity=0.10000000149011612] (16.,6.) -- (19.,6.) -- (18.,4.) -- cycle;
\fill[line width=0.7pt,color=qqwwzz,fill=qqwwzz,fill opacity=0.10000000149011612] (19.,6.) -- (22.,4.) -- (18.,4.) -- cycle;
\draw [line width=0.7pt,dashed,color=qqwwzz] (4.,6.)-- (7.,6.);
\draw [line width=0.7pt,color=qqwwzz] (6.,4.)-- (4.,6.);
\draw [line width=0.7pt,dashed,color=qqwwzz] (7.,6.)-- (10.,4.);
\draw [line width=0.7pt,color=qqwwzz] (10.,4.)-- (6.,4.);
\draw [line width=0.7pt,color=xfqqff] (6.,4.)-- (7.,6.);
\draw [line width=0.7pt] (4.,10.)-- (4.,6.);
\draw [line width=0.7pt,dashed] (4.,10.)-- (7.,6.);
\draw [line width=0.7pt] (6.,4.)-- (4.,10.);
\draw [line width=0.7pt,dashed] (7.,6.)-- (10.,8.);
\draw [line width=0.7pt] (10.,4.)-- (10.,8.);
\draw [line width=0.7pt] (10.,8.)-- (6.,4.);
\draw [line width=0.7pt] (4.,10.)-- (10.,8.);
\draw [line width=0.7pt,color=qqwwzz] (16.,6.)-- (19.,6.);
\draw [line width=0.7pt,color=qqwwzz] (18.,4.)-- (16.,6.);
\draw [line width=0.7pt,color=qqwwzz] (19.,6.)-- (22.,4.);
\draw [line width=0.7pt,color=qqwwzz] (22.,4.)-- (18.,4.);
\draw [line width=0.7pt,color=xfqqff] (18.,4.)-- (19.,6.);
\draw [line width=0.7pt,color=xfqqff] (31.,6.)-- (30.,4.);
\draw [-to,line width=0.7pt] (10.5,8.) -- (16.5,8.);
\draw [-to,line width=0.7pt] (22.5,8.) -- (28.5,8.);
\draw (6.5,10) node[anchor=north west] {$V\subset \R^3$};
\draw [color=qqwwzz](7,5.2) node[anchor=north west] {$\Xi_1 \subset \R^{3-1}$};
\draw [color=xfqqff](4.8,3.9) node[anchor=north west] {$\Xi_2 \subset \R^{3-2}$};
\draw [color=qqwwzz](7+12,5.2) node[anchor=north west] {$\Xi_1 \subset \R^{3-1}$};
\draw [color=xfqqff](4.8+12,3.9) node[anchor=north west] {$\Xi_2 \subset \R^{3-2}$};
\draw [color=xfqqff](4.8+24,3.9) node[anchor=north west] {$\Xi_2 \subset \R^{3-2}$};
\begin{scriptsize}
\draw (13.5,8) node[anchor=south] {$\mathrm{tr}_A: \mathcal{CG}^p(V) \to \mathcal{CG}^p(\Xi_1)$};
\draw (13.5,8) node[anchor=north] {$\mathrm{tr}^t_A: \underbrace{\mathcal{N}_{II}^{p-1}(V)}_{\supset \nabla \mathcal{CG}^{p}(V)} \to \underbrace{\mathcal{N}_{II}^{p-1}(\Xi_1)}_{\supset \nabla_t \CG^{p}(\Xi_1)}$};
\draw (25.5,8) node[anchor=south] {$\mathrm{tr}_s: \mathcal{CG}^p(\Xi_1) \to \mathcal{CG}^p(\Xi_2)$};
\draw (25.5,8) node[anchor=north] {$\mathrm{tr}^t_s: \underbrace{\mathcal{N}_{II}^{p-1}(\Xi_1)}_{\supset \nabla_t \CG^p(\Xi_1)} \to  \underbrace{[\mathcal{DG}^{p-1}(\Xi_2)]\vb{t}}_{=\nabla_t\CG^p(\Xi_2)}$};
\draw [fill=black] (4.,6.) circle (1.5pt);
\draw [fill=black] (7.,6.) circle (1.5pt);
\draw [fill=black] (6.,4.) circle (1.5pt);
\draw [fill=black] (10.,4.) circle (1.5pt);
\draw [fill=black] (4.,10.) circle (1.5pt);
\draw [fill=black] (10.,8.) circle (1.5pt);
\draw [fill=black] (16.,6.) circle (1.5pt);
\draw [fill=black] (19.,6.) circle (1.5pt);
\draw [fill=black] (22.,4.) circle (1.5pt);
\draw [fill=black] (18.,4.) circle (1.5pt);
\draw [fill=black] (31.,6.) circle (1.5pt);
\draw [fill=black] (30.,4.) circle (1.5pt);
\end{scriptsize}
\end{tikzpicture}
    \caption{Trace operators from a meshed three-dimensional domain $\vol \subset \R^3$ of three tetrahedral elements to domains of codimensions one $\Xi_1 \subset \R^{3-1}$ of two triangles and two $\Xi_2 \subset \R^{3-2}$ of a single line. The operators yield valid finite element spaces on the lower dimensional entities. Namely, restricting the continuous three-dimensional Lagrange space $\CG^p(\vol)$ onto $\Xi_1$ or $\Xi_2$ via $\tr_\surf$ or $\tr_s$ yields the equivalent Lagrange space on these lower dimensional domains, being $\CG^p(\Xi_1)$ and $\CG^p(\Xi_2)$. Analogously, the tangential trace operator $\tr_\surf^t$ yields a valid surface N\'ed\'elec finite element space $\Nedtwo^{p-1}(\Xi_1) \supset \nabla_t\CG^p(\Xi_1)$ from the volumetric N\'ed\'elec space $\Nedtwo^{p-1}(\vol) \supset \nabla \CG^p{\vol}$, and the application of the second tangential trace operator $\tr_s^t$ onto a curve is a valid discontinuous Lagrange space on it $[\DG^{p-1}(\Xi_2)]\vb{t} = \nabla_t \CG^p(\Xi_2)$.}
    \label{fig:hier} 
\end{figure}

The volumetric geometry of the domain is discretised using finite element meshing procedures. Correspondingly, the geometry of any embedded shell or beam in the domain is explicitly controlled by the meshing procedure as well. In other words, so-called conforming meshes are required, where the geometries of the lower dimensional shell and beams are clearly identified as faces and lines of three-dimensional polyhedra, compare \cref{fig:hier}.     

The displacement and infinitesimal rotation fields of the three-dimensional Cosserat model are defined using $\C^0(\vol)$-continuous Lagrange elements $\ud,\rtheta \in [\CG^p(\vol)]^3 \subset [\Hone(\vol)]^3$ over the volumetric domain $\vol \subset \R^3$. 
The coupling with the lower dimensional models is now achieved by restriction of the fields to codimensional domains using consistent Sobolev trace operators \cite{Hiptmair}, as depicted in \cref{fig:hier}. 
For the middle surface $\surf \subset \R^3$ of a shell we have
\begin{align}
    \tr_\surf \ud &= \ud \at_\surf \, , & \tr_\surf \rtheta &= \rtheta \at_\surf \, ,  
    \notag \\
    \tr_\surf^t (\D \ud) &= (\D \ud) \at_\surf \Pt = \D_t (\tr_\surf  \ud) = \D_t \ud \at_\surf \, ,
    & \tr_\surf^t (\D \rtheta) &= (\D \rtheta) \at_\surf \Pt = \D_t (\tr_\surf  \rtheta) = \D_t \rtheta \at_\surf \, ,
\end{align}
where $\Pt = \one - \vb{n} \otimes \vb{n}$, and analogously for the centroid line $\curv  \subset \R^3$ of beams
\begin{align}
    \tr_\curv \ud &= \ud \at_\curv \, , & \tr_\curv \rtheta &= \rtheta \at_\curv \, , 
    \notag \\
    \tr_\curv^t (\D \ud) &= (\D \ud) \at_\curv \Pt = \D_t (\tr_\curv  \ud) = \D_t \ud \at_\curv \, , 
    & \tr_\curv^t (\D \rtheta) & = (\D \rtheta) \at_\curv \Pt= \D_t (\tr_\curv  \rtheta) = \D_t \rtheta \at_\curv \, ,
\end{align}
where for curves we have $\Pt = \vb{t} \otimes \vb{t}$.
This approach is naturally consistent \textbf{on the finite element spaces since scalar products and norms remain well-defined and square integrable}. There hold the relations
\begin{align}
     &\tr_s u = \tr_s \tr_\surf u  \, , && \nabla_t\tr_\surf u = \tr_\surf^t \nabla u  \, , && \nabla_t \tr_s u = \nabla_t \tr_s \tr_\surf u = \tr_s^t \nabla_t \tr_\surf u = \tr_s^t \nabla u  \, ,
\end{align}
where $u:\vol \subset \R^3 \to \R$ represents one row of $\ud:\vol\subset \R^3 \to \R^3$.
In other words, one finds the commuting 
de Rham diagram \cite{PaulyDeRham}  
\begin{align}
    \boxed{
    \begin{matrix}
        &\CG^p(\vol) \cap \Hone(\vol) &\xrightarrow[]{\nabla} &\Nedtwo^{p-1}(\vol) \cap \Hc{,\vol} & \supset & \nabla \CG^p(\vol) 
        \\[0.5em]
        &\tr_\surf \bigg\downarrow & &\tr_\surf^t \bigg\downarrow
        \\[1em]
        &\CG^p(\surf) \cap \Hone(\surf) &\xrightarrow[]{\nabla_t} &\Nedtwo^{p-1}(\surf) \cap \Hct{,\surf} & \supset & \nabla_t \CG^p(\surf)
        \\[0.5em]
        &\tr_\curv \bigg\downarrow & &\tr_\curv^t \bigg\downarrow
        \\[1em]
        &\CG^p(\curv) \cap \Hone(\curv) &\xrightarrow[]{\nabla_t} & [\DG^{p-1}(\curv) \cap \Le(\curv)] \vb{t} & \supset & \nabla_t \CG^p(\curv)
    \end{matrix} \, ,
    }
\end{align}
where $\Nedtwo^{p-1}(\vol) \supset \nabla \CG^p(\vol)$, $\Nedtwo^{p-1}(\surf)\supset \nabla_t \CG^p(\surf)$ are the respective volume and surface N\'ed\'elec elements of the second type \cite{Ned2,sky_polytopal_2022,sky_higher_2023,sky_hybrid_2021}, and $[\DG^{p-1}(\curv)] \vb{t} = \nabla_t \CG^p(\curv)$ are discontinuous Lagrange elements on the curve.
Thus, the total mixed-dimensional energy functional can be naturally defined using the single discretisation $\{\ud,\rtheta\} \in [\CG^p(\vol)]^3 \times [\CG^p(\vol)]^3$ as
\begin{align}
   \boxed{
    \begin{aligned}
        I(\ud,\rtheta) = &I_\vol(\ud, \rtheta) + I_\surf(\tr_\surf \ud, \tr_\surf \rtheta) + I_\curv(\tr_\curv \ud, \tr_\curv \rtheta) \,,  &&  \ud, \rtheta \in [\CG^p(\vol)]^3 \, ,
    \end{aligned}
   }
\end{align}
where $I_\vol(\ud,\rtheta)$ is the energy of the volumetric Cosserat model, $I_\surf(\tr_\surf \ud, \tr_\surf \rtheta)$ is the energy of an embedded shell, and correspondingly $I_\curv(\tr_\curv \ud, \tr_\curv \rtheta)$ is the energy of an embedded beam. To clarify, the displacement field of the shell is simply set to $\vb{v} = \tr_\surf \ud$, and analogously for the beam we have $\vb{v} = \tr_s \ud$. 
Accordingly, the total external work is given by
\begin{align}
    \boxed{
    \begin{aligned}
        L(\ud,\rtheta) = &L_\vol(\ud, \rtheta) + L_\surf(\tr_\surf \ud, \tr_\surf \rtheta) + L_\curv(\tr_\curv \ud, \tr_\curv \rtheta) \,,  &&  \ud, \rtheta \in [\CG^p(\vol)]^3 \, ,
    \end{aligned}
    }
\end{align}
such that $L_\vol(\ud,\rtheta)$ is the volumetric work, $L_\surf(\tr_\surf \ud, \tr_\surf \rtheta)$ is work on an embedded shell, and $L_\curv(\tr_\curv \ud, \tr_\curv \rtheta)$ is the work of an embedded beam. Finally, the discrete problem reads 
\begin{align}
   \boxed{
   \begin{aligned}
       I(\ud,\rtheta) - L(\ud,\rtheta) \to  \min \quad \wrt \quad \{\ud,\rtheta\} \, , && \{\widetilde{\ud},\widetilde{\rtheta}\} = \argmin_{\ud,\rtheta} (I - L) \, , && \ud, \rtheta \in [\CG^p(\vol)]^3 \, .
   \end{aligned}
   }
\end{align}
\AS{We remark that, since traces of the kinematical bulk-fields are used to define the corresponding fields on codimensional domains, the Dirichlet boundary conditions of the fields on the codimensional domains are explicitly embedded in the finite element spaces of the bulk-fields.
Further, although the $\C^0$-regularity of the displacement and rotation fields, $\ud$ and $\rtheta$, as well as the tangential regularity of their gradients, $\D \ud$ and $\D \rtheta$, are guaranteed by construction, jumps in stiffness or material coefficients naturally lead to jumps in the resulting stress field $\jump{\bm{\sigma}}\neq 0$, see \cref{eq:cosseratstress}. This is not specific to the proposed method or the Cosserat model, and appears also in classical models of linear elasticity, compare \cite{SkyFormulae} for an example with the Reissner--Mindlin plate model.}

\section{\AS{Numerical examples}}
In the following we compute three examples using NGSolve with \textbf{cubic} order Lagrange polynomials $\CG^3(\vol)$, to mitigate potential locking effects \cite{NS21,sky2023reissnermindlin}. In the first example we consider the 3D-2D-coupling of a volume with shells, while the second example demonstrates the 3D-2D-1D coupling of a volume with a plate and beams. The last example showcases the 2D-1D coupling of shells and beams with intersections. \AS{The examples are available as open-source files under the link in Fn.~\ref{fn:github}}. 

\subsection{Reinforcement with a stiff shell}
\begin{figure}
    \centering
    \definecolor{asl}{rgb}{0.4980392156862745,0.,1.}
\definecolor{asb}{rgb}{0.,0.4,0.6}
\newcommand\x{1}
\begin{tikzpicture}[scale = 0.8, line cap=round,line join=round,>=triangle 45,x=1.0cm,y=1.0cm]
\clip(2.5,2.5) rectangle (13,8.5);

\draw [line width=0.7pt, dashed] (3.,7.5)-- (5.5,8.);
\draw [line width=0.7pt, dashed] (5.5,8.)-- (5.5,7.5);
\draw [line width=0.7pt, dashed] (5.5,7.5)-- (3.,7.);
\draw [line width=0.7pt, dashed] (3.,7.)-- (3.,7.5);
\draw [line width=0.7pt, dashed] (10.,3.5)-- (10.,3.);
\draw [line width=0.7pt, dashed] (10.,3.)-- (12.5,3.5);
\draw [line width=0.7pt, dashed] (12.5,3.5)-- (12.5,4.);
\draw [line width=0.7pt, dashed] (12.5,4.)-- (10.,3.5);

\draw [line width=0.7pt, dashed] plot [smooth] coordinates {(3,7) (4,6.8-\x/4) (6-\x,4.5-\x/4) (9-0.75*\x,3.2-\x/4) (10, 3)};
\draw [line width=0.7pt, dashed] plot [smooth] coordinates {(3,7+0.5) (4,6.8+0.5-\x/4) (6-\x,4.5+0.5-\x/4) (9-0.75*\x,3.2+0.5-\x/4) (10, 3+0.5)};
\begin{scope}
    \clip (3,10) rectangle (12.5, 2);
    \fill [opacity=0.1, color=asb] plot [smooth cycle] coordinates {(3,7) (4,6.8-\x/4) (6-\x,4.5-\x/4) (9-0.75*\x,3.2-\x/4) (10, 3) (10+2.5, 3+0.5) (9+2.5 + 0.75 * \x ,3.2+0.5+\x/4) (6+2.5 + \x ,4.5+0.5 + \x / 4) (4+2.5+0.4,6.8+0.5 + \x / 8) (3+2.5,7+0.5)};

    \fill[white!100]  (3,7-1.22) circle (1.2);
\end{scope}

\draw [line width=0.7pt, dashed] plot [smooth] coordinates {(3+2.5,7+0.5) (4+2.5+0.4,6.8+0.5 + \x / 8) (6+2.5 + \x ,4.5+0.5 + \x / 4) (9+2.5 + 0.75 * \x ,3.2+0.5+\x/4) (10+2.5, 3+0.5)};
\draw [line width=0.7pt, dashed] plot [smooth] coordinates {(3+2.5,7+0.5+0.5) (4+2.5+0.4,6.8+0.5 + \x / 8 + 0.5) (6+2.5 + \x ,4.5+0.5 + \x / 4+0.5) (9+2.5 + 0.75 * \x ,3.2+0.5+\x/4+0.5) (10+2.5, 3+0.5+0.5)};

\fill[fill=black,fill opacity=0.1] (3.,7.5) -- (5.5,8.) -- (5.5,7.5) -- (3.,7.) -- cycle;
\fill[fill=black,fill opacity=0.1] (10.,3.5) -- (10.,3.) -- (12.5,3.5) -- (12.5,4.) -- cycle;

\draw [-to, color=asl, line width=1pt] (7,6.5)-- (7,5.5);
\draw (7,6.1) node[color=asl, anchor=east] {$\vb{f}$};

\draw (8,4.75) node[color=asb , anchor=north west] {$\Xi$};
\draw (10,6) node[anchor=north west] {$V$};

\draw (11.5,3.35) node[anchor=north west] {$\surf_D^{2}$};
\draw (3.5,8.5) node[anchor=north west] {$\surf_D^{1}$};

\end{tikzpicture}
    \caption{Illustration of the volumetric domain $\vol$ with Dirichlet boundary conditions on $\surf_D^1$ and $\surf_D^2$ and body force $\vb{f}$. The domain is subsequently reinforced with a stiff shell defined on its lower surface $\Xi$.}
    \label{fig:shelldom}
\end{figure}
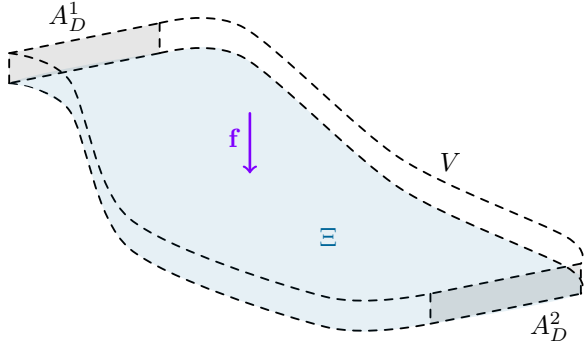

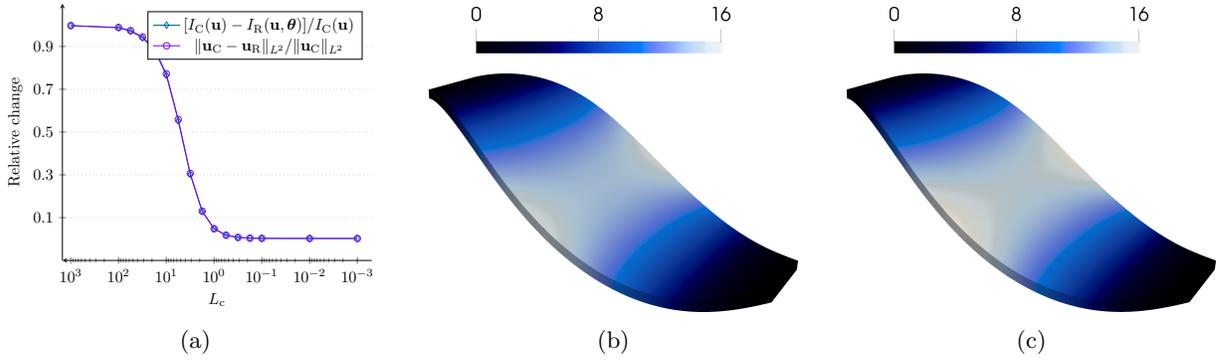
\begin{figure}
    	\centering
    	\begin{subfigure}{0.32\linewidth}
    		\centering
    		\definecolor{asb}{rgb}{0.,0.4,0.6}
\definecolor{asl}{rgb}{0.4980392156862745,0.,1.}
\begin{tikzpicture}[scale = 0.6]
			\begin{semilogxaxis}[
				/pgf/number format/1000 sep={},
				axis lines = left,
				xlabel={$\Lc$},
				ylabel={Relative change},
				xmin=0.5e-3, xmax=1.5e+3,
				ymin=-0.1, ymax=1.1,
				x dir=reverse,
				xtick={1e-3, 1e-2, 1e-1, 1, 1e+1, 1e+2, 1e+3},
				ytick={0.1, 0.3, 0.5, 0.7, 0.9},
				legend pos=north east,
				ymajorgrids=true,
				grid style=dotted,
				]
				\addplot[
				color=asb,
				mark=diamond,
				]
				coordinates {
                        ( 1000 , 0.9970902028979953 )
					( 100 , 0.9875014337770558 )
                        ( 56.23413251903491 , 0.9716393653822122 )
                        ( 31.622776601683793 , 0.9403492567136116 )
                        ( 17.78279410038923 , 0.8830282293874189 )
                        ( 10 , 0.7650818940376566 )
                        ( 5.623413251903491 , 0.5521831135995838 )
                        ( 3.1622776601683795 , 0.30274462642156835 )
                        ( 1.7782794100389228 , 0.1288385155768815 )
                        ( 1 , 0.047955994706616475 )
                        ( 0.5623413251903491 , 0.017862753228809283 )
                        ( 0.31622776601683794 , 0.0077153690924958735 )
                        ( 0.1778279410038923 , 0.0044291116444866285 )
                        ( 0.1 , 0.00338136512267824 )
                        ( 0.01 , 0.0029002299099233045 )
                        ( 0.001 , 0.0028954138440842026 )
				};
				\addlegendentry{$[I_\mathrm{C}(\ud) - I_\mathrm{R}(\ud,\rtheta)] / I_\mathrm{C}(\ud) $}
				\addplot[
				color=asl,
				mark=o,
				]
				coordinates {
                        ( 1000 , 0.9975154331389395 )
					( 100 , 0.9887816090672218 )
                        ( 56.23413251903491 , 0.9742079592334322 )
                        ( 31.622776601683793 , 0.9448587662531633 )
                        ( 17.78279410038923 , 0.8892283312911997 )
                        ( 10 , 0.7721095946673867 )
                        ( 5.623413251903491 , 0.5586770456030198 )
                        ( 3.1622776601683795 , 0.3068873726179885 )
                        ( 1.7782794100389228 , 0.13044510041089544 )
                        ( 1 , 0.048235438597977784 )
                        ( 0.5623413251903491 , 0.017663938535142934 )
                        ( 0.31622776601683794 , 0.007370787558694929 )
                        ( 0.1778279410038923 , 0.0040507298684588555 )
                        ( 0.1 , 0.0030000634487190618 )
                        ( 0.01 , 0.0025213026333444064 )
                        ( 0.001 , 0.0025165294192718256 )
				};
				\addlegendentry{$\norm{\ud_\mathrm{C} - \ud_\mathrm{R}}_{\Le} / \norm{\ud_\mathrm{C}}_{\Le} $}
			\end{semilogxaxis}
		\end{tikzpicture}
    		\caption{}
    	\end{subfigure}
    	\begin{subfigure}{0.32\linewidth}
    		\centering
    		\includegraphics[width=0.9\linewidth]{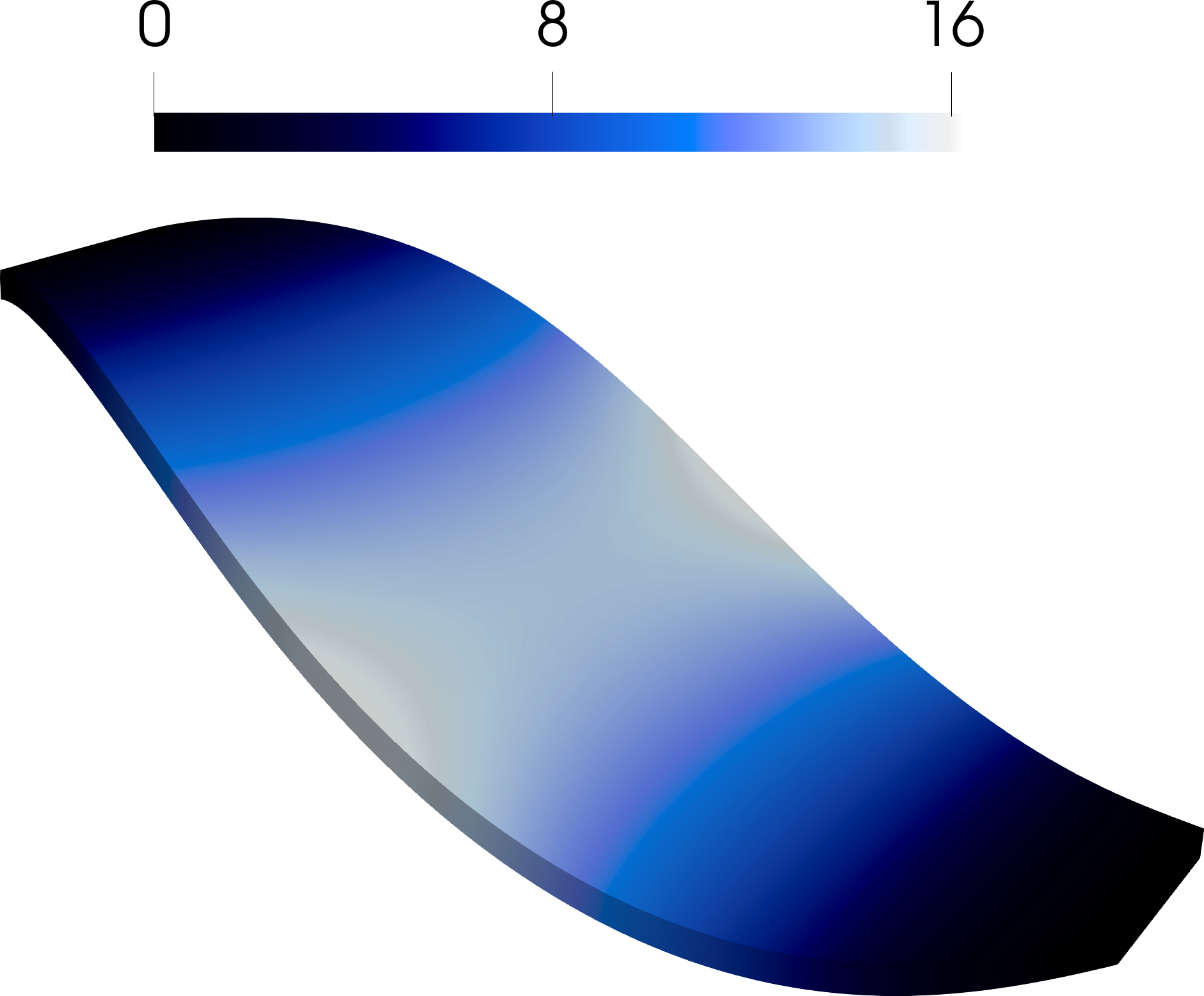}
    		\caption{}
    	\end{subfigure}
            \begin{subfigure}{0.32\linewidth}
    		\centering
    		\includegraphics[width=0.9\linewidth]{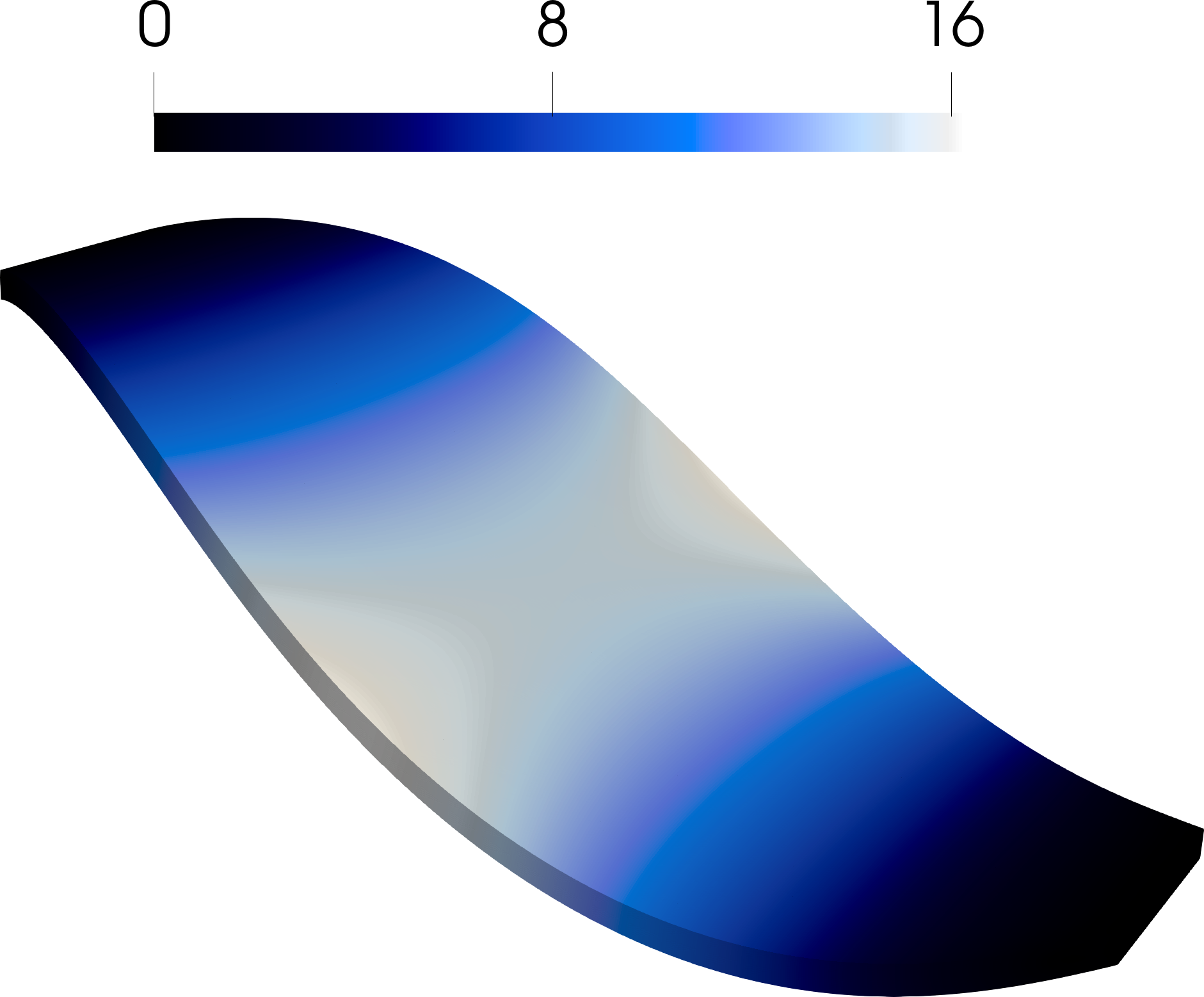}
    		\caption{}
    	\end{subfigure}
    	\caption{Convergence of the Cosserat model for $\Lc \to 0$ towards the equivalent Navier--Cauchy model of linear elasticity (a), and the corresponding displacement field $\ud$ for $\Lc = 1 \, \mm $ (b), and $\Lc = 10^{-2} \, \mm$ (c).}
    \label{fig:energylc}
\end{figure}

In the following we consider a three-dimensional domain made of silicone rubber and reinforce it with shells made of graphite. The domain is given by 
\begin{align}
    \overline{\vol} = \{ \bm{\xi} \in [0,1]^3 \; | \; \vb{x} = \begin{bmatrix}
        200(2\xi-1) & 40(2\eta-1)(3-2[2\xi-1]^2) & 10(\zeta - 7\sin[4\xi-2])
    \end{bmatrix}^T \} \, ,
\end{align}
such that its length is $400 \, \mm$, its minimal width is $80 \, \mm$, its maximal width is $240 \, \mm$ and its thickness is $10 \, \mm$, see \cref{fig:shelldom}.
The three-dimensional body is made of silicone rubber, whose Lam\'e parameters read 
\begin{align}
    &\lama \approx 5.328 \, \MPa \, , && \mua \approx 0.34 \, \MPa \, ,
\end{align}
such that its Poisson ratio is $\nu \approx 0.47$. Consequently, this type of rubber is still compressible and volumetric locking does not occur in numerical simulations. The volume force acting on the domain is 
\begin{align}
    &\vb{f} = -f \vb{e}_3 \, , && f = 10^{-6} \, \Nwtn / \mm^3 \, ,
\end{align}
pointing downwards. We start by fine-tuning the characteristic length-scale parameter $\Lc$ of the Cosserat model to capture an equivalent solution to that of the Navier--Cauchy model. We set 
\begin{align}
    &\muc = 1 \, , && \Lm (\Curl \Rtheta) = \id(\Curl \Rtheta) = \Curl \Rtheta \, ,
\end{align}
such that $\Lm$ is the identity map. The Dirichlet boundary of the displacements is given by 
\begin{align}
&\ud \at_{\surf_D^\vb{u}} = 0 \, , && \surf_D^\vb{u} = \surf_D^1 \cup \surf_D^2 \, ,
\end{align} 
where the boundary surfaces read
\begin{align}
    \surf_D^1 &= \{ (\eta,\zeta) \in [0,1]^2 \; | \; \vb{x} = \begin{bmatrix}
        -200 & 40(2\eta-1) & 10(\zeta - 7\sin[-2])
    \end{bmatrix}^T \} \, , \\
    \surf_D^2 &= \{ (\eta,\zeta) \in [0,1]^2 \; | \; \vb{x} = \begin{bmatrix}
        200 & 40(2\eta-1) & 10(\zeta - 7\sin[2])
    \end{bmatrix}^T \} \, ,
\end{align}
such that the displacement vanishes on $x = \pm 400$. For the rotations the entire boundary is homogeneous Neumann $\surf_N^{\Rtheta} = \partial \vol$.  
The resulting convergence in relative energy and displacements is given in \cref{fig:energylc}, where we compare with a classical Navier--Cauchy formulation of the same domain. 
We observe that for $\Lc \leq 1/\sqrt{10} \, \mm$ the convergence flattens towards a relative difference of $ < 1\%$ in both the energy and the displacements. Still, even for $\Lc = 1 \, \mm$ we observe $\max \norm{\ud} \approx 15.46 \, \mm$ in comparison to $\max \norm{\ud} \approx 16.28 \, \mm$ for $\Lc = 10^{-2} \, \mm$, which is equivalent to the result of the Navier--Cauchy model. We note that it is not possible to exactly capture the Navier--Cauchy model with the Cosserat model numerically, since $\Lc \to 0$ implies $\Rtheta = \skw \D \ud$, but these fields belong to different discrete spaces. Namely, $\skw \D \ud$ is discontinuous while $\Rtheta$ is continuous. Still, the deviation for $\Lc \leq 10^{-2} \, \mm$ is insignificantly small, such that we henceforth use $\Lc = 10^{-2} \, \mm$ in our computations. 
We remark that although an alternative computation with $\muc \to 0$ via $\muc = 10^{-7} \, \MPa \approx 0$ and $\Lc = 1 \, \mm$ leads to an equivalent result with $\max \norm{\ud} \approx 16.31 \, \mm$, \textbf{at the limit $\lim\muc \to 0$ this approach amounts to solving two independent problems}. Namely, $-\Di (\Cm_{\mathrm{M}} \sym \D \ud) = \vb{f}$ and $\mue \Lc^2 \Curl (\Lm \Curl \Rtheta) = 0$, resulting in $\skw \D \ud \neq \Rtheta \to 0$ evidenced by $\norm{\skw \D \ud}_{\Le} \approx 96.42 \neq \norm{\Rtheta}_{\Le} \approx 1.59$ for $\muc = 10^{-7} \, \MPa$, such that the coupling of the infinitesimal macro-rotation $\skw \D \ud$ and the infinitesimal micro-rotation $\Rtheta$ is lost. For $\muc = 0$ with no Dirichlet boundary $|\surf_D^{\Rtheta}| = 0$ the problem is not well-posed. 

\begin{figure}
    	\centering
    	\begin{subfigure}{0.32\linewidth}
    		\centering
                \includegraphics[width=0.9\linewidth]{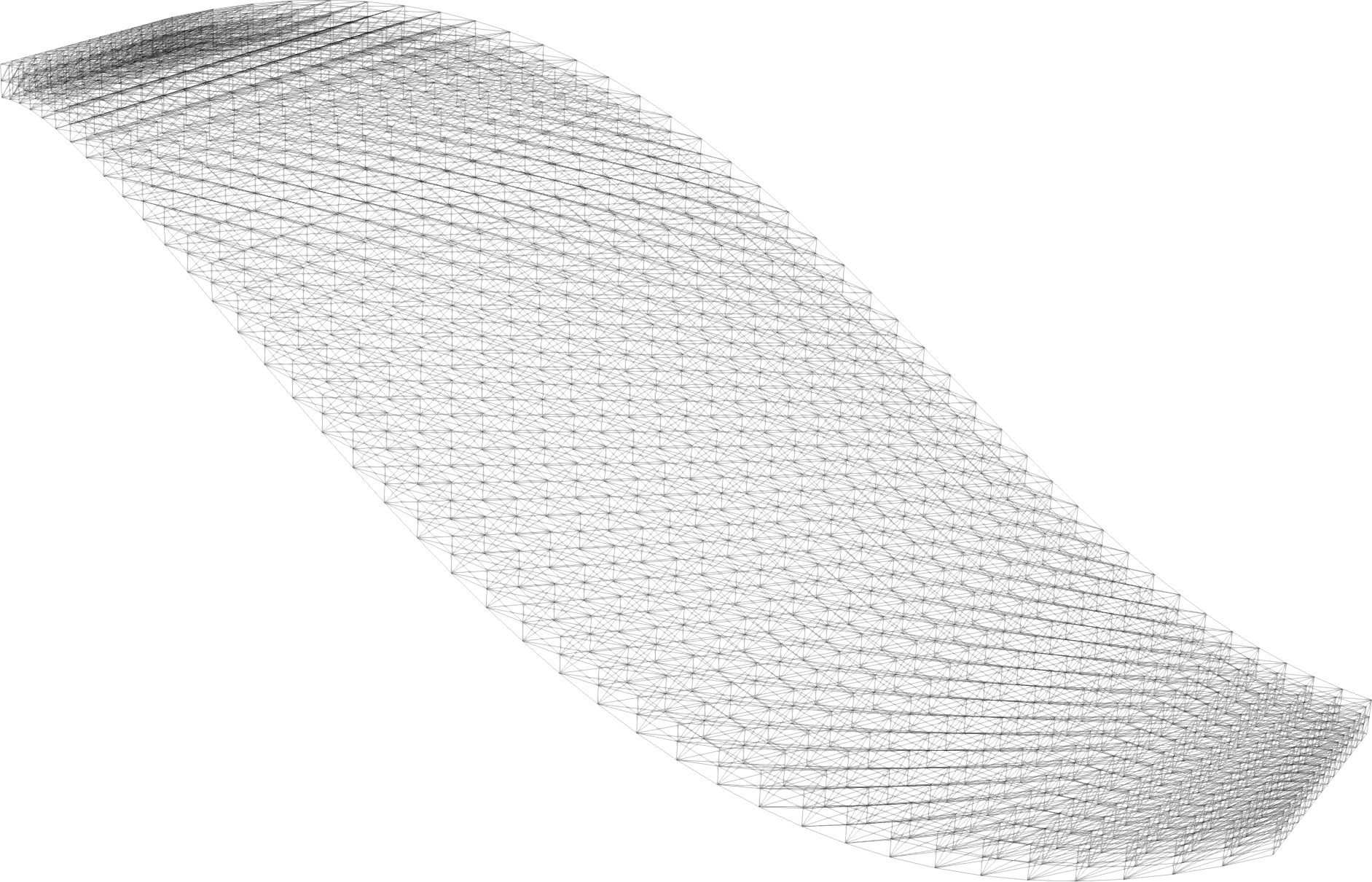}
    		\caption{}
    	\end{subfigure}
    	\begin{subfigure}{0.32\linewidth}
    		\centering
    		\includegraphics[width=0.9\linewidth]{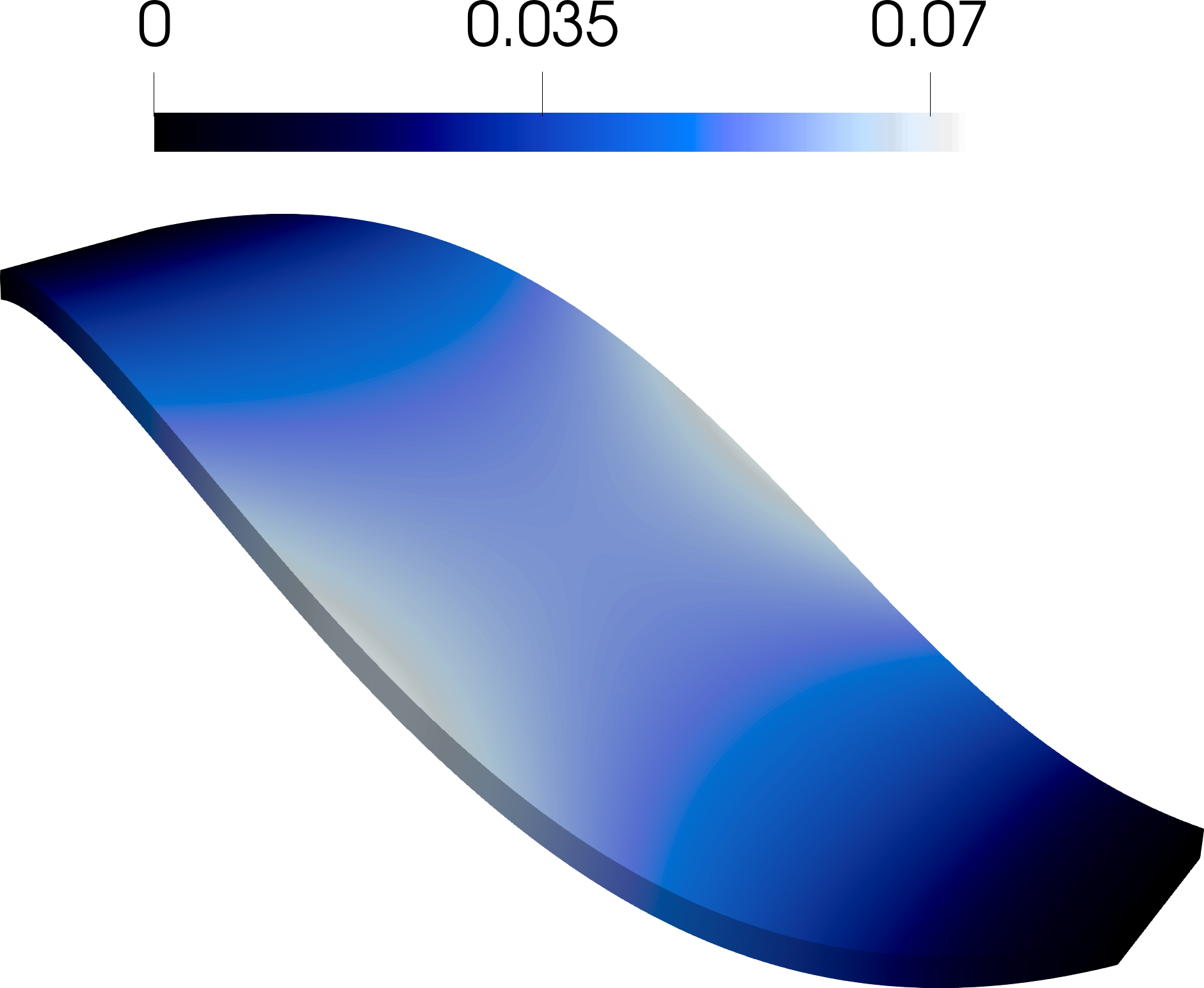}
    		\caption{}
    	\end{subfigure}
            \begin{subfigure}{0.32\linewidth}
    		\centering
    		\includegraphics[width=0.9\linewidth]{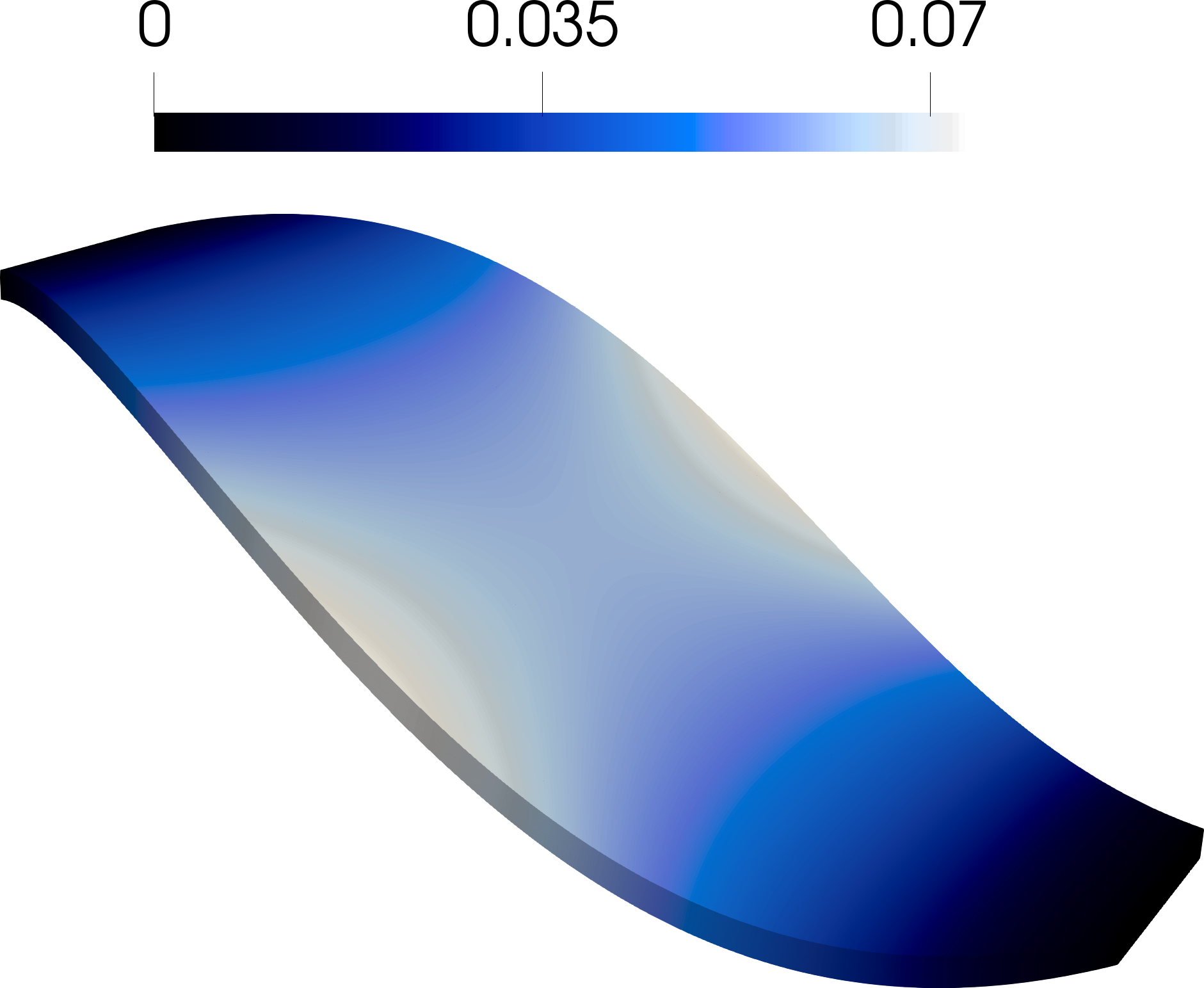}
    		\caption{}
    	\end{subfigure}
    	\caption{Finite element mesh with $7680$ cubic 
Lagrange elements, where the bottom surface $\Xi \subset \partial \vol$ is implicitly reinforced with shells via the energy $I_A(\tr_\surf \ud, \tr_\surf \rtheta)$ (a). Resulting displacement field of the shell (b) and  membrane-shell (c) formulations in $\mm$.}
    \label{fig:reinforced}
\end{figure}

The design in \cref{fig:energylc} is made of extremely soft silicone rubber, leading to deformations of $\approx 16 \, \mm$ even for the very small volume force of $10^{-6} \, \Nwtn / \mm^3$. The deformation can be substantially reduced by adding a thin reinforcement layer at the bottom surface of the shell
\begin{align}
    \Xi = \{ (\xi,\eta) \in [0,1]^2 \; | \; \vb{r} = \begin{bmatrix}
        200(2\xi-1) & 40(2\eta-1)(3-2[2\xi-1]^2) & - 70\sin(4\xi-2)
    \end{bmatrix}^T \} \, ,
\end{align}
which is simply the parametrisation of the domain evaluated with $\zeta = 0$.
We do so using the presented shell models composed of graphite. The material coefficients of the considered graphite @$1.6\, \mm$
(H237) material read
\begin{align}
    \lame &= 289.451 \, \MPa \, , & \mue &= 2122.64 \, \MPa \, , & \muc &\to + \infty \, , \notag \\  
    \mue\Lc a_1 &= 10867.9 \, \Nwtn \, , & \mue \Lc a_2 &= 122264 \, \Nwtn \, , &  \mue \Lc a_3 &= 0 \, ,
\end{align}
stemming from the experimental results in \cite{Lakes1995EXPERIMENTALMF}. Conversion formulae of the material coefficients between different forms of the Cosserat model are available in \cite{Ghiba}.
We note that the characterisation does not allow to determine $\Lc$, $a_1$, and $a_2$ individually, but can still be used in the simulation. For the infinitely large Cosserat couple modulus $\muc \to + \infty$ \cite{MADEO2016294} we take a value one order of magnitude higher than $\mue$ for simplicity $\muc = 10000 \, \MPa$, assuming it is sufficient to enforce the implied constraint
\begin{align}
    \muc \to +\infty : \qquad \Rtheta = \skw \D \ud \, , && \Curl \Rtheta = \Curl (\skw \D \ud) = - \dfrac{1}{2} \Curl (\D \ud)^T = -\dfrac{1}{2} (\D \curl \ud)^T \, ,
\end{align}
which must be satisfied for finite energies. We compare the result of the shell and membrane-shell formulations in \cref{fig:reinforced} setting the shell-thickness to $h = 1.6 \, \mm$, representing one layer of the graphite material.
The maximal displacement using the shell formulation from \cref{sec:shell} is $\max \norm{\ud} \approx 0.068 \, \mm$. In comparison, the membrane-shell formulation from \cref{sec:membrane-shell} yields $\max \norm{\ud} \approx 0.073 \, \mm$. Clearly, for both formulations we observe a significant reduction in the deformation from $\max\norm{\ud} \approx 16.28 \, \mm$, corroborating our shell-volume coupling approach. The relative difference between the two methods given by $100(0.073 - 0.068)/(0.073) \approx 6.8 \, \%$ suggests that for this small thickness the higher order energy terms relating to bending and curvature are negligible.

\subsection{A cantilever beam}
\begin{figure}
    	\centering
    	\begin{subfigure}{0.64\linewidth}
    		\centering
            \definecolor{uuuuuu}{rgb}{0.26666666666666666,0.26666666666666666,0.26666666666666666}
\definecolor{xfqqff}{rgb}{0.4980392156862745,0.,1.}
\definecolor{qqwwzz}{rgb}{0.,0.4,0.6}
\begin{tikzpicture}[scale = 0.7, line cap=round,line join=round,>=triangle 45,x=1.0cm,y=1.0cm]
\clip(2.5,3.5) rectangle (16.5,8.5);
\fill[line width=0.7pt,color=qqwwzz,fill=qqwwzz,fill opacity=0.10000000149011612] (3.,6.5) -- (14.,4.5) -- (16.,5.5) -- (5.,7.5) -- cycle;
\fill[line width=0.7pt,fill=black,fill opacity=0.05] (3.,7.) -- (5.,8.) -- (5.,6.5) -- (3.,5.5) -- cycle;
\draw [line width=0.7pt,color=xfqqff] (3.,6.5)-- (14.,4.5);
\draw [line width=0.7pt,dashed,color=qqwwzz] (14.,4.5)-- (16.,5.5);
\draw [line width=0.7pt,color=xfqqff] (16.,5.5)-- (5.,7.5);
\draw [line width=0.7pt,dashed,color=qqwwzz] (5.,7.5)-- (3.,6.5);
\draw [line width=0.7pt,dashed] (3.,5.5)-- (5.,6.5);
\draw [line width=0.7pt,dashed] (3.,5.5)-- (3.,7.);
\draw [line width=0.7pt,dashed] (5.,8.)-- (5.,6.5);
\draw [line width=0.7pt,dashed] (3.,7.)-- (5.,8.);
\draw [opacity=0.2, line width=0.7pt,dashed] (3.,7.)-- (14.,5.);
\draw [line width=0.7pt,dashed] (14.,5.)-- (16.,6.);
\draw [line width=0.7pt,dashed] (5.,8.)-- (16.,6.);
\draw [line width=0.7pt,dashed] (3.,5.5)-- (14.,3.5);
\draw [line width=0.7pt,dashed] (14.,3.5)-- (16.,4.5);
\draw [line width=0.7pt,dashed] (14.,5.)-- (14.,3.5);
\draw [line width=0.7pt,dashed] (16.,6.)-- (16.,4.5);
\draw [opacity=0.2, line width=0.7pt,dashed] (5.,6.5)-- (16.,4.5);
\draw [-to,line width=1pt] (3.,5.5) -- (3.,6.1);
\draw [-to,line width=1pt] (3.,5.5) -- (3.6,5.8);
\draw [-to,line width=1pt] (3.,5.5) -- (3.6,5.390909090909092);
\draw (11.3,7.572624886659496) node[anchor=north west] {$V$};
\draw [color=qqwwzz](8.980779397508933,6.5) node[anchor=north west] {$\Xi_1$};
\draw [color=xfqqff](10.50374992179453,5.075951896027373) node[anchor=north west] {$\Xi_2$};
\draw [color=xfqqff](13.075323102145623,6.1) node[anchor=north east] {$\Xi_2$};
\draw (3.1,8.25) node[anchor=north west] {$A_{D}^{\mathbf{u}}$};
\draw (3.55,5.4) node[anchor=north west] {$x$};
\draw (3.55,6) node[anchor=north west] {$y$};
\draw (2.45,6.55) node[anchor=north west] {$z$};
\begin{scriptsize}
\draw [fill=uuuuuu] (3.,5.5) circle (2.0pt);
\end{scriptsize}
\end{tikzpicture}
    		\caption{}
    	\end{subfigure}
            \begin{subfigure}{0.32\linewidth}
    		\centering
    		\includegraphics[width=0.9\linewidth]{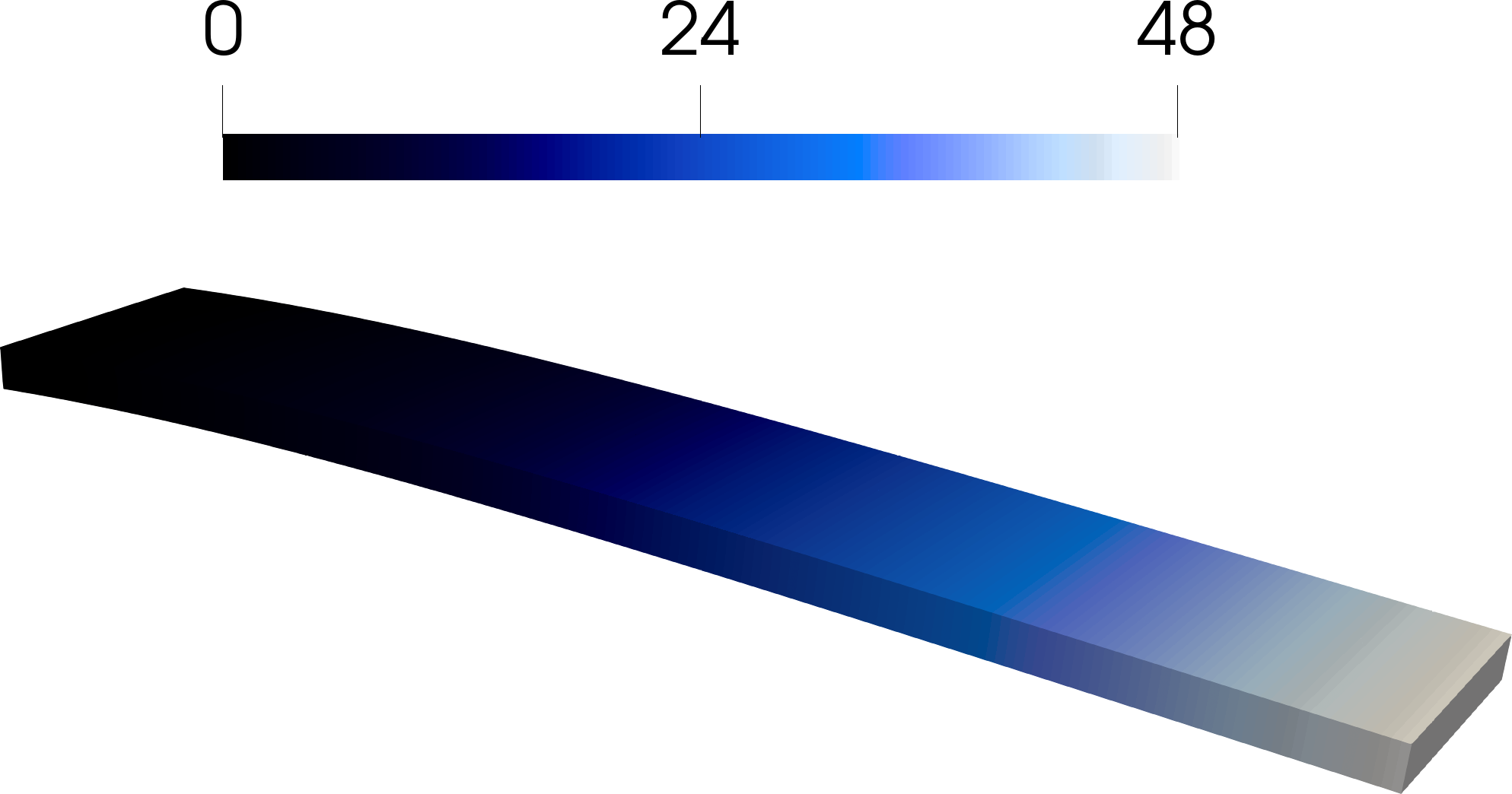}
    		\caption{}
    	\end{subfigure}
    	\caption{Depiction of the bulk domain $\vol$ with its Dirichlet boundary $\surf_D^{\ud}$, the codimensional domain of the embedded plate $\Xi_1$, and the centroid lines of the beams on both sides of the plate given by $\Xi_2$ (a). Deformation of the cantilever domain without any reinforcement (b) in $\mm$.}
    \label{fig:canti1}
\end{figure}
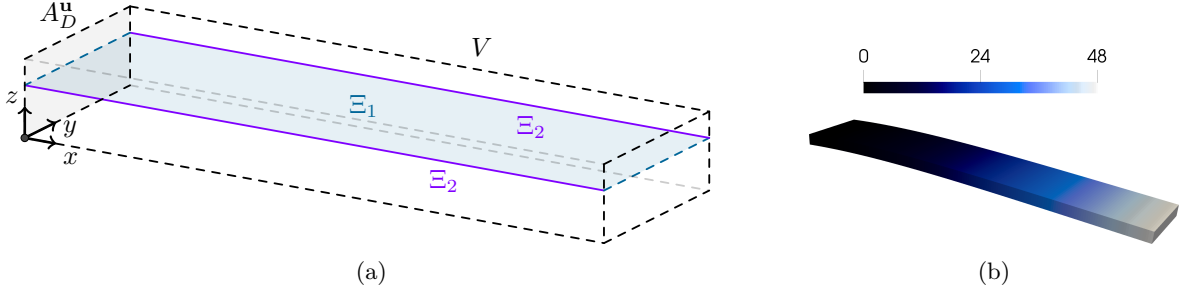

\begin{figure}
    	\centering
    	\begin{subfigure}{0.32\linewidth}
    		\centering
                \includegraphics[width=0.9\linewidth]{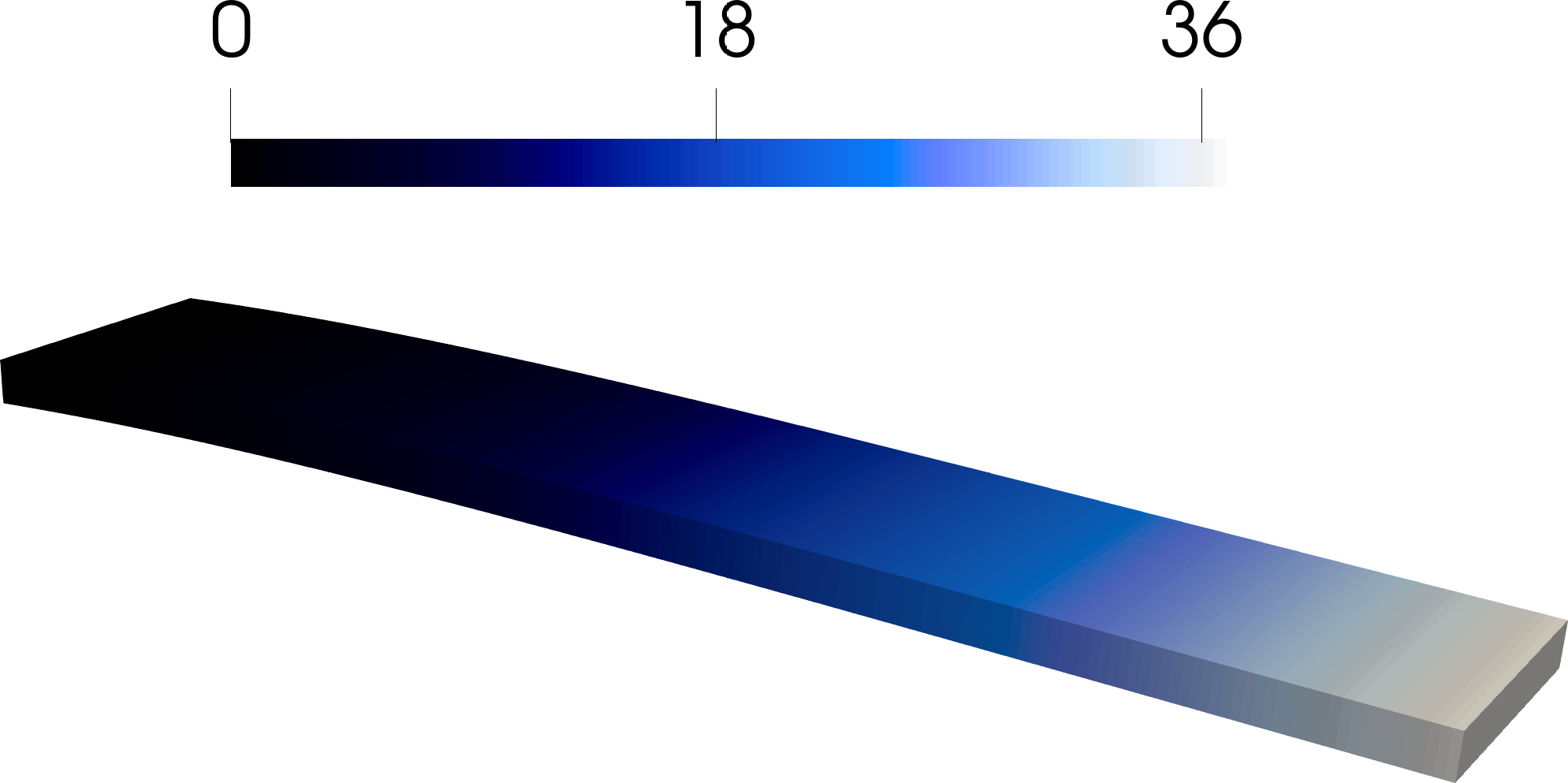}
    		\caption{}
    	\end{subfigure}
    	\begin{subfigure}{0.32\linewidth}
    		\centering
    		\includegraphics[width=0.9\linewidth]{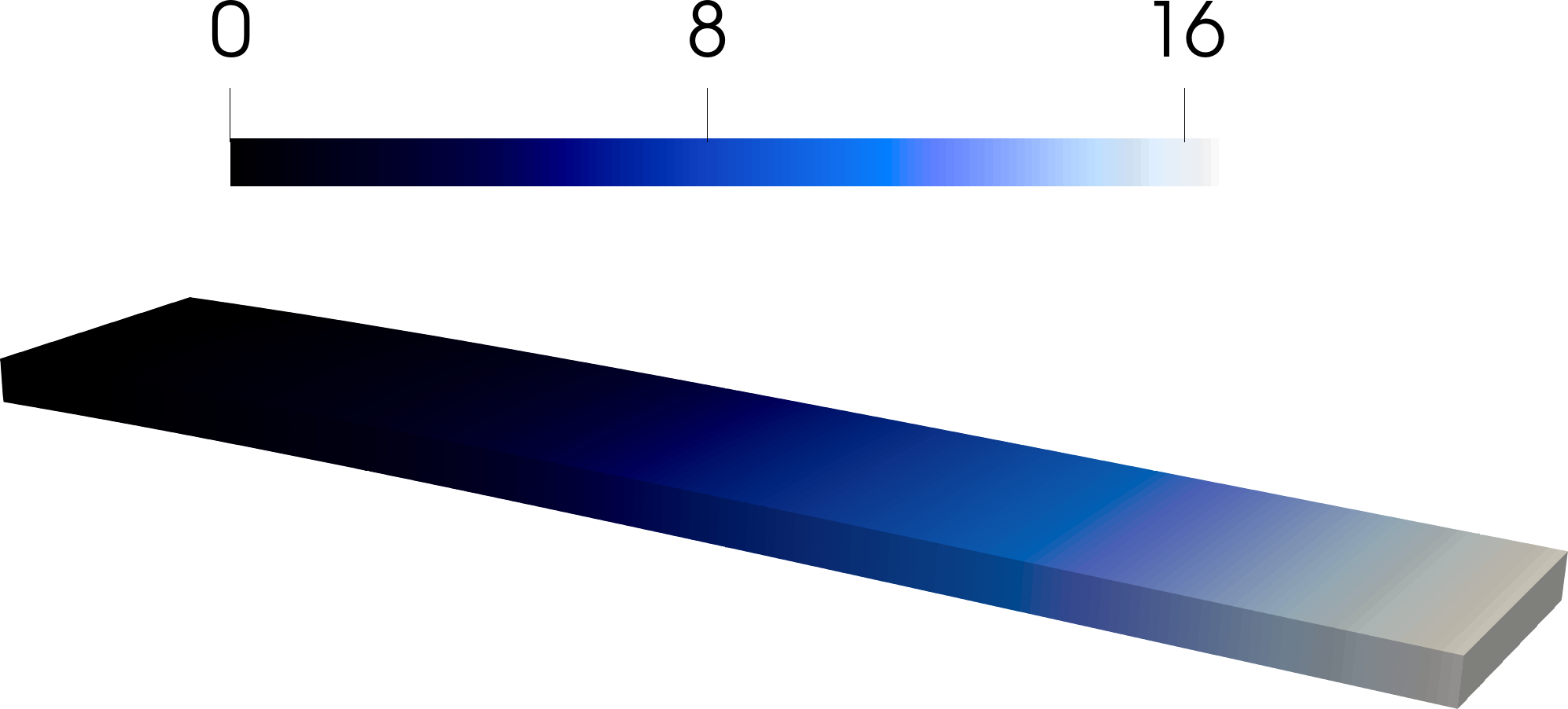}
    		\caption{}
    	\end{subfigure}
            \begin{subfigure}{0.32\linewidth}
    		\centering
    		\includegraphics[width=0.9\linewidth]{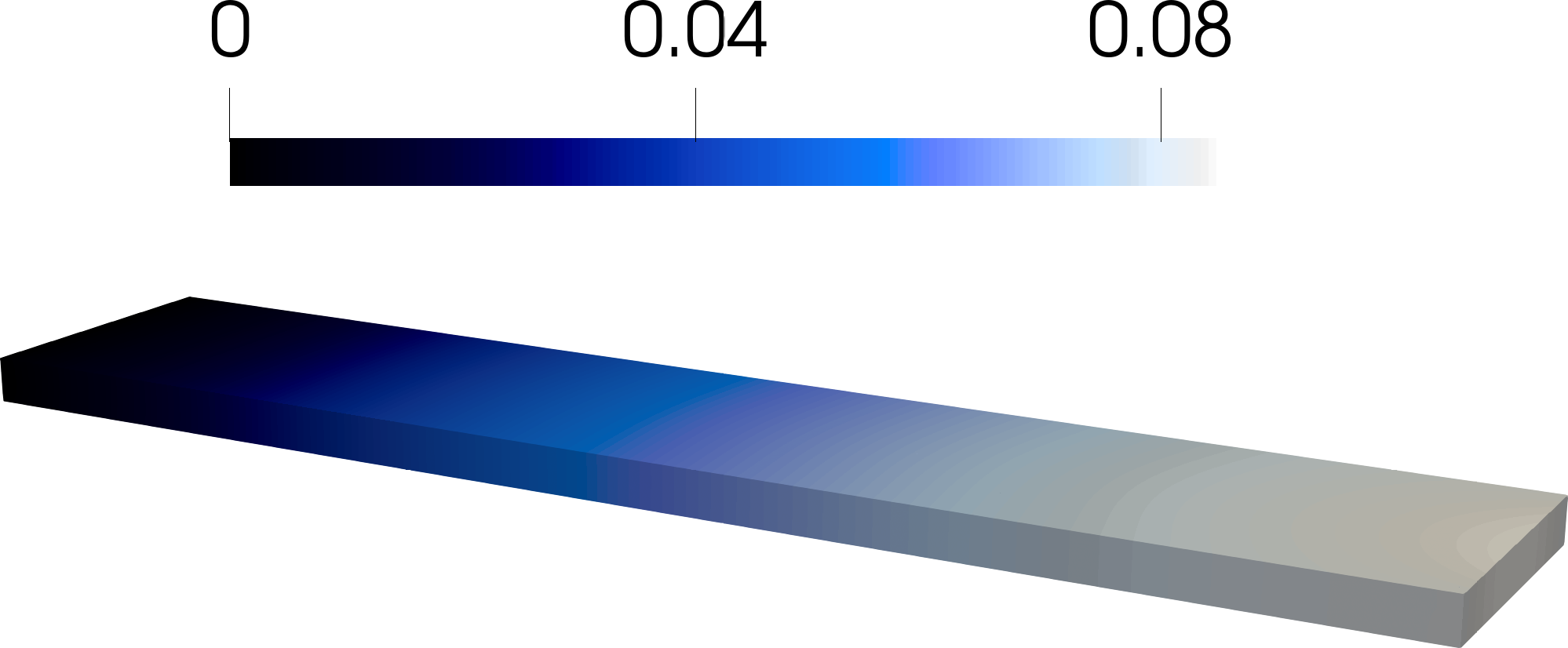}
    		\caption{}
    	\end{subfigure}
    	\caption{Deformations of the cantilever domain with beam reinforcements (a), plate reinforcement (b), and the combined plate and beam reinforcements (c) in $\mm$.}
    \label{fig:canti2}
\end{figure}

In the following we consider a simple beam shaped domain
\begin{align}
    \overline{\vol} = [0,500] \times [0,100] \times [0,15] \, ,
\end{align}
such that its length is $500 \, \mm$, its width is $100 \, \mm$ and its thickness is $15 \, \mm$.
The domain is loaded with the vertical volume force
\begin{align}
    &\vb{f} = -f \vb{e}_3 \, , && f = 10^{-5} \, \Nwtn / \mm^3 \, .
\end{align}
The bulk of the domain is set to be made of silver, such that its Lam\'e parameters read
\begin{align}
    \lama \approx 98.5 \, \MPa \, , && \mua \approx 30 \, \MPa \, ,
\end{align}
implying a Poisson ration of about $\nu \approx 0.38$. We define the Dirichlet boundary for the displacement as
\begin{align}
    \ud \at_{\surf_D^\ud} = 0 \, , && \surf_D^\ud = [0,0] \times [0,100] \times [0,15] \, ,
\end{align}
implying the kinematics of a cantilever beam. The complete boundary of the rotation field is set to be homogeneous Neumann $\surf_N^{\Rtheta} = \partial \vol$. Testing for various values of the characteristic $\Lc$ while setting $\muc = 1$ and $\Lm\Curl \Rtheta = \id(\Curl \Rtheta)=\Curl \Rtheta$, we find the maximal deformation for $\Lc = 10^{-2}$, such that lower values of $\Lc$ do not increase the deformation. As such, we conclude that for $\Lc = 10^{-2}$ we retrieve the equivalent Navier--Cauchy solution. We successively reinforce the bulk domain with a plate and beams made out the graphite @$1.6 \, \mm$ (H237) from the previous example. The thickness of the plate is again set to $h = 1.6 \, \mm$. The beam is defined with a full circular cross-section of radius $r = 0.8 \, \mm$, such that its surface and moments of inertia read
\begin{align}
    &A = \pi r^2 = 0.64 \pi \, \mm^2 \, , && I_\eta = I_\zeta = \dfrac{1}{4} \pi r^4 = 0.1024 \pi \, \mm^4 \, .
\end{align}
Due to its full symmetry the choice of its orientation given by the $\vb{n}$ and $\vb{c}$ vectors is inconsequential.
The plate reinforcement of the bulk is on codimension one. Its surface is defined as
\begin{align}
    \overline{\Xi}_1 = [0,500] \times [0,100] \times [10,10] \, .
\end{align}
Further, on two sides of $\Xi_1$ we introduce beam reinforcements
\begin{align}
    \overline{\Xi}_2 = \{[0,500] \times [0,0] \times [10,10]\} \cup \{[0,500] \times [100,100] \times [10,10]\} \, ,
\end{align}
on codimension one of the surface, representing codimension two of the bulk. Thus, we embed in the volume at $z = 10 \, \mm$ a plate, and two beams at the same depth with $y = 0$ and $y = 100$. The domain with its reinforcement is illustrated in \cref{fig:canti1}.
The maximal deformation without any reinforcement is also depicted in \cref{fig:canti1}, reading $\max \norm{\ud} \approx 48.49 \, \mm$. After reinforcement with two beams the maximal deformation reduces to $\max \norm{\ud} \approx 37.04 \, \mm$. Alternatively, reinforcing the bulk with the plate yields $\max \norm{\ud} \approx 16.71 \, \mm$. Finally, combining both reinforcements leads to the maximal deformation $\max \norm{\ud} \approx 0.08 \, \mm$. The respective deformation results are depicted in \cref{fig:canti2}.
We observe the expected kinematical behaviour. Firstly, the reinforcement with beams is less pronounced than the effect of reinforcing with the plate. Secondly, the combined reinforcements yield a stiffness-additive solution, which is the natural outcome for linear elastic mechanics. We note that we do not retrieve an absolute super-positional solution, which is presumably due to the stiffening effect of the rotations also in $x$- and $z$-directions.   

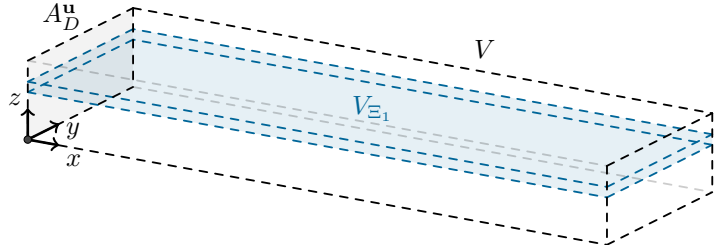
\begin{figure}
    \centering
    \definecolor{uuuuuu}{rgb}{0.26666666666666666,0.26666666666666666,0.26666666666666666}
\definecolor{xfqqff}{rgb}{0.4980392156862745,0.,1.}
\definecolor{qqwwzz}{rgb}{0.,0.4,0.6}
\begin{tikzpicture}[scale = 0.7, line cap=round,line join=round,>=triangle 45,x=1.0cm,y=1.0cm]
\clip(2.5,3.5) rectangle (16.5,8.5);

\fill[line width=0.7pt,color=qqwwzz,fill=qqwwzz,fill opacity=0.10000000149011612] (3.,6.5+0.1) -- (14.,4.5+0.1) -- (16.,5.5+0.1) -- (5.,7.5+0.1) -- cycle;
\draw [line width=0.7pt,dashed,color=qqwwzz] (3.,6.5+0.1)-- (14.,4.5+0.1);
\draw [line width=0.7pt,dashed,color=qqwwzz] (14.,4.5+0.1)-- (16.,5.5+0.1);
\draw [line width=0.7pt,dashed,color=qqwwzz] (16.,5.5+0.1)-- (5.,7.5+0.1);
\draw [line width=0.7pt,dashed,color=qqwwzz] (5.,7.5+0.1)-- (3.,6.5+0.1);

\fill[line width=0.7pt,color=qqwwzz,fill=qqwwzz,fill opacity=0.10000000149011612] (14.,4.5+0.1) -- (16.,5.5+0.1) -- (16.,5.5-0.1) -- (14.,4.5-0.1) -- cycle;

\fill[line width=0.7pt,color=qqwwzz,fill=qqwwzz,fill opacity=0.10000000149011612] (3.,6.5+0.1) -- (14.,4.5+0.1) -- (14.,4.5-0.1) -- (3.,6.5-0.1)  -- cycle;

\fill[line width=0.7pt,fill=black,fill opacity=0.05] (3.,7.) -- (5.,8.) -- (5.,6.5) -- (3.,5.5) -- cycle;
\draw [line width=0.7pt,dashed,color=qqwwzz] (3.,6.5-0.1)-- (14.,4.5-0.1);
\draw [line width=0.7pt,dashed,color=qqwwzz] (14.,4.5-0.1)-- (16.,5.5-0.1);
\draw [line width=0.7pt,dashed,color=qqwwzz] (16.,5.5-0.1)-- (5.,7.5-0.1);
\draw [line width=0.7pt,dashed,color=qqwwzz] (5.,7.5-0.1)-- (3.,6.5-0.1);

\draw [line width=0.7pt,dashed] (3.,5.5)-- (5.,6.5);
\draw [line width=0.7pt,dashed] (3.,5.5)-- (3.,7.);
\draw [line width=0.7pt,dashed] (5.,8.)-- (5.,6.5);
\draw [line width=0.7pt,dashed] (3.,7.)-- (5.,8.);
\draw [opacity=0.2, line width=0.7pt,dashed] (3.,7.)-- (14.,5.);
\draw [line width=0.7pt,dashed] (14.,5.)-- (16.,6.);
\draw [line width=0.7pt,dashed] (5.,8.)-- (16.,6.);
\draw [line width=0.7pt,dashed] (3.,5.5)-- (14.,3.5);
\draw [line width=0.7pt,dashed] (14.,3.5)-- (16.,4.5);
\draw [line width=0.7pt,dashed] (14.,5.)-- (14.,3.5);
\draw [line width=0.7pt,dashed] (16.,6.)-- (16.,4.5);
\draw [opacity=0.2, line width=0.7pt,dashed] (5.,6.5)-- (16.,4.5);
\draw [-to,line width=1pt] (3.,5.5) -- (3.,6.1);
\draw [-to,line width=1pt] (3.,5.5) -- (3.6,5.8);
\draw [-to,line width=1pt] (3.,5.5) -- (3.6,5.390909090909092);
\draw (11.3,7.572624886659496) node[anchor=north west] {$V$};
\draw [color=qqwwzz](8.980779397508933,6.5) node[anchor=north west] {$V_{\Xi_1}$};
\draw (3.1,8.25) node[anchor=north west] {$A_{D}^{\mathbf{u}}$};
\draw (3.55,5.4) node[anchor=north west] {$x$};
\draw (3.55,6) node[anchor=north west] {$y$};
\draw (2.45,6.55) node[anchor=north west] {$z$};
\begin{scriptsize}
\draw [fill=uuuuuu] (3.,5.5) circle (2.0pt);
\end{scriptsize}
\end{tikzpicture}
    \caption{Fully volumetric definition of the plate-reinforced bulk model with the explicit thickness $h$.}
    \label{fig:cantivol}
\end{figure}

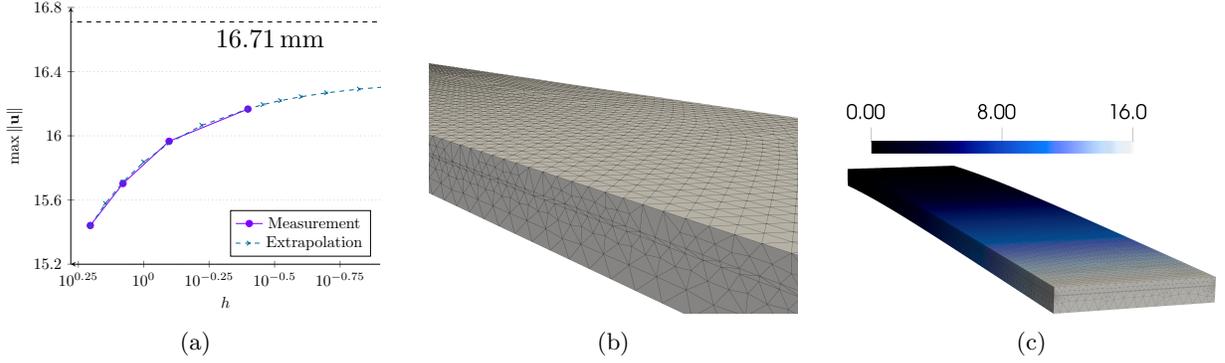
\begin{figure}
    \centering
    \begin{subfigure}{0.32\linewidth}
    		\centering
    		\definecolor{asb}{rgb}{0.,0.4,0.6}
\definecolor{asl}{rgb}{0.4980392156862745,0.,1.}
\begin{tikzpicture}[scale = 0.6]
	\begin{semilogxaxis}[
				/pgf/number format/1000 sep={},
				axis lines = left,
				xlabel={$h$},
				ylabel={$\max\norm{\ud}$},
				xmin=0.125, xmax=1.9, 
				ymin=15.2, ymax= 16.8, 
				x dir=reverse,
				xtick={0.1778279410038923, 0.31622776601683794, 0.5623413251903491, 1, 1.7782794100389228},
				ytick={15.2, 15.6, 16, 16.4, 16.8},
				legend pos=south east,
				ymajorgrids=true,
				grid style=dotted,
				]
				\addplot[
				color=asl,
				mark=*,
				]
				coordinates {
                        ( 1.6 , 15.44137668137552 )
                        ( 1.2 , 15.702873763963943 )
                        ( 0.8 , 15.966197152533859 )
                        ( 0.4 , 16.166876867190897 )
				};
				\addlegendentry{Measurement}
				\addplot[
				color=asb,
                dashed,
				mark=x,
				]
				coordinates {
                        ( 1.6 , 15.437883739730717 )
                        ( 1.4 , 15.579193017962087 )
                        ( 1.2 , 15.712155828458645 )
                        ( 1 , 15.837265157417198 )
                        ( 0.8 , 15.954984872681973 )
                        ( 0.6 , 16.06575144362738 )
                        ( 0.4 , 16.166876867190897 )
                        ( 0.35 , 16.195055126946855 )
                        ( 0.3 , 16.219755864408473 )
                        ( 0.25 , 16.24408349411568 )
                        ( 0.2 , 16.268043651908165 )
                        ( 0.15 , 16.29164188849554 )
                        ( 0.1 , 16.31488367074325 )
                        ( 0.05 , 16.337774382939045 )
                        ( 0.0001, 16.36031932804034 )
				};
				\addlegendentry{Extrapolation}
				\addplot[dashed,color=black, mark=none]
				coordinates {
					(0.1, 16.71)
					(2, 16.71)
				};
			\end{semilogxaxis}
			\draw (3,5.4) node[anchor=north west]{$16.71 \, \mm$};
		\end{tikzpicture}
    		\caption{}
    	\end{subfigure}
    	\begin{subfigure}{0.32\linewidth}
    		\centering
    		\includegraphics[width=0.9\linewidth]{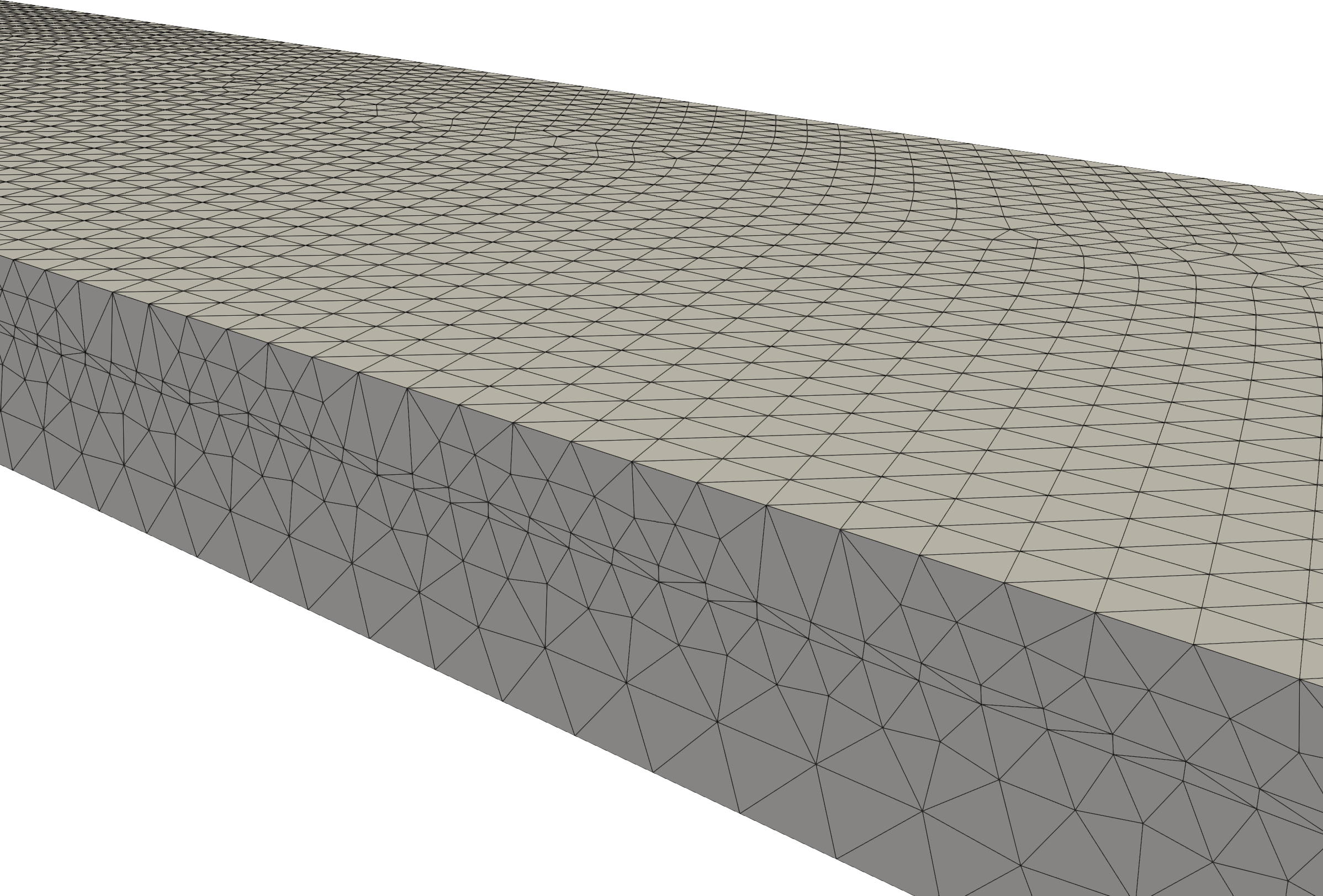}
    		\caption{}
    	\end{subfigure}
            \begin{subfigure}{0.32\linewidth}
    		\centering
    		\includegraphics[width=0.9\linewidth]{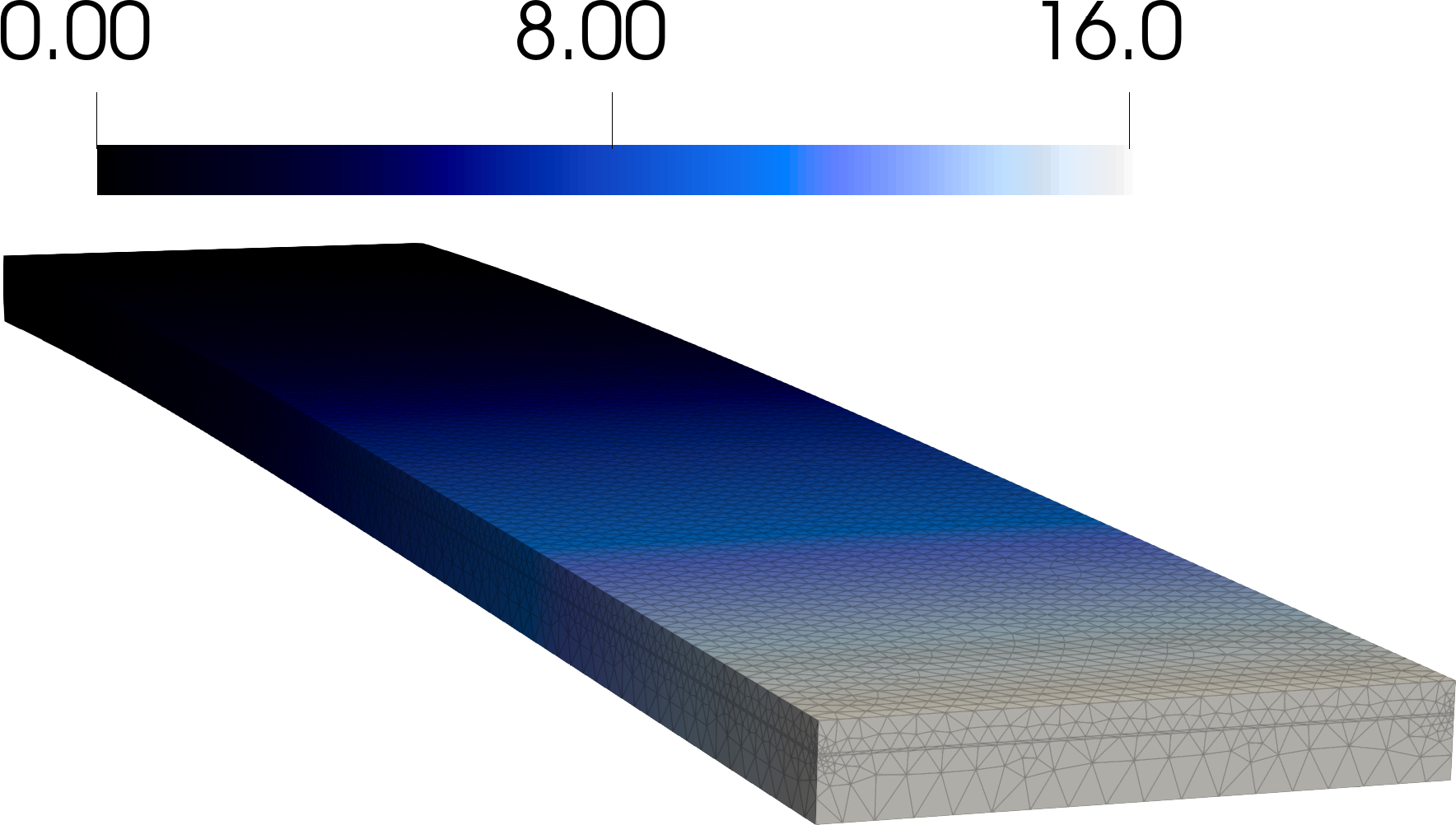}
    		\caption{}
    	\end{subfigure}
    	\caption{Convergence of the maximal displacement towards the limit of $\approx 16.36 \, \mm$ (a). Volumetric Mesh of the cantilever beam domain with $h = 0.8 \, \mm$, demonstrating the growing amount of elements which are need for scale-transition (b). Norm of the displacement field of the full volumetric model for $h = 0.4 \, \mm$ on a mesh of $129008$ quadratic elements with $1118430$ degrees of freedom (c).}
    \label{fig:conv}
\end{figure}

Due to the lack of analytical results or a sound mixed-dimensional error estimator it is difficult to undertake convergence studies in order to estimate the quality of the mixed-dimensional approximation. Comparisons with fully volumetric discretisations are also restricted, as the reduced models are an idealisation of the three-dimensional kinematics, incapable of capturing the full mechanical phenomenon. Nevertheless, we propose a weighted bulk benchmark to evaluate the approximation power of the coupled plate model to some extent. The plate model allows for a rudimentary comparison due to its simple geometry, in contrast to the additional errors induced by the non-smooth geometrical approximation of curved shells. We start by splitting the bulk domain across two volumes
\begin{align}
    &\overline{\vol} = \overline{\vol}_\vol \cup \overline{\vol}_{\Xi_1} \, ,   && \begin{aligned}
        \overline{\vol}_\vol &= \{[0,500] \times [0,100] \times [0,10-h/2] \} \cup \{[0,500] \times [0,100] \times [10+h/2,15] \} \, , \\ 
        \overline{\vol}_{\Xi_1} &= [0,500] \times [0,100] \times [10-h/2,10+h/2] \, ,
    \end{aligned}
\end{align}
such that $\overline{\vol}_\vol$ is the domain made of silver, and $\overline{\vol}_{\Xi_1}$ represents the domain of the plate with thickness $h$ made of graphite, see \cref{fig:cantivol}.
Clearly, the larger $h$ is, the less silver the domain $\overline{\vol}$ is composed of, and the stiffer the it becomes. However, in the mixed-dimensional model the thickness of the plate does not explicitly appear in the discretisation of the geometry, such that it does not affect the volume of the silver-made domain. At the same time, reducing $h$ in the fully volumetric model in order to account for this inconsistency reduces the overall stiffness of the model, although it remains constant in the corresponding mixed-dimensional model $h = 1.6 \, \mm$. Thus, we study the behaviour of the comparable fully volumetric model by decreasing the thickness $h$, while simultaneously compensating for this loss in stiffness with a linear scaling factor $k = 1.6/h \geq 1$. This multiplicative factor is applied to the material coefficients of the graphite-made domain. We use quadratic Lagrange elements for the computation with $h \in \{ 1.6, 1.2, 0.8, 0.4\} \, \mm$, such that the scaling factor reads $k \in \{ 1, 4/3, 2,  4\}$.    
The convergence, along with the mesh for $h = 0.8 \, \mm$ and the displacement for $h = 0.4 \, \mm$ are depicted in \cref{fig:conv}.
Clearly, taking the value $h = 0$ is impossible as it implies a zero volume and $k \to + \infty$. Even for $h \leq 0.3$ it already becomes difficult to generate satisfactory meshes that are manageable from a computational perspective. As such, we instead extrapolate the remaining values of the convergence curve with a curve-fitting based on a generic exponential function $f(h) = a e^{b h} + c$ with $a,b,c \in \R$. As seen in \cref{fig:conv}, the fitting perfectly matches the computed values and predicts a maximal displacement of $\max \norm{\ud} \approx 16.36 \, \mm$ for $h \to 0$. Thus, the relative difference in the maximal displacement between the fully volumetric model and the mixed-dimensional model is $100 (16.71 - 16.36) / 16.71 \approx 2.09 \, \%$ or $16.71 - 16.36 = 0.35 \, \mm$ in absolute values. This relatively small difference implies a good agreement between the two models, and demonstrates the accurate prediction power of the mixed-dimensional approach. Further, since the mixed-dimensional model slightly under-estimates the stiffness of the system and therefore over-estimates the deformation, the result is on the so-called safe side for subsequent design decisions.  

\subsection{An \reflectbox{S}-shaped reinforced shell}
The proposed coupling procedure is general across dimensions due to the degrees of freedom of the underlying linear Cosserat continuum. In particular, shells naturally include in-plane drill rotations. The fact that a three-dimensional rotation vector composed of three independent rotation fields is present in every model irrespective of dimensionality allows for seamless interaction \textbf{even across intersections}. This feature is greatly facilitated by the use of tangential differential calculus (TDC), since both the displacement field $\ud$ and the rotation field $\rtheta$ are defined on the Cartesian system. In other words, whether a rotation represents drill for one shell but bending for another is determined solely by the projections and no explicit transformations to the coordinate system of any hyper-surface are required. The same holds true analogously for intersecting beams with respect to bending and torsion. The final example serves to demonstrate this feature for 2D-1D couplings. 

\begin{figure}
    	\centering
    	\begin{subfigure}{0.32\linewidth}
    		\centering
                \input{figs/skel}
    		\caption{}
    	\end{subfigure}
    	\begin{subfigure}{0.32\linewidth}
    		\centering
    		\includegraphics[width=0.9\linewidth]{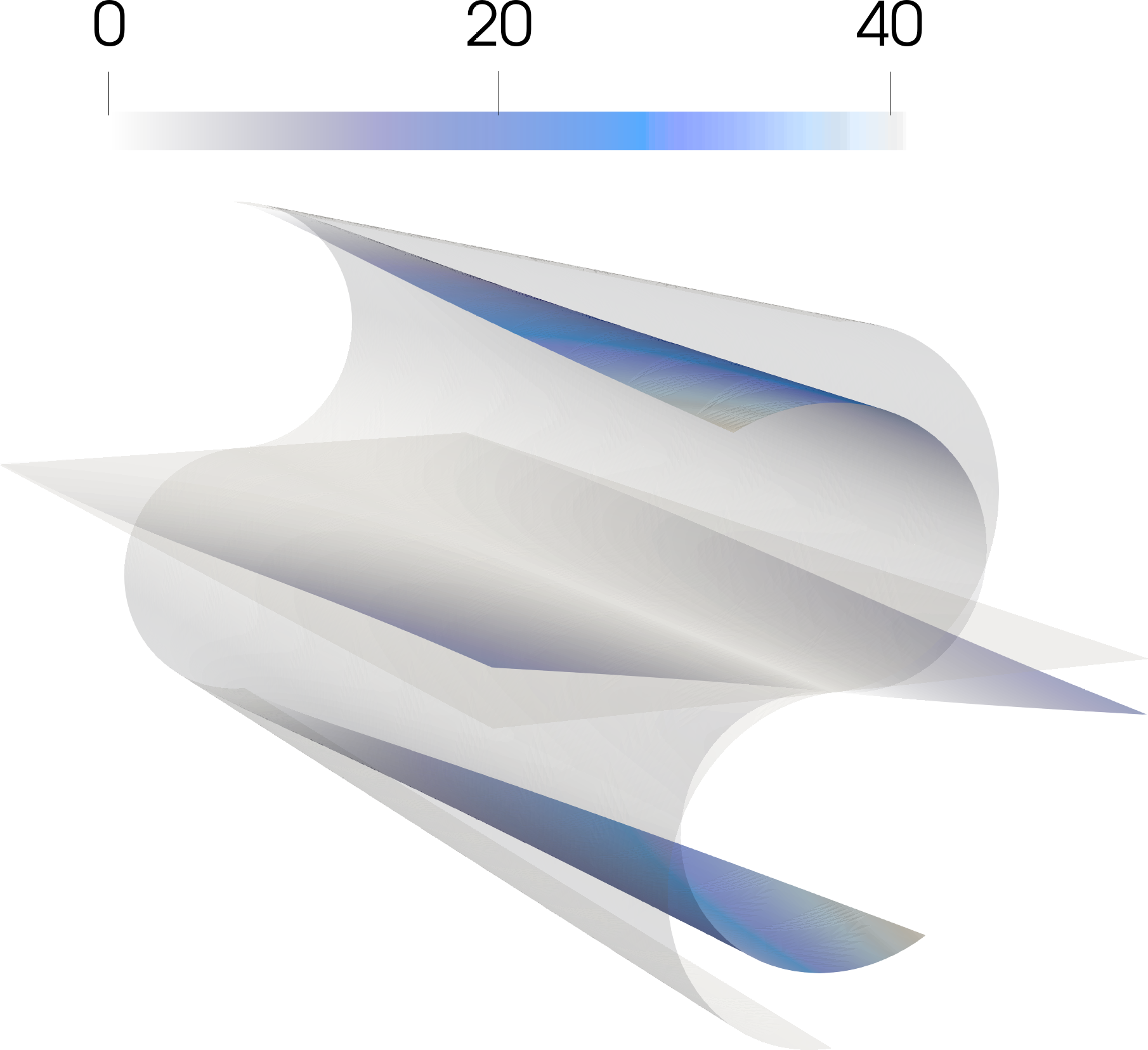}
    		\caption{}
    	\end{subfigure}
            \begin{subfigure}{0.32\linewidth}
    		\centering
    		\includegraphics[width=0.9\linewidth]{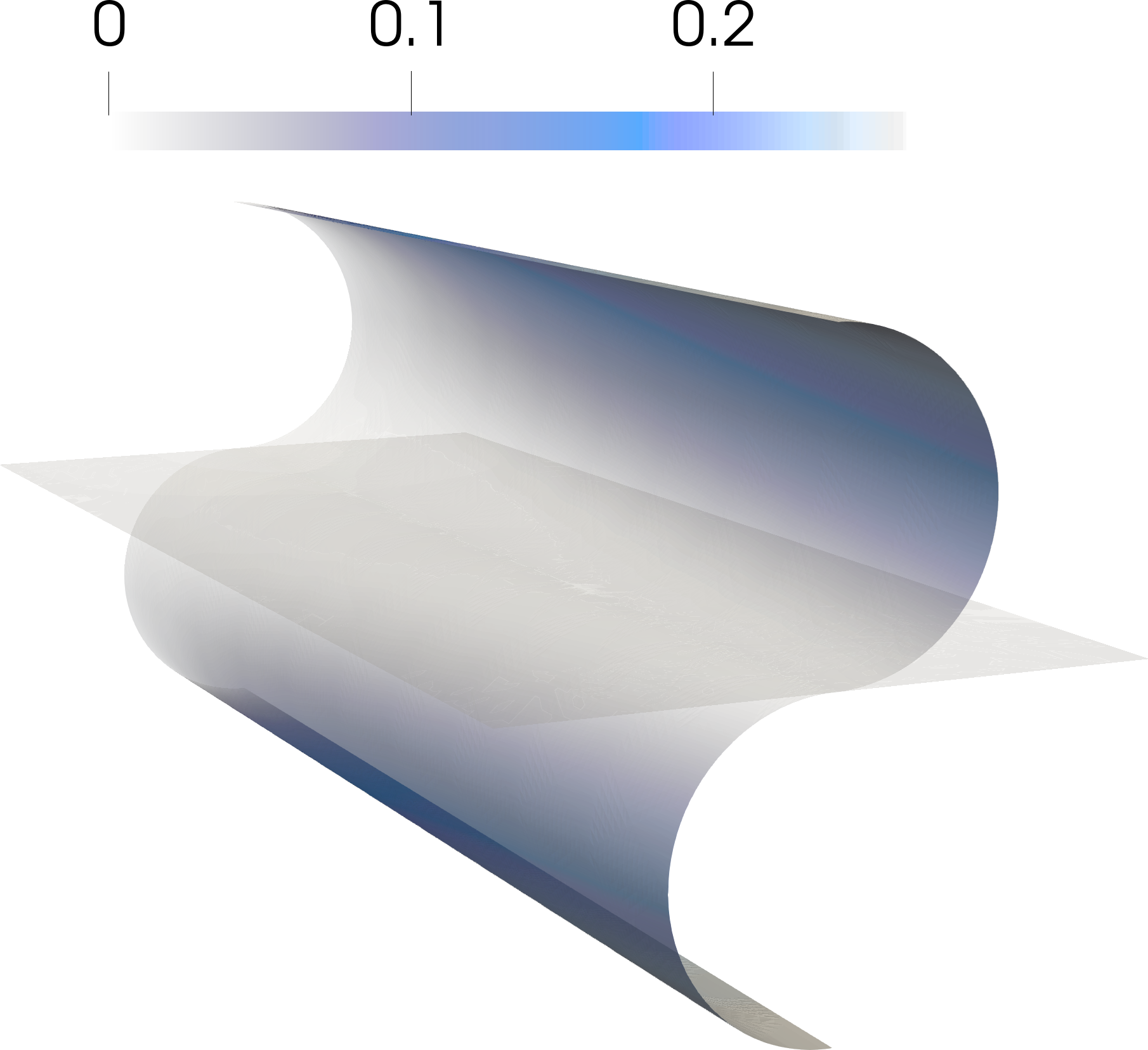}
    		\caption{}
    	\end{subfigure}
    	\caption{Domain of the shell $\surf$ with its codimensional reinforcement-frame $\Xi$ and Dirichlet boundary $\curv_D$ (a). Deformation of the shell without reinforcement (b), and with reinforcement (c).}
    \label{fig:S}
\end{figure}

We consider a reflected S-shaped (\reflectbox{S}-shape) domain constructed of two half-circles
\begin{align}
    \overline{\curv}_1 &= \{ \phi \in [\pi/2,  3\pi/4] \; | \; \vb{r} = \begin{bmatrix}
        0 & 50 \cos \phi  & 50 (\sin \phi - 1)    
    \end{bmatrix}^T \} \, , \notag \\
    \overline{\curv}_2 &= \{ \phi \in [-\pi/2,  \pi/2] \; | \; \vb{r} = \begin{bmatrix}
        0 & 50 \cos \phi  & 50 (\sin \phi + 1)    
    \end{bmatrix}^T \} \, ,
\end{align}
with a line crossing their intersection 
\begin{align}
    \overline{\curv}_3 = [0,-100,0] \times [0,100,0] \, .
\end{align}
The shell surface is then defined via the extrusion
\begin{align}
    \overline{\surf} = [0, 500] \times \{\overline{\curv}_1 \cup \overline{\curv}_2 \cup \overline{\curv}_3 \} \, .
\end{align}
Accordingly, the radius of the circles is $50 \, \mm$, the middle plate has the width of $200 \, \mm$, and the length of the structure is $500 \, \mm$. The Dirichlet boundary $\curv_D$ is defined on the edge curves at $x = 0$ and is applied such that both the displacements and rotations vanish 
\begin{align}
    \ud \at_{\curv_D} =  \rtheta \at_{\curv_D} = 0 \, .
\end{align}
The shell is set to be made of silver with thickness $h = 1.6 \, \mm$ and is subsequently reinforced with beams of a circular cross-section with $r = 0.8 \, \mm$ made of graphite. Curved beams are defined along the \reflectbox{S}-curves, whereas straight lines are defined with the straight beam model
\begin{align}
    \overline{\Xi} = & \{ (\eta, \phi) \in \{0,1,\dots,10\} \times [\pi/2,  3\pi/4] \; | \; \vb{r} = 50\begin{bmatrix}
        \eta & \cos \phi  & (\sin \phi - 1)    
    \end{bmatrix}^T  \} \notag \\
    &\cup \{ (\eta,\phi) \in \{0,1,\dots,10\} \times [-\pi/2,  \pi/2] \; | \; \vb{r} = 50\begin{bmatrix}
        \eta & \cos \phi  & (\sin \phi - 1)    
    \end{bmatrix}^T \} 
    \notag \\
    &\cup \{ (x,y) \in [0,500] \times \{-100,0,100\} \; | \; \vb{r} =\begin{bmatrix}
        x & y & 0
    \end{bmatrix}^T  \}
    \\
    &\cup \{ (x,y) \in \{0,500\} \times [-100,100]  \; | \; \vb{r} = \begin{bmatrix}
        x & y & 0
    \end{bmatrix}^T  \}
     \, , \notag
\end{align}
representing codimensional domains on the surface of the shell.
The corresponding material coefficients can be found in the previous examples. The force is applied as a line load in $\Nwtn/\mm$ to the upper and lower horizontal edges of the shell structure
\begin{align}
    \vb{q} = q \begin{bmatrix}
        0 \\ 1 \\ 0 
    \end{bmatrix} \, , && q = \left \{\begin{matrix}
        -10^{-4} \, z   & \text{for} & z = \pm 100  \\
        0 & \text{otherwise}
    \end{matrix} \right . \, ,
\end{align}
implying torsion.
The domain along with its reinforcing frame and the resulting displacements for $620$ cubic elements are depicted in \cref{fig:S}.
We clearly observe the resulting warp-torsion deformation of the non-reinforced thin structure with a maximal displacement of $\max \norm{\ud} \approx 41.14 \,  \mm$. Activating the beam-frame reinforcement reduces the maximal displacement to  $\max \norm{\ud} \approx 0.26 \,  \mm$. The intersection of multiple shells and beams with differing orientations poses no challenge to the proposed coupling approach, as evident by the symmetric displacement solutions. 

\section{Conclusions and outlook}
In this work we employed the linear Cosserat micropolar model in dislocation form to derive corresponding shell, plate, and beam models via kinematical and dimensional reduction. Further, we have shown that for $\Lc \to 0$ it is possible to recover the Navier--Cauchy response while using the Cosserat model.
The reduced dimensional models were derived by standard engineering assumptions for the kinematics and splits in the integration over the volumetric domain.  
For the derivation we made use of the modern framework of tangential differential calculus. Thus, we were able to derive the corresponding energy functionals without the need for curvilinear coordinates or Christoffel symbols. Another major advantage of this approach is its applicability to automated solvers of partial differential equations, where the energy functional could be written as is, and solved directly without further treatment. 
All the various dimensional models intrinsically share the same type of degrees of freedom, namely translations and rotations. As a consequence, we were able to define coupled systems using merely consistent Sobolev trace operators, without the need for intermediate finite elements or Lagrange multipliers. The natural $\C^0(\body)$-regularity of continuous Lagrange elements ensured that its Sobolev traces on codimensional domains remain well-defined and square integrable also on these lower dimensional manifolds. We were thus able to validate our approach with three numerical examples computed in NGSolve (\AS{see Fn.~\ref{fn:github}}), where we coupled a volume with a shell, a volume and with a plate and beams, or a shell with beams.  

The approach presented in this work has three advantages. Firstly, it is simple as no special treatment is required for the coupling. Instead, the coupling is done by defining additional energy functionals for Sobolev traces of the displacement and rotation fields on codimensional domains. 
Secondly, it is arbitrarily valid across dimensions, even for intersections.
Thirdly, although not explicitly discussed in this work, the coupling modulus $\muc$ can be understood as a penalty term or a rotational spring. In other words, it could be used to control how strong the coupling is, allowing for deltas on interfaces of mixed-dimensional domains. Collectively, we surmise that our approach is relevant for engineering designs with mixed-dimensional parts, such as sandwich elements in aerospace, or general composite and fibre-reinforced materials.  

Notwithstanding, the presented approach does entail one disadvantage, being its additional rotational degrees of freedom on all volumetric elements, making it computationally more expensive than the comparable Navier--Cauchy model for isotropic Cauchy materials.   

\section*{Acknowledgements}
Adam Sky is grateful for the technical support by Christopher Lackner (cerbsim, Technical University of Vienna) of the NGSolve community. Patrizio Neff acknowledges support in the framework of the DFG-Priority Programme
2256 “Variational Methods for Predicting Complex Phenomena in Engineering Structures and Materials”, Neff 902/10-1, Project-No. 440935806. 

\bibliographystyle{spmpsci}   

\footnotesize
\bibliography{ref}   

\begin{thebibliography}{10}
\providecommand{\url}[1]{{#1}}
\providecommand{\urlprefix}{URL }
\expandafter\ifx\csname urlstyle\endcsname\relax
  \providecommand{\doi}[1]{DOI~\discretionary{}{}{}#1}\else
  \providecommand{\doi}{DOI~\discretionary{}{}{}\begingroup \urlstyle{rm}\Url}\fi

\bibitem{ARF2023116198}
Arf, J., Reichle, M., Klinkel, S., Simeon, B.: Scaled boundary isogeometric analysis with ${C}^1$ coupling for {Kirchhoff} plate theory.
\newblock Computer Methods in Applied Mechanics and Engineering \textbf{415}, 116198 (2023)

\bibitem{Banerjee}
Banerjee, J., Kennedy, D., Elishakoff, I.: Further insights into the {Timoshenko–Ehrenfest} beam theory.
\newblock Journal of Vibration and Acoustics \textbf{144}(6), 061011 (2022)

\bibitem{baratta_2023_10447666}
Baratta, I.A., Dean, J.P., Dokken, J.S., Habera, M., Hale, J.S., Richardson, C.N., Rognes, M.E., Scroggs, M.W., Sime, N., Wells, G.N.: {DOLFINx: The next generation FEniCS problem solving environment} (2023)

\bibitem{boon2024mixed}
Boon, W.M., Duran, O., Nordbotten, J.M.: Mixed finite element methods for linear {C}osserat equations.
\newblock arXiv  (2024).
\newblock \urlprefix\url{https://arxiv.org/abs/2403.15136v2}

\bibitem{BORKOVIC2023115848}
Borković, A., Gfrerer, M., Marussig, B.: Geometrically exact isogeometric {Bernoulli–Euler} beam based on the {Frenet–Serret} frame.
\newblock Computer Methods in Applied Mechanics and Engineering \textbf{405}, 115848 (2023)

\bibitem{BOURNIVAL2010838}
Bournival, S., Cuillière, J.C., François, V.: A mesh-geometry based method for coupling {1D} and {3D} elements.
\newblock Advances in Engineering Software \textbf{41}(6), 838--858 (2010)

\bibitem{BUCK2015159}
Buck, F., Brylka, B., Müller, V., Müller, T., Weidenmann, K.A., Hrymak, A.N., Henning, F., Böhlke, T.: Two-scale structural mechanical modeling of long fiber reinforced thermoplastics.
\newblock Composites Science and Technology \textbf{117}, 159--167 (2015)

\bibitem{BURMAN201851}
Burman, E., Hansbo, P., Larson, M.G.: A simple approach for finite element simulation of reinforced plates.
\newblock Finite Elements in Analysis and Design \textbf{142}, 51--60 (2018)

\bibitem{Mircea}
Bîrsan, M., Ghiba, I.D., Martin, R.J., Neff, P.: Refined dimensional reduction for isotropic elastic {Cosserat} shells with initial curvature.
\newblock Mathematics and Mechanics of Solids \textbf{24}(12), 4000--4019 (2019)

\bibitem{Choi2022}
Choi, M.J., Klinkel, S., Sauer, R.A.: An isogeometric finite element formulation for frictionless contact of {Cosserat} rods with unconstrained directors.
\newblock Computational Mechanics \textbf{70}(6), 1107--1144 (2022)

\bibitem{Cosserat}
Cosserat, E., Cosserat, F.: Th{\'e}orie des corps d{\'e}formables.
\newblock Nature \textbf{81}(2072), 67--67 (1909).
\newblock \urlprefix\url{https://www.uni-due.de/imperia/md/content/mathematik/ag_neff/cosserat09_eng.pdf}

\bibitem{Delfour}
Delfour, M., Zolsio, J.P.: Shapes and Geometries: Metrics, Analysis, Differential Calculus, and Optimization.
\newblock Society for Industrial and Applied Mathematics (2010)

\bibitem{DICARA2024117856}
{Di Cara}, G., D’Ottavio, M., Polit, O.: Variable kinematics finite plate elements for the buckling analysis of sandwich composite panels.
\newblock Composite Structures \textbf{330}, 117856 (2024)

\bibitem{Eringen1999}
Eringen, A.: Microcontinuum Field Theories. I. Foundations and Solids.
\newblock Springer-Verlag New York (1999)

\bibitem{Forest2006}
Forest, S., Sievert, R.: Nonlinear microstrain theories.
\newblock International Journal of Solids and Structures \textbf{43}(24), 7224--7245 (2006)

\bibitem{FRIES2023116223}
Fries, T.P., Kaiser, M.W.: On the simultaneous solution of structural membranes on all level sets within a bulk domain.
\newblock Computer Methods in Applied Mechanics and Engineering \textbf{415}, 116223 (2023)

\bibitem{Gangl2021}
Gangl, P., Sturm, K., Neunteufel, M., Sch{\"o}berl, J.: Fully and semi-automated shape differentiation in {NGSolve}.
\newblock Structural and Multidisciplinary Optimization \textbf{63}(3), 1579--1607 (2021)

\bibitem{Ghiba2020I}
Ghiba, I.D., B{\^i}rsan, M., Lewintan, P., Neff, P.: The isotropic {Cosserat} shell model including terms up to $\mathcal{O}(h^{5})$. {Part I}: Derivation in matrix notation.
\newblock Journal of Elasticity \textbf{142}(2), 201--262 (2020)

\bibitem{Ghiba2020II}
Ghiba, I.D., B{\^i}rsan, M., Lewintan, P., Neff, P.: The isotropic {Cosserat} shell model including terms up to $\mathcal{O}(h^{5})$. {Part II}: Existence of minimizers.
\newblock Journal of Elasticity \textbf{142}(2), 263--290 (2020)

\bibitem{Ghiba2021}
Ghiba, I.D., B{\^i}rsan, M., Lewintan, P., Neff, P.: A constrained {Cosserat} shell model up to order $\mathcal{O}(h^{5})$: Modelling, existence of minimizers, relations to classical shell models and scaling invariance of the bending tensor.
\newblock Journal of Elasticity \textbf{146}(1), 83--141 (2021)

\bibitem{Ghiba2023}
Ghiba, I.D., B{\^i}rsan, M., Neff, P.: A linear isotropic {Cosserat} shell model including terms up to $\mathcal{O}(h^{5})$. {E}xistence and uniqueness.
\newblock Journal of Elasticity \textbf{154}(1), 579--605 (2023)

\bibitem{ghiba2023essaydeformationmeasuresisotropic}
Ghiba, I.D., Lewintan, P., Sky, A., Neff, P.: An essay on deformation measures in isotropic thin shell theories. {Bending} versus curvature.
\newblock Mathematics and Mechanics of Solids  (2024)

\bibitem{Ghiba}
Ghiba, I.D., Rizzi, G., Madeo, A., Neff, P.: Cosserat micropolar elasticity: classical {Eringen} vs. dislocation form.
\newblock Journal of Mechanics of Materials and Structures \textbf{18}, 93--123 (2023)

\bibitem{GOURGIOTIS2024112700}
Gourgiotis, P., Rizzi, G., Lewintan, P., Bernardini, D., Sky, A., Madeo, A., Neff, P.: Green’s functions for the isotropic planar relaxed micromorphic model — concentrated force and concentrated couple.
\newblock International Journal of Solids and Structures \textbf{292}, 112700 (2024)

\bibitem{Gruttmann}
Gruttmann, F., Sauer, R., Wagner, W.: Shear stresses in prismatic beams with arbitrary cross-sections.
\newblock International Journal for Numerical Methods in Engineering \textbf{45}(7), 865--889 (1999)

\bibitem{Gurtin1975}
Gurtin, M.E., Ian~Murdoch, A.: A continuum theory of elastic material surfaces.
\newblock Archive for Rational Mechanics and Analysis \textbf{57}(4), 291--323 (1975)

\bibitem{GORTHOFER2020109456}
Görthofer, J., Schneider, M., Ospald, F., Hrymak, A., Böhlke, T.: Computational homogenization of sheet molding compound composites based on high fidelity representative volume elements.
\newblock Computational Materials Science \textbf{174}, 109456 (2020)

\bibitem{hale_simple_2018}
Hale, J.S., Brunetti, M., Bordas, S.P., Maurini, C.: Simple and extensible plate and shell finite element models through automatic code generation tools.
\newblock Computers \& Structures \textbf{209}, 163--181 (2018)

\bibitem{FiredrakeUserManual}
Ham, D.A., Kelly, P.H.J., Mitchell, L., Cotter, C.J., Kirby, R.C., Sagiyama, K., Bouziani, N., Vorderwuelbecke, S., Gregory, T.J., Betteridge, J., Shapero, D.R., Nixon-Hill, R.W., Ward, C.J., Farrell, P.E., Brubeck, P.D., Marsden, I., Gibson, T.H., Homolya, M., Sun, T., McRae, A.T.T., Luporini, F., Gregory, A., Lange, M., Funke, S.W., Rathgeber, F., Bercea, G.T., Markall, G.R.: Firedrake User Manual.
\newblock Imperial College London and University of Oxford and Baylor University and University of Washington, first edition edn. (2023)

\bibitem{HANSBO20141}
Hansbo, P., Larson, M.G.: Finite element modeling of a linear membrane shell problem using tangential differential calculus.
\newblock Computer Methods in Applied Mechanics and Engineering \textbf{270}, 1--14 (2014)

\bibitem{HANSBO2022114707}
Hansbo, P., Larson, M.G.: Nitsche’s finite element method for model coupling in elasticity.
\newblock Computer Methods in Applied Mechanics and Engineering \textbf{392}, 114707 (2022)

\bibitem{Hansbo2015}
Hansbo, P., Larson, M.G., Larsson, F.: Tangential differential calculus and the finite element modeling of a large deformation elastic membrane problem.
\newblock Computational Mechanics \textbf{56}(1), 87--95 (2015)

\bibitem{Hansbo2014}
Hansbo, P., Larson, M.G., Larsson, K.: Variational formulation of curved beams in global coordinates.
\newblock Computational Mechanics \textbf{53}(4), 611--623 (2014)

\bibitem{Harsch}
Harsch, J., Sailer, S., Eugster, S.R.: A total {Lagrangian}, objective and intrinsically locking-free {Petrov–Galerkin} {$\mathrm{SE}(3)$} {Cosserat} rod finite element formulation.
\newblock International Journal for Numerical Methods in Engineering \textbf{124}(13), 2965--2994 (2023)

\bibitem{Bence}
Hauck, B., Szekrényes, A.: Enhanced beam and plate finite elements with shear stress continuity for compressible sandwich structures.
\newblock Mathematics and Mechanics of Solids \textbf{29}(7), 1325--1363 (2024)

\bibitem{Hiptmair}
Hiptmair, R., Pauly, D., Schulz, E.: Traces for {H}ilbert complexes.
\newblock Journal of Functional Analysis \textbf{284}(10), 109905 (2023)

\bibitem{HuBordas2020}
Hu, Q., Xia, Y., Natarajan, S., Zilian, A., Hu, P., Bordas, S.P.A.: Isogeometric analysis of thin {Reissner--Mindlin} shells: locking phenomena and {B}-bar method.
\newblock Computational Mechanics \textbf{65}(5), 1323--1341 (2020)

\bibitem{Itskov}
Itskov, M.: Tensor Algebra and Tensor Analysis for Engineers.
\newblock Springer Cham (2019)

\bibitem{Jeong1}
Jeong, J., Neff, P.: Existence, uniqueness and stability in linear {C}osserat elasticity for weakest curvature conditions.
\newblock Mathematics and Mechanics of Solids \textbf{15}(1), 78--95 (2010)

\bibitem{Jeong2}
Jeong, J., Ramézani, H., Münch, I., Neff, P.: A numerical study for linear isotropic {C}osserat elasticity with conformally invariant curvature.
\newblock Zeitschrift für Angewandte Mathematik und Mechanik \textbf{89}(7), 552--569 (2009)

\bibitem{KaiserKirch}
Kaiser, M.W., Fries, T.P.: Curved, linear {Kirchhoff} beams formulated using tangential differential calculus and lagrange multipliers.
\newblock Proc. Appl. Math. Mech. \textbf{22}(1), e202200042 (2023)

\bibitem{Kaiser}
Kaiser, M.W., Fries, T.P.: Simultaneous analysis of continuously embedded {Reissner--Mindlin} shells in {3D} bulk domains.
\newblock International Journal for Numerical Methods in Engineering p. e7495 (2024)

\bibitem{Klarmann2022}
Klarmann, S., Wackerfu{\ss}, J., Klinkel, S.: Coupling {2D} continuum and beam elements: a mixed formulation for avoiding spurious stresses.
\newblock Computational Mechanics \textbf{70}(6), 1145--1166 (2022)

\bibitem{Kuchta}
Kuchta, M.: Assembly of multiscale linear pde operators.
\newblock In: F.J. Vermolen, C.~Vuik (eds.) Numerical Mathematics and Advanced Applications ENUMATH 2019, pp. 641--650. Springer International Publishing, Cham (2021)

\bibitem{Lakes1995EXPERIMENTALMF}
Lakes, R.S.: Experimental methods for study of {C}osserat elastic solids and other generalized elastic continua.
\newblock In: H.B. M\"uhlhaus, editor, Continuum Models for Materials with Microstructure, Wiley, pp. 1--25 (1995).
\newblock \urlprefix\url{http://silver.neep.wisc.edu/∼lakes/CossRv.pdf}

\bibitem{Lewintan2021}
Lewintan, P., M{\"u}ller, S., Neff, P.: Korn inequalities for incompatible tensor fields in three space dimensions with conformally invariant dislocation energy.
\newblock Calculus of Variations and Partial Differential Equations \textbf{60}(4), 150 (2021)

\bibitem{MADEO2016294}
Madeo, A., Ghiba, I.D., Neff, P., Münch, I.: A new view on boundary conditions in the {Grioli–Koiter–Mindlin–Toupin} indeterminate couple stress model.
\newblock European Journal of Mechanics - A/Solids \textbf{59}, 294--322 (2016)

\bibitem{Meier2019}
Meier, C., Popp, A., Wall, W.A.: Geometrically exact finite element formulations for slender beams: {Kirchhoff--Love} theory versus {Simo--Reissner} theory.
\newblock Archives of Computational Methods in Engineering \textbf{26}(1), 163--243 (2019)

\bibitem{Mindlin1964}
Mindlin, R.: Micro-structure in linear elasticity.
\newblock Archive for Rational Mechanics and Analysis \textbf{16}, 51--78 (1964)

\bibitem{SaemDrill}
Mohammadi~Saem, M., Lewintan, P., Neff, P.: On in-plane drill rotations for {C}osserat surfaces.
\newblock Proceedings of the Royal Society A: Mathematical, Physical and Engineering Sciences \textbf{477}(2252), 20210158 (2021)

\bibitem{Munch2011}
M{\"u}nch, I., Neff, P., Wagner, W.: Transversely isotropic material: nonlinear {C}osserat versus classical approach.
\newblock Continuum Mechanics and Thermodynamics \textbf{23}(1), 27--34 (2011)

\bibitem{MULLER201636}
Müller, V., Böhlke, T.: Prediction of effective elastic properties of fiber reinforced composites using fiber orientation tensors.
\newblock Composites Science and Technology \textbf{130}, 36--45 (2016)

\bibitem{NEBEL2023116309}
Nebel, L.J., Sander, O., Bîrsan, M., Neff, P.: A geometrically nonlinear {C}osserat shell model for orientable and non-orientable surfaces: {D}iscretization with geometric finite elements.
\newblock Computer Methods in Applied Mechanics and Engineering \textbf{416}, 116309 (2023)

\bibitem{Ned2}
N{\'e}d{\'e}lec, J.C.: A new family of mixed finite elements in $\mathbb{R}^3$.
\newblock Numerische Mathematik \textbf{50}(1), 57--81 (1986)

\bibitem{Neff2015}
Neff, P., B{\^i}rsan, M., Osterbrink, F.: Existence theorem for geometrically nonlinear {Cosserat} micropolar model under uniform convexity requirements.
\newblock Journal of Elasticity \textbf{121}(1), 119--141 (2015)

\bibitem{NeffWieners}
Neff, P., Che\l{}mi\'{n}ski, K., M\"{u}ller, W., Wieners, C.: A numerical solution method for an infinitesimal elasto-plastic {C}osserat model.
\newblock Mathematical Models and Methods in Applied Sciences \textbf{17}(08), 1211--1239 (2007)

\bibitem{Neff2014}
Neff, P., Ghiba, I.D., Madeo, A., Placidi, L., Rosi, G.: A unifying perspective: the relaxed linear micromorphic continuum.
\newblock Continuum Mechanics and Thermodynamics \textbf{26}(5), 639--681 (2014)

\bibitem{NeffReissner}
Neff, P., Hong, K.I., Jeong, J.: The {R}eissner–{M}indlin plate is the {$\Gamma$}-limit of {C}osserat elasticity.
\newblock Mathematical Models and Methods in Applied Sciences \textbf{20}(09), 1553--1590 (2010)

\bibitem{Jeong3}
Neff, P., Jeong, J.: A new paradigm: the linear isotropic {C}osserat model with conformally invariant curvature energy.
\newblock Zeitschrift für Angewandte Mathematik und Mechanik \textbf{89}(2), 107--122 (2009)

\bibitem{Neff2009}
Neff, P., M{\"u}nch, I.: Simple shear in nonlinear {Cosserat} elasticity: bifurcation and induced microstructure.
\newblock Continuum Mechanics and Thermodynamics \textbf{21}(3), 195--221 (2009)

\bibitem{refId0}
{Neff, Patrizio}, {Münch, Ingo}: Curl bounds grad on {$\mathrm{SO}(3)$}.
\newblock ESAIM: COCV \textbf{14}(1), 148--159 (2008)

\bibitem{MichaelThesis}
Neunteufel, M.: {Mixed Finite Element Methods for Nonlinear Continuum Mechanics and Shells}.
\newblock Ph.D. thesis, Technische Universität Wien (2021).
\newblock \urlprefix\url{http://hdl.handle.net/20.500.12708/17043}

\bibitem{NS21}
Neunteufel, M., Sch\"oberl, J.: Avoiding membrane locking with {R}egge interpolation.
\newblock Computer Methods in Applied Mechanics and Engineering \textbf{373}, 113524 (2021)

\bibitem{neunteufel_hellanherrmannjohnson_2019}
Neunteufel, M., Schöberl, J.: The {Hellan}–{Herrmann}–{Johnson} method for nonlinear shells.
\newblock Computers \& Structures \textbf{225}, 106109 (2019)

\bibitem{Nguyen2021}
Nguyen, D.T.A., Li, L., Ji, H.: Stable and accurate numerical methods for generalized {Kirchhoff--Love} plates.
\newblock Journal of Engineering Mathematics \textbf{130}(1), 6 (2021)

\bibitem{nguyena2013nitschesmethodmethodmixed}
Nguyena, V.P., Kerfriden, P., Clausb, S., Bordas, S.P.A.: Nitsche's method method for mixed dimensional analysis: conforming and non-conforming continuum-beam and continuum-plate coupling.
\newblock arXiv  (2013).
\newblock \urlprefix\url{https://arxiv.org/abs/1308.2910}

\bibitem{PaulyDeRham}
Pauly, D., Schomburg, M.: Hilbert complexes with mixed boundary conditions -- {Part} 1: de {Rham} complex.
\newblock Mathematical Methods in the Applied Sciences \textbf{45}(5), 2465--2507 (2022)

\bibitem{pechstein_tdnns_2017}
Pechstein, A.S., Schöberl, J.: The {TDNNS} method for {Reissner}–{Mindlin} plates.
\newblock Numerische Mathematik \textbf{137}(3), 713--740 (2017)

\bibitem{Pi}
Pi, Y.L., Bradford, M.A., Uy, B.: A spatially curved-beam element with warping and {Wagner} effects.
\newblock International Journal for Numerical Methods in Engineering \textbf{63}(9), 1342--1369 (2005)

\bibitem{Alessio}
Rubino, A., Accornero, F., Carpinteri, A.: Fracture mechanics approach to minimum reinforcement design of fibre-reinforced and hybrid-reinforced concrete beams.
\newblock International Journal of Damage Mechanics \textbf{0}(0), 10567895241245865 (0)

\bibitem{Sch1997}
Sch{\"o}berl, J.: {NETGEN} an advancing front 2{D}/3{D}-mesh generator based on abstract rules.
\newblock Computing and Visualization in Science \textbf{1}(1), 41--52 (1997)

\bibitem{Sch2014}
Sch{\"o}berl, J.: C++ 11 implementation of finite elements in {NGS}olve.
\newblock Institute for Analysis and Scientific Computing, Vienna University of Technology  (2014).
\newblock \urlprefix\url{https://www.asc.tuwien.ac.at/~schoeberl/wiki/publications/ngs-cpp11.pdf}

\bibitem{Schollhammer2019}
Sch{\"o}llhammer, D., Fries, T.P.: {Kirchhoff--Love} shell theory based on tangential differential calculus.
\newblock Computational Mechanics \textbf{64}(1), 113--131 (2019)

\bibitem{SCHOLLHAMMER2019172}
Schöllhammer, D., Fries, T.: {Reissner–Mindlin} shell theory based on tangential differential calculus.
\newblock Computer Methods in Applied Mechanics and Engineering \textbf{352}, 172--188 (2019)

\bibitem{Sevilla}
Sevilla, R., Rees, L., Hassan, O.: The generation of triangular meshes for {NURBS}-enhanced {FEM}.
\newblock International Journal for Numerical Methods in Engineering \textbf{108}(8), 941--968 (2016)

\bibitem{Shim2002}
Shim, K.W., Monaghan, D.J., Armstrong, C.G.: Mixed dimensional coupling in finite element stress analysis.
\newblock Engineering with Computers \textbf{18}(3), 241--252 (2002)

\bibitem{Shirani}
Shirani, M., Steigmann, D.J., Neff, P.: {The Legendre–Hadamard condition in Cosserat elasticity theory}.
\newblock The Quarterly Journal of Mechanics and Applied Mathematics \textbf{73}(4), 293--303 (2020)

\bibitem{SIMO1991371}
Simo, J., Vu-Quoc, L.: A geometrically-exact rod model incorporating shear and torsion-warping deformation.
\newblock International Journal of Solids and Structures \textbf{27}(3), 371--393 (1991)

\bibitem{sky_polytopal_2022}
Sky, A., Muench, I.: Polytopal templates for semi-continuous vectorial finite elements of arbitrary order on triangulations and tetrahedralizations.
\newblock Finite Elements in Analysis and Design \textbf{236}, 104155 (2024)

\bibitem{sky_higher_2023}
Sky, A., Muench, I., Rizzi, G., Neff, P.: Higher order {Bernstein–Bézier} and {Nédélec} finite elements for the relaxed micromorphic model.
\newblock Journal of Computational and Applied Mathematics \textbf{438}, 115568 (2024)

\bibitem{sky2023reissnermindlin}
Sky, A., Neunteufel, M., Hale, J.S., Zilian, A.: A {R}eissner–{M}indlin plate formulation using symmetric {H}u-{Z}hang elements via polytopal transformations.
\newblock Computer Methods in Applied Mechanics and Engineering \textbf{416}, 116291 (2023)

\bibitem{SkyFormulae}
Sky, A., Neunteufel, M., Hale, J.S., Zilian, A.: Formulae and transformations for simplicial tensorial finite elements via polytopal templates.
\newblock arXiv  (2024).
\newblock \urlprefix\url{https://arxiv.org/abs/2405.10402}

\bibitem{SKYNOVEL}
Sky, A., Neunteufel, M., Lewintan, P., Zilian, A., Neff, P.: Novel {H(symCurl)-conforming} finite elements for the relaxed micromorphic sequence.
\newblock Computer Methods in Applied Mechanics and Engineering \textbf{418}, 116494 (2024)

\bibitem{SKY2022115298}
Sky, A., Neunteufel, M., Muench, I., Schöberl, J., Neff, P.: Primal and mixed finite element formulations for the relaxed micromorphic model.
\newblock Computer Methods in Applied Mechanics and Engineering \textbf{399}, 115298 (2022)

\bibitem{sky_hybrid_2021}
Sky, A., Neunteufel, M., Münch, I., Schöberl, J., Neff, P.: A hybrid $\mathit{H}^1 \times \mathit{H}(\mathrm{curl})$ finite element formulation for a relaxed micromorphic continuum model of antiplane shear.
\newblock Computational Mechanics \textbf{68}(1), 1--24 (2021)

\bibitem{Song}
Song, H., Hodges, D.: Rigorous joining of advanced reduced-dimensional beam models to {2-D} finite element models.
\newblock 51st AIAA/ASME/ASCE/AHS/ASC Structures, Structural Dynamics, and Materials and Co-located Conferences. American Institute of Aeronautics and Astronautics (2010)

\bibitem{Steinbrecher2020}
Steinbrecher, I., Mayr, M., Grill, M.J., Kremheller, J., Meier, C., Popp, A.: A mortar-type finite element approach for embedding {1D} beams into {3D} solid volumes.
\newblock Computational Mechanics \textbf{66}(6), 1377--1398 (2020)

\bibitem{Steinbrecher2022}
Steinbrecher, I., Popp, A., Meier, C.: Consistent coupling of positions and rotations for embedding {1D} {Cosserat} beams into {3D} solid volumes.
\newblock Computational Mechanics \textbf{69}(3), 701--732 (2022)

\bibitem{VO2022114883}
Vo, D., Nanakorn, P., Bui, T.Q., Rungamornrat, J.: On invariance of spatial isogeometric {Timoshenko–Ehrenfest} beam formulations for static analysis.
\newblock Computer Methods in Applied Mechanics and Engineering \textbf{394}, 114883 (2022)

\bibitem{WEEGER2017100}
Weeger, O., Yeung, S.K., Dunn, M.L.: Isogeometric collocation methods for {Cosserat} rods and rod structures.
\newblock Computer Methods in Applied Mechanics and Engineering \textbf{316}, 100--122 (2017).
\newblock Special Issue on Isogeometric Analysis: Progress and Challenges

\bibitem{Yamamoto2019}
Yamamoto, T., Yamada, T., Matsui, K.: Numerical procedure to couple shell to solid elements by using {N}itsche's method.
\newblock Computational Mechanics \textbf{63}(1), 69--98 (2019)

\bibitem{YANG2024105545}
Yang, Y., Lu, P., Liu, Z., Dong, L., Lin, J., Yang, T., Ren, Q., Wu, C.: Effect of steel fibre with different orientations on mechanical properties of {3D}-printed steel-fibre reinforced concrete: {M}esoscale finite element analysis.
\newblock Cement and Concrete Composites \textbf{150}, 105545 (2024)

\bibitem{Yuan}
Yuan, J., Mu, Z., Elishakoff, I.: Novel modification to the {Timoshenko–Ehrenfest} theory for inhomogeneous and nonuniform beams.
\newblock AIAA Journal \textbf{58}(2), 939--948 (2020)

\bibitem{ZOU2024103653}
Zou, X., Lo, S.B., Sevilla, R., Hassan, O., Morgan, K.: The generation of {3D} surface meshes for {NURBS-Enhanced FEM}.
\newblock Computer-Aided Design \textbf{168}, 103653 (2024)

\bibitem{Zou1}
Zou, X., Lo, S.B., Sevilla, R., Hassan, O., Morgan, K.: Towards a volume mesh generator tailored for {NEFEM}.
\newblock In: E.~Ruiz-Giron{\'e}s, R.~Sevilla, D.~Moxey (eds.) SIAM International Meshing Roundtable 2023, pp. 397--418. Springer Nature Switzerland, Cham (2024)

\end{thebibliography}

\normalsize
\appendix

\section{Tensor calculus}\label{ap:tensorcal}
Let $\Omega \subset \R^3$ define some reference three-dimensional body and $\vol \subset \R^3$ a physical three-dimensional body, we define the mapping 
\begin{align}
    &\vb{x} = \vb{x}(\xi,\eta,\zeta) \, , && \vb{x}:\Omega \to \vol \, .
\end{align}
The covariant tangent vectors in the physical body read
\begin{align}
    \vb{g}_i = \partial_i^\xi \vb{x} = \pder{}{\xi^i} \vb{x} = \pder{x^j}{\xi^i} \vb{e}_j = \vb{x}_{,i} \, .
\end{align}
Let the mapping $\vb{x}$ be invertible $\bm{\xi} = \bm{\xi}(\vb{x}) = \vb{x}^{-1}(\bm{\xi})$, using the chain rule we can directly deduce
\begin{align}
    \pder{x^j}{\xi^i} \pder{\xi^k}{x^l} \con{\vb{e}_j}{\vb{e}_l} = \pder{x^j}{\xi^i} \pder{\xi^k}{x^l} \delta_{jl} = \pder{x^j}{\xi^i} \pder{\xi^k}{x^j}
    =  \pder{\xi^k}{x^j}\pder{x^j}{\xi^i} = \pder{\xi^k}{\xi^i} = \delta_{ki} \, .
\end{align}
implying the definition of so called contravariant vectors
\begin{align}
   &\vb{g}^j = \pder{\xi^j}{x^i} \vb{e}_i \,, && \con{\vb{g}_i}{\vb{g}^j} = \delta_{ij} \, . 
\end{align}
Inversely, we can retrieve the Cartesian basis vectors form the co- and contravariant basis
\begin{align}
    \vb{e}_j = \pder{\xi^i}{x^j} \vb{g}_i  \, , && \vb{e}_i = \pder{x^i}{\xi^j} \vb{g}^j \, . 
\end{align}
The covariant vectors define the metric tensor of the three-dimensional body
\begin{align}
    &\bm{G} = g_{ij} \vb{e}_i \otimes \vb{e}_j = \con{\vb{g}_i}{\vb{g}_j} \, \vb{e}_i \otimes \vb{e}_j = \pder{x^k}{\xi^i} \pder{x^k}{\xi^j} \vb{e}_i \otimes \vb{e}_j  \, .
\end{align}
We observe 
\begin{align}
    (g_{ij} \vb{e}_i \otimes \vb{e}_j) (\con{\vb{g}^k}{\vb{g}^l} \, \vb{e}_k \otimes \vb{e}_l) = (g_{ij} \vb{e}_i \otimes \vb{e}_j) (g^{kl} \vb{e}_k \otimes \vb{e}_l) = g_{ij} g^{kl} \delta_{jk} \vb{e}_i \otimes \vb{e}_l = g_{ik} g^{kl} \vb{e}_i \otimes \vb{e}_l = \delta_{il}  \vb{e}_i \otimes \vb{e}_l \, ,
\end{align}
since 
\begin{align}
    g_{ik} g^{kl} =  \pder{x^q}{\xi^i} \pder{x^q}{\xi^k} \pder{\xi^k}{x^r} \pder{\xi^l}{x^r} = \pder{x^q}{\xi^i} \pder{x^q}{x^r} \pder{\xi^l}{x^r} = \delta_{qr} \pder{x^q}{\xi^i} \pder{\xi^l}{x^r} = \pder{x^r}{\xi^i} \pder{\xi^l}{x^r} = \pder{\xi^l}{\xi^i} = \delta_{li} \, .
\end{align}
Consequently, the inverse of the metric tensor is
\begin{align}
    \bm{G}^{-1} = g^{ij} \vb{e}_i \otimes \vb{e}_j = \con{\vb{g}^i}{\vb{g}^j} \, \vb{e}_i \otimes \vb{e}_j \, .
\end{align}
A covariant vector is mapped to a contravariant vector via
\begin{align}
    \vb{g}_i g^{ij} = \pder{x^k}{\xi^i} \vb{e}_k   \pder{\xi^i}{x^l} \pder{\xi^j}{x^l} = \vb{e}_k \pder{x^k}{x^l} \pder{\xi^j}{x^l} = \vb{e}_k \delta_{kl} \pder{\xi^j}{x^l} = \vb{e}_k  \pder{\xi^j}{x^k} = \vb{g}^j \, ,  \label{eq:changeind1}
\end{align}
and vice versa
\begin{align}
    \vb{g}_j = \vb{g}^i g_{ij} \, . \label{eq:changeind2}
\end{align}
With this basic machinery in place, we can now express second order derivatives.
Second order derivatives of the mapping can be expressed using Christoffel symbols
\begin{align}
    \partial_j^\xi \vb{g}_i = \pder{}{\xi^j} \vb{g}_i = \vb{g}_{i,j} = \Gamma_{ij}^{\phantom{ij}k} \vb{g}_k = \Gamma_{ijk} \vb{g}^k \, .
\end{align}
The components of the Christoffel symbols $\Gamma_{ij}^{\phantom{ij}k}$ and $\Gamma_{ijk}$ are identified using the orthogonality of the co- and contravariant vectors
\begin{align}
    &\con{\vb{g}_{i,j}}{\vb{g}^l} = \Gamma_{ij}^{\phantom{ij}k} \con{\vb{g}_k}{\vb{g}^l} = \Gamma_{ij}^{\phantom{ij}k} \delta_{k}^{\phantom{k}l} = \Gamma_{ij}^{\phantom{ij}l} \, , && \con{\vb{g}_{i,j}}{\vb{g}_l} = \Gamma_{ijk} \con{\vb{g}^k}{\vb{g}_l} = \Gamma_{ijk} \delta^{k}_{\phantom{k}l} = \Gamma_{ijl} \, .
\end{align}

\section{Tangential differential calculus}\label{ap:tdc}
The identity tensor reads
\begin{align}
    \one = \delta_{ij} \vb{e}_i \otimes \vb{e}_j = \vb{e}_i \otimes \vb{e}_i \, .
\end{align}
It can also be expressed as a mixed tensor of co- and contravariant basis vectors
\begin{align}
    \one = \vb{e}_i \otimes \vb{e}_i = \pder{\xi^j}{x^i} \vb{g}_j \otimes \pder{x^i}{\xi^k} \vb{g}^k = \pder{\xi^j}{\xi^k} \vb{g}_j \otimes \vb{g}^k = \delta_{jk} \vb{g}_j \otimes \vb{g}^k = \vb{g}_j \otimes \vb{g}^j \, ,  
\end{align}
and analogously as $\one = \vb{g}^j \otimes \vb{g}_j$. 
If we define some hyper-surface $\surf \subset \R^3$ by a mapping $\vb{r}:\omega \subset \R^2 \to \surf$, whose tangent vectors are $\vb{g}_1$ and $\vb{g}_2$, then its normal unit vector $\vb{g}_3$ can be defined as
\begin{align}
    \vb{g}_3 = \dfrac{\vb{g}_1 \times \vb{g}_2}{\norm{\vb{g}_1 \times \vb{g}_2}} \, .
\end{align}
We observe that
\begin{align}
    &\con{\vb{g}_3}{\vb{g}_3} = 1 \, , && \con{\vb{g}_\alpha}{\vb{g}_3} = 0 \, , 
\end{align}
implying that the co- and contravariant normal vectors are the same $\vb{g}_3 = \vb{g}^3$. Since $\vb{g}_3$ is the unit normal vector to the surface, we can define a corresponding tangential projection operator
\begin{align}
    &\Pt = \one - \vb{n} \otimes \vb{n} = \one - \vb{g}_3 \otimes \vb{g}^3 =  \vb{g}_i \otimes \vb{g}^i - \vb{g}_3 \otimes \vb{g}^3 =  \vb{g}_\alpha \otimes \vb{g}^\alpha \, , && \Pt \, \vb{v} =  (\vb{g}_\alpha \otimes \vb{g}^\alpha) \,  v^i \vb{g}_i = v^\alpha \vb{g}_\alpha \, . 
\end{align}
In other words, the tensor eliminates any non-tangential components. Analogously, we can define the normal projection operator
\begin{align}
    &&\Pn = \one - \Pt = \vb{g}_3 \otimes \vb{g}^3 = \vb{n} \otimes \vb{n} \, , && \Pn \, \vb{v} =  (\vb{g}_3 \otimes \vb{g}^3) \,  v^i \vb{g}_i = v^3 \vb{g}_3 \, . 
\end{align}
Now, let some function $\lambda$ depend only on two parameters of a reference surface $(\lambda \circ \vb{r})(\xi,\eta)$, then its gradient reads
\begin{align}
    \nabla \lambda = \vb{e}^i \partial_i^x \lambda = \vb{e}_i \pder{}{x^i} \lambda = \vb{e}_i \pder{}{\xi^\alpha} \pder{\xi^\alpha}{x^i} \lambda = \vb{e}_i \pder{\xi^\alpha}{x^i} \pder{}{\xi^\alpha}  \lambda = \vb{g}^\alpha \partial_\alpha^\xi \lambda = \lambda_{,\alpha} \, \vb{g}^\alpha \, .
\end{align}
Evidently, $\nabla \lambda : \omega \to \tansurf \surf$ is a tangential vector in this case such that 
\begin{align}
    \Pt \, \nabla \lambda = (\vb{g}^\beta \otimes \vb{g}_\beta) \lambda_{,\alpha} \, \vb{g}^\alpha =  \lambda_{,\beta} \, \vb{g}^\beta \, ,
\end{align}
leaves the vector unchanged. We call the projected gradient the tangential gradient and define it as
\begin{align}
    \Pt \, \nabla \lambda = \nabla_t \lambda \, . 
\end{align}
If $\lambda$ is a function of a three-dimensional system $(\lambda \circ \vb{x})(\xi,\eta,\zeta)$, then the tangential gradient eliminates the out-of-plane component
\begin{align}
    \nabla_t \lambda = \Pt \, \vb{e}^i \partial_i^x \lambda = \Pt \, \vb{g}^i \partial_i^\xi \lambda = \vb{g}^\alpha \partial_\alpha^\xi \lambda =  \lambda_{,\alpha} \, \vb{g}^\alpha  \, .
\end{align}
Analogously, we can define the same operator for vectors
\begin{align}
    \D_t \vb{v} =  (\D \vb{v}) \Pt = \vb{v}_{,\alpha} \otimes \vb{g}^\alpha \, ,
\end{align}
by applying the projection row-wise.
If we restrict also the vectorial basis on the left of the tensor to the tangential plane we retrieve the so-called covariant gradient
\begin{align}
    \Dcov \vb{v} = \Pt \D_t \vb{v} = \Pt (\D \vb{v}) \Pt \, .
\end{align}
As the name suggests, the latter simply expresses derivatives within the tangential coordinate system. 
With the vectorial gradient defined, the tangential divergence is naturally 
\begin{align}
    \di_t \vb{v} = \tr (\D_t \vb{v}) = \con{(\D \vb{v})\Pt}{\one} = \con{\D \vb{v}}{\Pt} = \con{\vb{v}_{,i} \otimes \vb{g}^i}{\vb{g}^\alpha\otimes\vb{g}_\alpha} = \con{\vb{v}_{,i}}{\vb{g}^\alpha} \delta_{i\alpha} = \con{\vb{v}_{,\alpha}}{\vb{g}^\alpha}  \, . 
    \label{eq:divt}
\end{align}
We note that the trace of the covariant gradient yields the same result
\begin{align}
     \tr (\Dcov \vb{v}) = \con{\Pt(\D \vb{v})\Pt}{\one} = \con{\D \vb{v}}{\Pt} = \di_t \vb{v} \, , 
\end{align}
seeing as $\Pt \Pt = \Pt$ by its projection property. Accordingly, the tensorial tangential divergence reads
\begin{align}
    \Di_t \bm{T} = (\bm{T}_{,\alpha}) \,  \vb{g}^\alpha = (\D \bm{T}) \Pt \, ,
\end{align}
where the latter implies a double-contraction. The tangential gradient also allows to define the surface curl operator
\begin{align}
    \curl_t \vb{v} = \con{\D_t \vb{v}}{\Anti \vb{n}} = \con{\vb{v}_{,\alpha}\otimes \vb{g}^\alpha}{\vb{n}\times\one} = \con{\vb{v}_{,\alpha}\otimes \vb{g}^\alpha}{\vb{n}\times\vb{g}^i\otimes\vb{g}_i} = \con{\vb{v}_{,\alpha}}{\vb{n}\times\vb{g}^i}\delta_{\alpha i }  = \con{\vb{g}^\alpha\times\vb{v}_{,\alpha}}{\vb{n}} \, ,
\end{align}
where we applied the circular shift in the last step. If the vector is a function of the plane $\vb{v} = (\vb{v} \circ \vb{r})(\xi,\eta)$, then the definition coincides with the analogous formula of the full gradient   
\begin{align}
    \curl_t \vb{v} = \con{\D \vb{v}}{\Anti \vb{n}} = \con{\vb{v}_{,j}\otimes \vb{g}^j}{\vb{n}\times\one}  = \con{\vb{v}_{,j}}{\vb{n}\times\vb{g}^j} = \con{\vb{v}_{,\alpha}}{\vb{n}\times\vb{g}^\alpha} = \con{\vb{g}^\alpha\times\vb{v}_{,\alpha}}{\vb{n}} \, ,
\end{align}
since $\vb{n} \times \vb{g}^3 = \vb{n} \times \vb{n} = 0$. Finally, the tangential curl can also be related to the covariant skew-symmetric gradient of a vector field $\vb{v} = (\vb{v} \circ \vb{r})(\xi,\eta)$ via
\begin{align}
    \skw \Dcov \vb{v} = \Pt (\skw \D \vb{v}) \Pt = \dfrac{1}{2} \Pt (\Anti \curl \vb{v}) \Pt = \dfrac{1}{2} (\curl_t \vb{v})(\Anti \vb{n}) \, ,   
    \label{eq:curlDcov}
\end{align}
where any in-plane vectors are eliminated due to 
\begin{align}
    \Pt (\Anti \vb{g}_\alpha ) \Pt = \Pt [\vb{g}_\alpha \times (\vb{g}_\beta \otimes \vb{g}^\beta)] = \dfrac{1}{\norm{\vb{g}_1 \times \vb{g}_2}} \Pt (\vb{n} \otimes \varepsilon_{\alpha \beta \gamma}\vb{g}^\gamma) = 0 \, . 
\end{align}

\section{The Weingarten curvature tensor}\label{ap:weingarten}
Let a hyper-surface $\surf \subset \R^3$ in a three dimensional space $\R^3$ be mapped from some flat surface $\omega \subset \R^2$ via $\vb{r}:\omega \subset \R^2 \to \surf$, the tangent vectors on the hyper-surface read
\begin{align}
    &\vb{t}_1 = \vb{g}_1 = \pder{}{\xi} \vb{r} = \vb{r}_{,\xi} \, , && \vb{t}_2 = \vb{g}_2 = \pder{}{\eta} \vb{r} = \vb{r}_{,\eta} \, .
\end{align}
Using the latter two one can define a unit normal vector on the hyper-surface
\begin{align}
    &\vb{n} = \dfrac{\vb{t}_1  \times \vb{t}_2}{\norm{\vb{t}_1  \times \vb{t}_2}} \, , && \norm{\vb{n}} = 1 \, .
    \label{eq:surfnormal}
\end{align}
Clearly, the tangential vectors $\vb{t}_\alpha$ are orthogonal to $\vb{n}$ by the very construction $\vb{t}_\alpha \perp \vb{n}$, such that
\begin{align}
    &\con{\vb{t}_\alpha}{\vb{n}} = 0 \, ,&& \partial_\beta^\xi \con{\vb{t}_\alpha}{\vb{n}} = \pder{}{\xi^\beta} \con{\vb{t}_\alpha}{\vb{n}} = \con{\vb{t}_{\alpha,\beta}}{\vb{n}} + \con{\vb{t}_{\alpha}}{\vb{n_{,\beta}}} = 0 \, , 
\end{align}
implying the equality
\begin{align}
    \con{\vb{t}_{\alpha,\beta}}{\vb{n}} = - \con{\vb{t}_{\alpha}}{\vb{n_{,\beta}}} \, .
    \label{eq:wid1}
\end{align}
Further, since $\vb{n}$ is a unit vector there holds
\begin{align}
     &\con{\vb{n}}{\vb{n}} = 1 \, , && \partial_\beta^\xi \con{\vb{n}}{\vb{n}} = \pder{}{\xi^\beta} \con{\vb{n}}{\vb{n}} = 2 \con{\vb{n}_{,\beta}}{\vb{n}} = 0 \,,
     \label{eq:wid2}
\end{align}
implying that infinitesimal changes in $\vb{n}$ with respect to $\{\xi,\eta\}$ are orthogonal to it $\vb{n}_{,\beta} \perp \vb{n}$.
Since $\vb{n}$ is a unit vector normal to the hyper-surface $\surf$, any change in its orientation is a measure of curvature. 
Using the Christoffel symbols we express derivatives of the tangent vectors as
\begin{align}
    \vb{t}_{\alpha,\beta} = \Gamma_{\alpha\beta}^{\phantom{\alpha\beta}i} \vb{g}_i = \Gamma_{\alpha\beta}^{\phantom{\alpha\beta}\gamma} \vb{t}_\gamma + \Gamma_{\alpha\beta}^{\phantom{\alpha\beta}3} \vb{g}_3 =  \Gamma_{\alpha\beta}^{\phantom{\alpha\beta}\gamma} \vb{t}_\gamma + W_{\alpha\beta} \vb{n} \, , = \Gamma_{\alpha\beta\gamma} \vb{t}^\gamma + W_{\alpha\beta} \vb{n} \, , 
\end{align}
where we define the covariant components of the so called Weingarten map as
\begin{align}
    W_{\alpha\beta} = \Gamma_{\alpha\beta}^{\phantom{\alpha\beta}3} = \con{\vb{t}_{\alpha,\beta} }{\vb{n}} = \con{\vb{t}_{\alpha,\beta} }{\vb{g}^3}= \con{\vb{t}_{\alpha,\beta} }{\vb{g}_3} = \Gamma_{\alpha\beta3} \, .
\end{align}
Doing the same for derivatives of the normal vector we find
\begin{align}
    \vb{n}_{,\alpha} = \Gamma_{3\alpha}^{\phantom{3\alpha}k} \vb{g}_k = \Gamma_{3\alpha}^{\phantom{3\alpha}\beta} \vb{t}_\beta + \Gamma_{3\alpha}^{\phantom{3\alpha}3} \vb{n}  \, .
\end{align}
Now, using \cref{eq:wid1} and \cref{eq:wid2} we find
\begin{align}
    \Gamma_{3\alpha\beta} = \con{\vb{n}_{,\alpha}}{\vb{t}_\beta} = -\con{\vb{t}_{\beta,\alpha}}{\vb{n}} = -W_{\alpha \beta} \, , && \Gamma_{3\alpha}^{\phantom{3\alpha}3} = \con{\vb{n}_{,\alpha}}{\vb{n}} = 0 \, , 
\end{align}
such that any infinitesimal change in the unit normal vector of the surface is captured by the Weingarten map
\begin{align}
    \vb{n}_{,\alpha} = -W_{\alpha \beta}\vb{g}^\beta = -W_\alpha^{\phantom{\alpha}\beta} \vb{t}_\beta \, , 
    \label{eq:dern}
\end{align}
where we used \cref{eq:changeind1} and \cref{eq:changeind2} to switch to a mixed variant definition. Clearly, the components of the Weingarten tensor are curvature measures, such that the tensor itself is known as the curvature tensor. Its components can be identified using  
\begin{align}
    W_{\alpha}^{\phantom{\alpha}\beta} = -\con{\vb{n}_{,\alpha}}{\vb{g}^\beta} \, .
\end{align}
Consequently, the tensor itself can be defined via
\begin{align}
    \bm{W} = -\D \vb{n} = -\vb{n}_{,\alpha} \otimes \vb{g}^\alpha =  W_{\alpha}^{\phantom{\alpha}\beta} \vb{t}_\beta \otimes \vb{g}^\alpha \, .
\end{align}
It is clear that the tensor is tangential, such we can also write
\begin{align}
    \bm{W} = -\D_t \vb{n} \, ,
\end{align}
using tangential differential calculus. Let $\kappa_1$ and $\kappa_2$ be the eigenvalues of the tensor, its invariants read
\begin{align}
    &H = \dfrac{1}{2} \tr \bm{W} = \dfrac{1}{2} W_\alpha^{\phantom{\alpha}\alpha} = \dfrac{1}{2} (W_1^{\phantom{1}1} + W_2^{\phantom{2}2}) = \dfrac{1}{2}(\kappa_1 + \kappa_2) \, , && K = \det \bm{W} = W_1^{\phantom{1}1} W_2^{\phantom{2}2} - W_1^{\phantom{1}2} W_2^{\phantom{1}1} = \kappa_1 \kappa_2   \, , 
\end{align}
representing the mean- and Gaussian curvature measures, respectively.

\section{The shell-shifter}\label{ap:shellshifter}
Let the middle hyper-surface of the shell $\surf$ be mapped from some flat reference surface $\omega$
\begin{align}
    &\vb{r} = \vb{r}(\xi,\eta) \, , && \vb{r}:\omega \subset \R^2 \to \surf \subset \R^3 \, , 
\end{align}
one can define the complete volume of the shell via
\begin{align}
    &\vb{x} = \vb{x}(\xi,\eta,\zeta) \, , \vb{x} = \vb{r} + \zeta \vb{n} \,, && \vb{x}:\Omega \subset \R^3 \to \vol \subset \R^3 \, ,
    \label{eq:shellvolpar}
\end{align}
where $\zeta \in [-h/2, h/2]$ is the thickness parameter of the shell, and the normal vector $\vb{n}$ is defined using the middle surface as per \cref{eq:surfnormal}. The volume of the shell is therefore given by
\begin{align}
    \vol = \surf \times [-h/2, h/2] \subset \R^3 \, .
\end{align}
The infinitesimal tangent vectors of the shell read
\begin{align}
    &\pder{\vb{x}}{\xi}\dd \xi = \vb{g}_1 \dd \xi = (\vb{t}_1 + \zeta \vb{n}_{,\xi}) \dd \xi \, , && \pder{\vb{x}}{\eta}\dd \eta = \vb{g}_2 \dd \eta =  (\vb{t}_2 + \zeta \vb{n}_{,\eta}) \dd \eta \, , && \pder{\vb{x}}{\zeta}\dd \zeta = \vb{g}_3 \dd \zeta = \vb{n} \dd \zeta \, .
\end{align}
As the surface is parameterised by $\{\xi,\eta\}$, its infinitesimal tangent vectors reads
\begin{align}
    &\pder{\vb{r}}{\xi}\dd \xi = \vb{t}_1 \dd \xi \, , && \pder{\vb{r}}{\eta}\dd \xi = \vb{t}_2 \dd \eta \, . 
\end{align}
Consequently, an infinitesimal surface element of the middle surface reads
\begin{align}
    \dd \surf = \norm{\vb{t}_1 \dd \xi \times \vb{t}_2 \dd \eta} = \norm{\vb{t}_1  \times \vb{t}_2 } \dd \xi \dd \eta = \norm{\vb{t}_1  \times \vb{t}_2 } \dd \omega \, . 
\end{align}
We find the infinitesimal volume element of the shell using the triple vector product
\begin{align}
    \dd \vol = \con{\vb{x}_{,\xi}\dd \xi \times \vb{x}_{,\eta} \dd \eta}{\vb{x}_{,\zeta} \dd \zeta} = \con{\vb{x}_{,\xi} \times \vb{x}_{,\eta}}{\vb{x}_{,\zeta} } \dd \zeta \dd\omega  = \con{(\vb{t}_1 + \zeta \vb{n}_{,\xi}) \times (\vb{t}_2 + \zeta \vb{n}_{,\eta})}{\vb{n}} \dd \zeta \dd\omega  \, .
\end{align}
The term can be decomposed into the additive parts
\begin{align}
    &\con{\vb{t}_1 \times \vb{t}_2}{\vb{n}} \dd \zeta \dd\omega  \, , && \con{\vb{t}_1 \times \zeta\vb{n}_{,\eta} + \zeta\vb{n}_{,\xi} \times \vb{t}_2 }{\vb{n}} \dd \zeta \dd\omega \, , && \con{\zeta\vb{n}_{,\xi} \times \zeta\vb{n}_{,\eta}}{\vb{n}} \dd \zeta \dd\omega \, . 
\end{align}
The first term is clearly
\begin{align}
    \con{\vb{t}_1 \times \vb{t}_2}{\vb{n}} \dd \zeta \dd\omega = \dd \zeta \dd \surf \, ,
\end{align}
due to $\vb{n} \parallel \vb{t}_1 \times \vb{t}_2$ and $\norm{\vb{n}}=1$. Expanding the second part while applying the curvature tensor $\bm{W} = -\D_t \vb{n}$ to express derivatives of the normal vector we find
\begin{align}
    &\vb{t}_1 \times \zeta\vb{n}_{,\eta} = -\zeta \vb{t}_1 \times W_2^{\phantom{2}2}\vb{t}_2 \, , && \zeta\vb{n}_{,\xi} \times \vb{t}_2 = -\zeta W_1^{\phantom{1}1} \vb{t}_1 \times \vb{t}_2 \, , 
\end{align}
yielding together
\begin{align}
    \con{-\zeta \vb{t}_1 \times W_2^{\phantom{2}2}\vb{t}_2 -\zeta W_1^{\phantom{1}1} \vb{t}_1 \times \vb{t}_2}{\vb{n}} \dd \zeta \dd\omega =  -\zeta  ( W_1^{\phantom{1}1} + W_2^{\phantom{2}2}) \con{\vb{t}_1 \times\vb{t}_2}{\vb{n}} \dd \zeta \dd\omega = -2H \, \zeta \, \dd \zeta  \dd \surf \, . 
\end{align}
Using the curvature tensor $\bm{W}$ again for the third term we find
\begin{align}
    \zeta\vb{n}_{,\xi} \times \zeta\vb{n}_{,\eta} = \zeta^2 (W_1^{\phantom{1}1}\vb{t}_1 + W_1^{\phantom{1}2}\vb{t}_2) \times (W_2^{\phantom{2}1}\vb{t}_1 + W_2^{\phantom{2}2}\vb{t}_2) = \zeta^2(\vb{t}_1 \times \vb{t}_2)(W_1^{\phantom{1}1}W_2^{\phantom{2}2} - W_1^{\phantom{1}2}W_2^{\phantom{1}1}) = K \, \zeta^2(\vb{t}_1 \times \vb{t}_2) \, , 
\end{align}
such the term reads
\begin{align}
    \con{\zeta\vb{n}_{,\xi} \times \zeta\vb{n}_{,\eta}}{\vb{n}} \dd \zeta \dd\omega = \kappa \, \zeta^2 \con{\vb{t}_1 \times \vb{t}_2}{\vb{n}} \dd \zeta \dd\omega = K \, \zeta^2 \, \dd \zeta  \dd \surf \, .
\end{align}
Putting it all together, we find Steiner's formula for the change in volume
\begin{align}
    \dd \vol = (1 -2H \, \zeta + K \, \zeta^2) \, \dd \zeta  \dd \surf \, ,
\end{align}
where $1 -2H \, \zeta + K \, \zeta^2$ is called the shell-shifter.

\section{Tangential gradients on shells} \label{ap:tangradshell}
In general, for fields on the shell $\lambda = (\lambda \circ \vb{x})(\xi,\eta,\zeta)= (\lambda \circ \vb{x})(\bm{\xi})$ it is possible to determine corresponding gradients given by tangential differential calculus. In other words, the gradients can be characterised by the parametrisation of the middle surface. We observe that the covariant basis of the parametrisation $\vb{x}(\xi,\eta,\zeta)$ from \cref{eq:shellvolpar} can be expanded as
\begin{align}
    \vb{g}_1 = \vb{t}_1 + \zeta \vb{n}_{,\xi} = \vb{t}_1 - \zeta (W_1^{\phantom{1}1} \vb{t}_1 + W_1^{\phantom{1}2} \vb{t}_2) \, , &&
    \vb{g}_2 = \vb{t}_2 + \zeta \vb{n}_{,\eta} = \vb{t}_2 - \zeta (W_2^{\phantom{1}1} \vb{t}_1 + W_2^{\phantom{1}2} \vb{t}_2) \, , &&
    \vb{g}_3 = \vb{n} \, ,
\end{align}
using \cref{eq:dern}. Clearly, $\vb{g}_1$ and $\vb{g}_2$ are tangential to the middle surface $\vb{g}_\alpha \perp \vb{n}$, such that we immediately obtain $\vb{g}^3 = \vb{g}_3 = \vb{n}$ and $\vb{g}^\alpha \perp \vb{n}$ for the contravariant basis. Now, let $\vb{t}^1$ and $\vb{t}^2$ be the contravariant basis of the parametrisation of the middle surface $\vb{r} = \vb{r}(\xi,\eta)$, the vectors satisfy $\con{\vb{t}^\alpha}{\vb{t}_\beta} = \delta_{\alpha\beta}$ and $\vb{t}^\alpha \perp \vb{n}$, such that we can use them to the define the contravariant basis of $\vb{x}(\xi,\eta,\zeta)$ implicitly as
\begin{align}
    &\vb{g}^1 = c_{11} \vb{t}^1 + c_{12} \vb{t}^2 \, , && \vb{g}^2 = c_{21} \vb{t}^1 + c_{22} \vb{t}^2 \, .
\end{align}
The contravariant basis must satisfy $\con{\vb{g}^\alpha}{\vb{g}_\beta} = \delta_{\alpha \beta}$, leading to following system of equations
\begin{subequations}
    \begin{align}
    \con{\vb{g}^1}{\vb{g}_1} &= (1-\zeta W_1^{\phantom{1}1}) c_{11} - \zeta W_1^{\phantom{1}2} c_{12} = 1 \, , \\  
    \con{\vb{g}^1}{\vb{g}_2} &= -\zeta W_2^{\phantom{1}1} c_{11} + (1-\zeta W_2^{\phantom{1}2}) c_{12} = 0 \, , \\
    \con{\vb{g}^2}{\vb{g}_1} &= (1-\zeta W_1^{\phantom{1}1}) c_{21} - \zeta W_1^{\phantom{1}2} c_{22} = 0 \, , \\  
    \con{\vb{g}^2}{\vb{g}_2} &= -\zeta W_2^{\phantom{1}1} c_{21} + (1-\zeta W_2^{\phantom{1}2}) c_{22} = 1 \, .
\end{align}
\end{subequations}
Solving the system leads to 
\begin{align}
    &c_{11} = \dfrac{1-\zeta W_2^{\phantom{1}2}}{1 -2H \, \zeta + K \, \zeta^2} \, , &&
    c_{12} = \dfrac{\zeta W_2^{\phantom{1}1}}{1 -2H \, \zeta + K \, \zeta^2} \, , &&
    c_{21} = \dfrac{\zeta W_1^{\phantom{1}2}}{1 -2H \, \zeta + K \, \zeta^2} \, , &&
    c_{22} = \dfrac{1-\zeta W_1^{\phantom{1}1}}{1 -2H \, \zeta + K \, \zeta^2} \, .
\end{align}
Thus, the contravariant basis reads
\begin{align}
    &\vb{g}^1 = \dfrac{1}{1 -2H \, \zeta + K \, \zeta^2}[(1-\zeta W_2^{\phantom{1}2})\vb{t}^1 + \zeta W_2^{\phantom{1}1} \vb{t}^2] \, ,&&
    \vb{g}^2 = \dfrac{1}{1 -2H \, \zeta + K \, \zeta^2}[\zeta W_1^{\phantom{1}2} \vb{t}^1 + (1-\zeta W_1^{\phantom{1}1})\vb{t}^2] \, .
\end{align}
With the contravariant basis at hand we can write gradients of fields on $\vb{x}(\xi,\eta,\zeta)$ explicitly as
\begin{align}
    &\nabla \lambda = \vb{g}^i \lambda_{,i} = \lambda_{,\alpha} \vb{g}^\alpha + \lambda_{,\zeta}\vb{n} \, , && \lambda = (\lambda \circ \vb{x})(\bm{\xi}) \, .
\end{align}
As the tangential gradient of a field $\nabla_t \lambda = \lambda_{,\alpha} \vb{t}^\alpha$ is defined with respect to the parametrisation of the middle surface $\vb{r}(\xi,\eta)$ we find
\begin{align}
    \lambda_{,\alpha} \vb{g}^\alpha = \dfrac{1}{1 -2H \, \zeta + K \, \zeta^2}[\lambda_{,\xi} ([1-\zeta W_2^{\phantom{1}2}]\vb{t}^1 + \zeta W_2^{\phantom{1}1} \vb{t}^2) + \lambda_{,\eta}(\zeta W_1^{\phantom{1}2} \vb{t}^1 + [1-\zeta W_1^{\phantom{1}1}]\vb{t}^2)] \, ,
\end{align}
being the product of 
\begin{align}
    \dfrac{1}{1 -2H \, \zeta + K \, \zeta^2} [(1-\zeta W_2^{\phantom{1}2})\vb{t}^1 \otimes \vb{t}_1 + \zeta W_1^{\phantom{1}2} \vb{t}^1 \otimes \vb{t}_2 + \zeta W_2^{\phantom{1}1} \vb{t}^2 \otimes \vb{t}_1 + (1-\zeta W_1^{\phantom{1}1})\vb{t}^2 \otimes \vb{t}_2] (\lambda_{,\xi}\vb{t}^1 + \lambda_{,\eta}\vb{t}^2) \, .
\end{align}
Now, applying the skew-symmetric transformation tensor
\begin{align}
    \bm{T} = \vb{t}_1 \otimes \vb{t}^2 - \vb{t}_2 \otimes \vb{t}^1  \, ,
\end{align}
to $\bm{W} = W_{\alpha}^{\phantom{\alpha}\beta} \vb{t}_\beta \otimes \vb{t}^\alpha$, we finally identify 
\begin{align}
    \lambda_{,\alpha}\vb{g}^\alpha = \dfrac{1}{1 -2H \, \zeta + K \, \zeta^2} [\Pt - \zeta (\bm{T} \bm{W}\bm{T})^T ] \nabla_t \lambda = \dfrac{1}{1 -2H \, \zeta + K \, \zeta^2} (\Pt + \zeta \bm{T} \bm{W}\bm{T}^T ) \nabla_t \lambda \, ,
    \label{eq:fulltander}
\end{align}
such that a gradient with respect to $\vb{x}(\xi,\eta,\zeta)$ can be written as
\begin{align}
    \nabla \lambda = \dfrac{1}{1 -2H \, \zeta + K \, \zeta^2}(\Pt + \zeta \bm{T} \bm{W} \bm{T}^T)\nabla_t \lambda + \lambda_{,\zeta}\vb{n} \, . 
    \label{eq:shelltangrad}
\end{align}
At this point it is important to note that \textbf{if the explicit structure of the function $\lambda$ is known}, it may be possible to \textbf{significantly simplify its gradient}. For example, given the following function
\begin{align}
    &\lambda(\vb{x}) = \lambda_0(\vb{r}) + \lambda_1(\zeta) \lambda_2(\vb{r}) \, , && \vb{r}:\omega \subset \R^2 \to \surf \subset \R^3   \,,
\end{align}
using the product rule we find
\begin{align}
    \nabla \lambda = \nabla \lambda_0 + \lambda_2 \nabla \lambda_1 + \lambda_1 \nabla \lambda_2 = \nabla_t \lambda_0 + \lambda_2 \lambda_{1,\zeta} \vb{n} + \lambda_1 \nabla_t \lambda_2 \, , 
    \label{eq:splittangrad}
\end{align}
since $\lambda_0$ and $\lambda_2$ are functions of the Cartesian coordinates of \textbf{solely the middle surface}, and $\nabla_t \lambda_1 = 0$. Thus, for certain explicit structures it is possible to circumvent the need for \cref{eq:fulltander}.  

\section{Curved beams in three dimensions}\label{ap:curves}
A curve in three-dimensional space can be mapped from some parametric space $\gamma \subset \R$ via
\begin{align}
    \vb{r} = \vb{r}(\xi) \, , && \vb{r}:\gamma \subset \R \to \R^3 \, .
\end{align}
A unit tangent vector to the curve is given by 
\begin{align}
    &\vb{t} = \dfrac{\vb{r}_{,\xi}}{\norm{\vb{r}_{,\xi}}}  \, , && \vb{r}_{,\xi} = \norm{\vb{r}_{,\xi}} \vb{t} \, .
\end{align}
Accordingly, the infinitesimal curve element is given by
\begin{align}
    \dd \curv = \norm{\dd \vb{r}} = \norm{\vb{r}_{,\xi}} \, \dd \xi  \, , && \der{s}{\xi} = \norm{\vb{r}_{,\xi}} \, , && \der{\xi}{s} = \norm{\vb{r}_{,\xi}}^{-1} \, ,
\end{align}
such that 
\begin{align}
    \vb{t} = \der{}{s} \vb{r} = \pder{\vb{r}}{\xi}\pder{\xi}{s} = \norm{\vb{r}_{,\xi}} \norm{\vb{r}_{,\xi}}^{-1} \vb{t}\,   \, .
\end{align}
We use the unit tangent vector to define the tangential and normal projection operators
\begin{align}
    &\Pt = \vb{t} \otimes \vb{t} \, , && \Pn = \one - \Pt = \one - \vb{t} \otimes \vb{t} \, .
\end{align}
Further, we re-parameterise the curve $\vb{r} = \vb{r}(s)$ using the arc-length parameter
\begin{align}
    s(\xi) = \int_0^\xi \dd \curv = \int_0^\xi \norm{\vb{r}_{,\xi}} \, \dd \xi \, , && \xi(s) = \int_0^s \dd \xi = \int_0^s \norm{\vb{r}_{,\xi}}^{-1} \, \dd \curv \, ,  && s \in [0,l] \, ,
\end{align}
where we assume that $\xi$ and $s$ start at zero for simplicity. 
Thus, the tangential gradient of a function $\lambda = \lambda(s)$ with respect to the curve reads
\begin{align}
    \nabla \lambda = \lambda_{,s} \vb{g}^s = \lambda_{,s} \vb{t} = \Pt \nabla \lambda = \nabla_t \lambda \, , 
\end{align}
where $\vb{g}^s$ is the contravariant pseudo-inverse of the covariant tangent vector $\vb{g}_s = \vb{t}$. For the divergence of some vector $\vb{v} = \vb{v}(s)$ we find
\begin{align}
    \di_t \vb{v} = \tr(\D_t \vb{v}) = \con{\D \vb{v}}{\Pt} = \con{\vb{v}_{,s} \otimes \vb{t}}{\vb{t}\otimes \vb{t}} = \con{\vb{v}_{,s}}{\vb{t}} \, .  
\end{align}
Next, by defining two orthogonal unit vectors 
\begin{align}
   \vb{n} = \vb{n}(s) \, , &&  \vb{c}(s) = \vb{t} \times \vb{n} \, , && \con{\vb{t}}{\vb{n}} = \con{\vb{t}}{\vb{c}} = \con{\vb{n}}{\vb{c}} = 0 \, , && \norm{\vb{n}} = \norm{\vb{c}} = 1 \, ,
\end{align}
we can construct a moving coordinate system along the curve. 
We observe that 
\begin{align}
    \der{}{s}\con{\vb{t}}{\vb{t}} =  2\con{\vb{t}_{,s}}{\vb{t}} = 0 \, ,  
\end{align}
implying $\vb{t}_{,s} \perp \vb{t}$. Thus, we can define
\begin{align}
    \vb{t}_{,s} = \kappa_n \vb{n} + \kappa_c \vb{c} \, , && \con{\vb{t}_{,s}}{\vb{n}} = \kappa_n \, , && \con{\vb{t}_{,s}}{\vb{c}} = \kappa_c \, , 
\end{align}
such that $\kappa_n$ and $\kappa_c$ are curvature measures.
There holds
\begin{align}
    &\der{}{s} \con{\vb{n}}{\vb{n}} = 2\con{\vb{n}_{,s}}{\vb{n}} = 0  \, , && \der{}{s} \con{\vb{t}}{\vb{n}} = \con{\vb{t}_{,s}}{\vb{n}} + \con{\vb{t}}{\vb{n}_{,s}} = 0 \, , 
\end{align}
such that we can define
\begin{align}
    &\vb{n}_{,s} = - \kappa_n \vb{t} + \tau \vb{c} \, , && \con{\vb{n}_{,s}}{\vb{t}} = -\con{\vb{t}_{,s}}{\vb{n}} = - \kappa_n \, , && \con{\vb{n}_{,s}}{\vb{c}} = \tau \, ,
\end{align}
where $\tau$ is called the torsion of the curve. Finally, we find
\begin{align}
    \vb{c}_{,s} = \der{}{s}(\vb{t} \times \vb{n}) = \vb{t}_{,s} \times \vb{n} + \vb{t} \times \vb{n}_{,s} = -\kappa_c \vb{t} - \tau \vb{n} \, .
\end{align}

Next, we define a beam in three-dimensional space as a curve with a thickness given by some surface mapping
\begin{align}
    &\vb{x}(s,\eta,\zeta) = \vb{r} + \eta \vb{n} + \zeta \vb{c} \, , && \vb{x}: [0,l]\times\omega \subset \R^3 \to \vol \subset \R^3 \, ,
\end{align}
where $\eta$ and $\zeta$ span the surface $\omega \subset \R^2$ of the cross-section of the beam and $l$ is its length. The covariant basis of the map reads
\begin{align}
    \vb{g}_1 = \vb{x}_{,s} = (1- \kappa_n\eta - \kappa_c\zeta)\vb{t} - \tau \zeta \vb{n} + \tau \eta \vb{c} \, , && \vb{g}_2 = \vb{x}_{,\eta} = \vb{n} \, , && \vb{g}_3 = \vb{x}_{,\zeta} = \vb{c} \, .
\end{align}
Thus, we immediately get that the contravariant basis satisfies $\vb{g}^1 \perp \vb{n}$, $\vb{g}^1 \perp \vb{c}$, $\vb{g}^2 \perp \vb{c}$ and $\vb{g}^3 \perp \vb{n}$ by the Kronecker delta property $\con{\vb{g}^i}{\vb{g}_j} = \delta_{ij}$, such that
\begin{align}
    &\vb{g}^1 = \dfrac{1}{1- \kappa_n\eta - \kappa_c\zeta} \vb{t} \, , &&
    \vb{g}^2 = \dfrac{\tau \zeta}{1- \kappa_n\eta - \kappa_c\zeta} \vb{t} + \vb{n} \, ,  && 
    \vb{g}^3 = \dfrac{-\tau \eta}{1- \kappa_n\eta - \kappa_c\zeta} \vb{t} + \vb{c} \, .
\end{align}
Consequently, gradients of some function $\lambda = \lambda(\vb{x})$ with respect to the mapping read
\begin{align}
    \nabla \lambda &= \lambda_{,i}\vb{g}^i = \dfrac{1}{1- \kappa_n\eta - \kappa_c\zeta} (\lambda_{,s} + \tau \zeta \lambda_{,\eta} - \tau \eta \lambda_{,\zeta}) \vb{t} + \lambda_{,\eta} \vb{n} + \lambda_{,\zeta} \vb{c} \notag \\
    & = \dfrac{1}{1- \kappa_n\eta - \kappa_c\zeta} \nabla_t \lambda + \dfrac{1}{1- \kappa_n\eta - \kappa_c\zeta}(\cof\bm{T} + [1- \kappa_n\eta - \kappa_c\zeta]\Pn) (\lambda_{,\eta} \vb{n} + \lambda_{,\zeta} \vb{c}) \, ,
\end{align}
where $\bm{T}$ is defined as the torsion tensor of the beam
\begin{align}
    &\bm{T} = \tau(\eta\vb{c} -  \zeta \vb{n})  \otimes \vb{t} \, , && \bm{T}\vb{t} = \tau(\eta\vb{c} -  \zeta \vb{n}) \, , && \cof \bm{T} = \tau \vb{t} \otimes (\zeta \vb{n} - \eta\vb{c}) \, .   
\end{align}
If the structure of $\lambda$ is known, it is possible to exploit it to reduce the complexity of the gradient
\begin{align}
    &\lambda(\vb{x}) = \lambda_0(s) + \lambda_1(\eta,\zeta)  \, , && \nabla \lambda = \nabla_t \lambda_0 +  \dfrac{\tau \zeta \lambda_{1,\eta} - \tau \eta \lambda_{1,\zeta}}{1 - \kappa_n\eta - \kappa_c\zeta} \vb{t} + \lambda_{1,\eta}\vb{n} +  \lambda_{1,\zeta} \vb{c} \, ,
\end{align}
Clearly, if the cross-section of the beam does not twist $\tau = 0$, the gradient simplifies even further to
\begin{align}
    &\nabla \lambda = \nabla_t \lambda_0 + \lambda_{1,\eta}\vb{n} +  \lambda_{1,\zeta} \vb{c} = \nabla_t \lambda_0 + \Pn\nabla \lambda_1  = \nabla_t \lambda_1 + \nabla_n \lambda_1 \, , && \nabla_n \lambda = \Pn \nabla \lambda \, ,
\end{align}
motivating the definition of the normal gradient $\nabla_n(\cdot)$.
Finally, the mapping of an infinitesimal beam volume element reads
\begin{align}
    \dd \vol = \con{\vb{x}_{,\eta} \dd \eta \times \vb{x}_{,\zeta} \dd \zeta}{\vb{x}_{,s}\dd s } =  \con{\vb{t}}{(1- \kappa_n\eta - \kappa_c\zeta)\vb{t} - \tau \zeta \vb{n} + \tau \eta \vb{c}} \, \dd \eta \dd \zeta \, \dd s = (1- \kappa_n\eta - \kappa_c\zeta) \, \dd \omega \, \dd s \, .
    \label{eq:beamvoldecomp}
\end{align}
The term $1- \kappa_n\eta - \kappa_c\zeta$ represents the beam-analogue of the shell-shifter.

\section{Asymptotic analysis for shells}\label{ap:asym}
Let every term in an integral over $\zeta$ be a quadratic form
\begin{align}
    &(a + b \zeta)^2 = a^2 + 2 a  b\, \zeta + b^2 \zeta^2 \, , && a \neq a(\zeta) \, , && b \neq b(\zeta) \, , 
\end{align}
composed of a constant part and a linear part \textbf{with respect to} $\zeta$, we find the following integration formulae
\begin{subequations}
\label{eq:asymanalform}
    \begin{align}
        \int_{-h/2}^{h/2} (a + b \zeta)^2c \, \dd \zeta &= ch \, a^2 + \dfrac{ch^3}{12}b^2 \, , \\
        \int_{-h/2}^{h/2} (a + b \zeta)^2 (c\,\zeta) \, \dd \zeta &= \dfrac{ch^3}{6}ab   \, , \\
        \int_{-h/2}^{h/2} (a + b \zeta)^2 (c\,\zeta^2) \, \dd \zeta &= \dfrac{ch^3}{12}a^2  \, , 
    \end{align}
\end{subequations}
where we omitted higher order terms $\O(\xi^3)$, $c \neq c(\zeta)$, and terms of the form $(\cdot)\zeta$ vanish by the symmetry of integration $\int_{-h/2}^{h/2}(\cdot)\zeta\,\dd \zeta = 0$.

\end{document}